\documentclass[12pt]{amsart}
\usepackage{graphicx,amssymb,latexsym,eepic,booktabs}
\usepackage[dvips,centering,includehead,width=15.1cm,height=22.7cm]{geometry}
\usepackage[labelfont=sc,labelsep=period,margin=0pt,belowskip=0pt]{caption}
\usepackage[obeyspaces,spaces]{url}
\DeclareUrlCommand\arxiv{\urlstyle{tt}%
\def\UrlBreaks{\do\@\do\\\do\/\do\!\do\_\do\|\do\;\do\>\do\]%
 \do\)\do\,\do\?\do\'\do\+\do\=\do\#}%
}

\numberwithin{equation}{section}

\newtheorem{theorem}{Theorem}[section]
\newtheorem{proposition}[theorem]{Proposition}
\newtheorem{conjecture}[theorem]{Conjecture}
\newtheorem{corollary}[theorem]{Corollary}
\newtheorem{lemma}[theorem]{Lemma}

\theoremstyle{definition}

\newtheorem{remark}[theorem]{Remark}
\newtheorem{example}[theorem]{Example}

\newtheorem{definition}[theorem]{Definition}
\newtheorem{problem}[theorem]{Problem}

\newcommand{\margin}[1]{\marginpar{\tiny #1}}

\newcommand{\Arcs}{\mathbf{A}}
\newcommand{\Exch}{\mathbf{E}}
\newcommand{\Mark}{\mathbf{M}}

\newcommand{\Surf}{\mathbf{S}}

\newcommand{\SM}{(\Surf,\Mark)}
\newcommand{\SMprime}{(\Surf',\Mark')}

\newcommand{\DDSM}{\Delta^{\hspace{-1pt}\circ}(\Surf,\Mark)}
\newcommand{\DTSM}{\Delta^{\hspace{-1pt}\notch}(\Surf,\Mark)}
\newcommand{\EESM}{\Exch^{\circ\hspace{-1pt}}(\Surf,\Mark)}
\newcommand{\ETSM}{\Exch^{\notch}(\Surf,\Mark)}
\newcommand{\AASM}{\Arcs^{\hspace{-1pt}\circ}(\Surf,\Mark)}
\newcommand{\ATSM}{\Arcs^{\hspace{-1pt}\notch}(\Surf,\Mark)}

\newcommand{\Acal}{\mathcal{A}}
\newcommand{\Bcal}{\mathcal{B}}
\newcommand{\Fcal}{\mathcal{F}}
\newcommand{\Rcal}{\mathcal{R}}
\newcommand{\Scal}{\mathcal{S}}
\newcommand{\Xcal}{\mathcal{X}}
\newcommand{\MCG}{\mathcal{MCG}}

\newcommand{\CC}{\mathbb{C}}
\newcommand{\PP}{\mathbb{P}}
\newcommand{\ZZ}{\mathbb{Z}}
\newcommand{\RR}{\mathbb{R}}

\newcommand{\SL}{\mathit{SL}}

\newcommand{\notch}{\scriptstyle\bowtie} 

\DeclareMathOperator{\rank}{rank}

\def\dd{\mathbf{d}} 
 
\def\pp{\mathbf{p}} 
\def\xx{\mathbf{x}}

\begin{document}

\title[Cluster algebras and triangulated surfaces]
{Cluster algebras and triangulated surfaces \\[.1in]
Part~I: Cluster complexes}

\author{Sergey Fomin}
\address{\hspace{-.3in} Department of Mathematics, University of Michigan,
Ann Arbor, MI 48109, USA} 
\email{fomin@umich.edu}

\author{Michael Shapiro}
\address{\hspace{-.3in} Department of Mathematics, Michigan State University,
East Lansing, MI 48824,~USA} 
\email{mshapiro@math.msu.edu}

\author{Dylan Thurston}
\address{\hspace{-.3in} Department of Mathematics, Barnard College, Columbia University, 
New York, NY 10027, USA}
\email{dpt@math.columbia.edu}

\date{\today 
}

\thanks{This work was partially supported by NSF grants DMS-0245385 and
  DMS-0555880 (S.~F.),
DMS-0401178 and PHY-0555346 (M.~S.), DMS-0071550
(D.~T.), BSF grant 2002375 (M.~S.), and an A.~P.~Sloan Research
Fellowship (D.~T.).}

\subjclass[2000]{
Primary 
16S99, 
Secondary 
05E99, 
57N05, 
57M50. 
}

\keywords{Cluster algebra, cluster complex,  exchange graph,
triangulated surface, flip, arc complex, matrix mutation}


\begin{abstract}
We establish basic properties of cluster algebras associated with
oriented bordered surfaces with marked points. 
In particular, we show that the underlying cluster complex of such a cluster algebra
does not depend on the choice of coefficients, describe this
complex explicitly  in terms of ``tagged triangulations'' of the
surface, and determine its homotopy type and its growth rate.


\end{abstract}

\maketitle

\ \vspace{-.3in}

\tableofcontents


\section{Introduction}

Cluster algebras are a class of commutative rings endowed with an
additional combinatorial structure, which involves a set of distinguished
generators (cluster variables) grouped into overlapping subsets (clusters) of
the same cardinality. The original motivation for cluster
algebras came from representation theory, specifically the study of
dual canonical bases and total positivity phenomena in semisimple Lie
groups. (See \cite{cdm} for an introduction to cluster algebras in the
historical order of development.) In subsequent years, constructions
involving cluster algebras were discovered in various mathematical
contexts, including a topological/geometric one that is the main focus of this
paper. 

Shortly after the introduction of cluster algebras~\cite{ca1}, 
cluster-algebraic structures
in Teichm\"uller theory were defined and studied by
M.~Gekhtman, M.~Shapiro, and A.~Vainshtein~\cite{gsv2} and, in a
more general setting, by
V.~Fock and A.~Goncharov~\cite{fock-goncharov1, fock-goncharov2},
making use of foundational work by R.~Penner~\cite{penner-decorated,penner-bordered} 
and V.~Fock~\cite{fock-dual}. 
These cluster algebras were built using the following crucial observation: 
in hyperbolic geometry, an analogue of the classical Ptolemy
Theorem 
holds for 
the \emph{lambda-lengths} of geodesic arcs triangulating 
a bordered and/or punctured surface with geodesic sides and cusps at the
vertices (the arcs connect horocycles drawn around cusps and punctures);
these lambda-lengths are the \emph{Penner coordinates} on the corresponding
\emph{decorated Teichm\"uller space}~\cite{penner-decorated, penner-universal}.
The constructions based on this observation made it possible to
transport the notions, results, and insights from the theory of 
cluster algebras (cluster monomials, Laurent phenomenon, finite type
classification, canonical bases, etc.) into geometric settings. 

This  paper and its sequel~\cite{fst-hyper} exploit the aforementioned
connection in the opposite direction:
we use ideas and techniques of combinatorial topology and 
hyperbolic geometry for the benefit of cluster algebra theory. 
We undertake a systematic study of
the wide class of cluster algebras whose underlying combinatorial
structures are determined and governed by the topology of
bordered surfaces.
We demonstrate that the pivotal role in the
construction of a cluster algebra associated with a bordered surface
is played by its \emph{tagged arc complex}, a certain simplicial
complex generalizing the classical notion of an arc complex. 

\smallskip

Recall that a cluster algebra is built around a combinatorial
scaffolding formed by 
\emph{exchange matrices}, 
related to each other by \emph{matrix mutations}. 
The basic observation in \cite{fock-goncharov1,fock-goncharov2,gsv2} 
is that this kind of structure arises when one considers 
\emph{signed adjacency matrices} associated with triangulations 
of an oriented surface, 
with vertices at a fixed set of \emph{marked points}. 
More specifically, 
matrix mutations arise as transformations of signed adjacency matrices 
that correspond to \emph{flips} of triangulations.

Here is a sketch of the main construction studied in this paper, a
class of cluster algebras
associated with 2-dimensional surfaces, aimed at a reader familiar
with the basics of cluster algebra theory (see, e.g., \cite{cdm}).  
Let $\Surf$ be a connected oriented ($2$-dimensional) Riemann surface with
boundary. 
Fix a finite set $\Mark$ of \emph{marked points} on~$\Surf$ that includes 
at least one marked point on each boundary component, plus 
possibly some interior points. 
For technical reasons, we need to exclude the possibility that 
$\Surf$ is closed, with exactly $2$ marked points.  (This restriction
will be removed in the sequel~\cite{fst-hyper} to this paper.)
Consider a triangulation $T$ of $\Surf$ whose vertices are
located at marked points in~$\Mark$ and whose edges are simple
pairwise non-intersecting curves (\emph{arcs}), 
considered up to isotopy, connecting some of the marked points. 
The triangulation $T$ defines a skew-symmetric matrix~$B(T)$  
defined by looking at the (signed) adjacencies between the edges
of~$T$. 
This matrix serves as an \emph{exchange
  matrix} of a cluster algebra~$\Acal$---actually, a whole class of
such algebras, depending on the choice of \emph{coefficients}. 
This class of cluster algebras is in fact independent of the choice 
of a triangulation~$T$; it only depends on the underlying
surface~$\SM$ with marked points. 
This is because any two triangulations are connected by \emph{flips},
and a flip corresponds to a \emph{mutation} of the corresponding
matrices~$B(T)$. 

The construction is more complicated that it may seem at first
glance because not every arc in any triangulation can be flipped,
whereas the cluster algebra axioms do not allow exceptions:
every cluster variable in every cluster is exchangeable. 
Consequently, not every seed in $\Acal$ is naturally associated with a
triangulation, and not every cluster variable is naturally labeled by
an arc. 

\pagebreak[3]

In this paper, we resolve these problems,
and proceed further to establish the basic properties of cluster
algebras associated with oriented bordered surfaces with marked points. 
In particular, we obtain the following results: 
\begin{itemize}
\item
\emph{The seeds in $\Acal$ are uniquely determined by their clusters} 
[Theorem~\ref{th:cluster-complex-top}]. 
\item
\emph{The exchange graph and the cluster complex of $\Acal$ are uniquely
determined by the surface~$\SM$; i.e., they do not depend on the choice of
coefficients} 
[Theorem~\ref{th:cluster-complex-top}]. 
\item
\emph{The seeds containing a particular cluster variable form a connected
subgraph of the exchange graph} 
[Theorem~\ref{th:cluster-complex-top}].
\item
\emph{There are finitely many different exchange matrices appearing in the
seeds of~$\Acal$}
[Corollary~\ref{pr:mut-class-ideal}]. 
\item 
\emph{The class of exchange matrices of cluster algebras of topological
  origin (i.e., integer matrices that can be realized as signed adjacency matrices of
some triangulated bordered surface) can be given a concrete combinatorial 
description} 
[Theorem~\ref{th:B-criterion}].
\item
\emph{The exponents appearing in denominators of Laurent expansions
of cluster variables can be identified as certain intersection numbers}
[Theorem~\ref{th:denoms-as-inters-numbers}]. 
\item
\emph{The cluster complex of $\Acal$ has a concrete combinatorial
description in terms of ``tagged arcs''}
[Theorem~\ref{th:cluster-cpx=arc-cpx}]. 
\item
\emph{The cluster complex of $\Acal$ is a flag complex, i.e., it is the clique 
complex for its own $1$-skeleton} 
[Theorem~\ref{th:cluster-complex-top}].  
\item
\emph{The cluster complex of $\Acal$ is either contractible, or homotopy
  equivalent to a sphere}
[Theorem~\ref{th:cluster-complex-contractible-or-sphere}]. 
\item
\emph{The exchange graphs of such cluster algebras can be completely
  classified  according to their growth (polynomial vs.\ exponential)}
[Section~\ref{sec:poly-growth}]. 
\end{itemize}
The results above in particular establish, for the class of
cluster algebras under consideration, several conjectures which are
expected to hold for \emph{any} cluster algebra
(\cite[Conjecture~4.14, parts (1)--(3)]{cdm}, \cite[Conjecture~7.4]{ca4}, and
Conjecture~\ref{conj:flag-complex} of this paper).




\smallskip

Our main constructions depend only on the combinatorics of
triangulated surfaces, and could in principle be formulated in entirely
combinatorial language, with no reference to topology or geometry. 
In the sequel~\cite{fst-hyper} to this paper, we are going to 
take a complementary, geometric approach to the construction and study
of cluster algebras associated with bordered surfaces; 
this approach is based on an interpretation of cluster variables 
as certain generalized lambda-lengths. 


\smallskip

Our results can be extended, with modifications, to 
cluster algebras whose exchange matrices are not necessarily skew
symmetric, such as for example the finite and affine types $B$ and~$C$.
However, in this first paper we confine ourselves to the skew symmetric
case.

\smallskip

\pagebreak[3]

The paper is organized as follows.
Requisite background on triangulated surfaces 
and cluster algebras is presented, respectively, 
in Sections~\ref{sec:surfaces}--\ref{sec:adj-mat}
and Sections~\ref{sec:cluster-triangul}--\ref{sec:types-diagrams}.
In particular, we introduce, in Section~\ref{sec:cluster-triangul},
the class of cluster algebras $\Acal$ 
whose exchange matrices are associated
with a bordered surface with marked points $\SM$ as explained above.
The structural properties of any such algebra~$\Acal$ are stated
in Theorem~\ref{th:cluster-complex-top}. 
In Section~\ref{sec:tagged-arcs}, we introduce our main new
combinatorial construction of a \emph{tagged arc complex}~$\DTSM$. 
This simplicial complex is a pseudomanifold \linebreak[3]
(Theorem~\ref{th:tagged-pseudo}).
Our main result (Theorem~\ref{th:cluster-cpx=arc-cpx}) asserts that 
$\DTSM$ is in fact isomorphic to the
cluster complex~$\Delta\SM$ of~$\Acal$, provided $\SM$
is not a closed surface with at most two marked points.
Thus, the cluster variables in~$\Acal$ are naturally labeled by the
tagged arcs in~$\SM$, and they form a cluster if and only if the
corresponding tagged arcs form a maximal simplex in~$\DTSM$. 

In Section~\ref{sec:intersection-pairing}, we introduce the 
notion of a (generalized) intersection number for tagged arcs, 
and show that these numbers satisfy the \emph{tropical} version of the
exchange relation recurrences. 
This leads to an interpretation of these intersection numbers 
as exponents in the denominators of Laurent expansions of
cluster variables (see Theorem~\ref{th:denoms-as-inters-numbers}). 

The proof of Theorems~\ref{th:cluster-cpx=arc-cpx}
and~\ref{th:cluster-complex-top} 
presented in Section~\ref{sec:proofs} 
is based on a generalization to the tagged setup of a classical result
from combinatorial topology (see Theorem~\ref{th:flip-cycles})
asserting that the fundamental group of the graph of flips 
is generated by cycles of length $4$ and~$5$ (commuting flips plus the
pentagon relation). 
In the tagged case, this property fails for closed surfaces with $2$
marked points, which explains our need to exclude these surfaces in
the present paper.

In Section~\ref{sec:cluster-cplx-conj}, we determine, up to homotopy, 
the topology of cluster complexes (or tagged arc complexes) 
associated with bordered surfaces.
Specifically, we show that the tagged arc complex is homotopy equivalent to
a sphere for polygons with $\le1$ puncture and for closed surfaces
(with punctures); in all other cases, the cluster complex is
contractible. 

In Section~\ref{sec:poly-growth}, we resolve the dichotomy of
polynomial vs.\ exponential growth for the exchange
graphs of the cluster algebras in question. 

The mutation equivalence class formed by exchange matrices of a
cluster algebra associated with a triangulated surface is finite. 
This leads to the discussion, in Section~\ref{sec:finit-mut-classes},  of 
the general problem of recognizing (or classifying) 
finite mutation classes.
Some examples of this kind 
(e.g., finite and affine types $E_6$ through $E_8$, Grassmannians
$\operatorname{Gr}_{3,9}$ and $\operatorname{Gr}_{4,8}$) 
do \emph{not} come from triangulated
surfaces. 



The main result of Section~\ref{sec:block-decomp} is a purely
combinatorial description of the class of matrices $B$ that arise as
exchange matrices of cluster algebras associated with bordered
surfaces. 
In other words, we characterize those matrices $B$ for which there
exists a bordered surface with marked points $\SM$ such that $B=B(T)$
for some tagged (equivalently, ordinary---see
Proposition~\ref{pr:mut-class-ideal}) 
triangulation~$T$ of~$\SM$. 

The problem of recovering the topology of $\Surf$ from the matrix~$B(T)$
associated with a triangulation~$T$ of~$\SM$ is discussed in 
Section~\ref{sec:recover-topology}.  For closed surfaces with
punctures, all topological information can be recovered from~$B(T)$;
in general, it cannot.
In Theorem~\ref{th:corank}, we compute the rank of $B(T)$ for an
arbitrary (ideal) triangulation~$T$.
Specifically, the corank of $B(T)$ equals the number of punctures plus
the number of boundary components with an even number of marked
points. 

\pagebreak[3]

\section{Bordered surfaces with marked points and ideal triangulations}
\label{sec:surfaces}

This section reviews some standard 
background from combinatorial topology of surfaces
(c.f.\ e.g., \cite{chekhov-penner,fg-dual-teich,hatcher,ivanov}), 
setting  up terminology and notation to be used throughout the paper.


\begin{definition}[\emph{Bordered surfaces with marked points}] 
\label{def:ciliated}

Let $\Surf$ be a connected oriented $2$-dimen\-sion\-al Riemann surface with
boundary. Fix a finite set $\Mark$
of \emph{marked points} in the closure of~$\Surf$.
Marked points in the interior of~$\Surf$ are called \emph{punctures}. 
In the terminology of \cite{fg-dual-teich, fock-rosly}, $(\Surf,\Mark)$~is a
\emph{ciliated surface}\footnote{There is a slight difference between
  our setup and the one employed by V.~Fock and A.~Goncharov~\cite{fg-dual-teich}:
we use punctures where they use boundary components without marked points.}. 

We will work with triangulations of
$\Surf$ whose vertices are located at the marked points in~$\Mark$. 
We exclude situations in which such triangulations cannot be
constructed or there is only one triangulation.
Specifically, we assume that $\Mark$ is non-empty, that 
there is at least one marked point on each connected component of the
boundary of~$\Surf$, and that $\SM$ is none of the following: 
\begin{itemize}
\item
a sphere with one or two punctures; 
\item
an unpunctured or once-punctured monogon; 
\item
an unpunctured digon; or
\item
an unpunctured triangle. 
\end{itemize}
(In this paper, an \emph{$m$-gon} is a disk with $m$ marked points on
the boundary.) 
It will also be convenient to exclude, for technical reasons, 
\begin{itemize}
\item
a sphere with three punctures. 
\end{itemize}
A \emph{bordered surface with marked points} is an
object $\SM$ 
satisfying the aforementioned constraints. 
All of our subsequent constructions depend on a choice of such an
object. 
\end{definition}

Up to homeomorphism, $\SM$ is defined by the following data:
\begin{itemize}
\item
the genus $g$ of the original Riemann surface;
\item
the number $b$ of boundary components; 
\item
the integer partition with $b$ parts describing the number of marked
points on each boundary component; 
\item
the number $p$ of punctures.
\end{itemize}

\begin{definition}[\emph{Arcs}]
\label{def:arcs}
A (simple) \emph{arc} $\gamma$ in $\SM$ is a curve
in~$\Surf$ such that
\begin{itemize}
\item
the endpoints of $\gamma$ are marked points in~$\Mark$; 
\item
$\gamma$ does not intersect itself, except that its endpoints may
  coincide;
\item
except for the endpoints, $\gamma$ is disjoint from~$\Mark$ and from
the boundary of~$\Surf$;
\item
$\gamma$ 
does not cut out 
an unpunctured monogon or an unpunctured digon.  
(In other words, $\gamma$
is not contractible into $\Mark$ or onto the boundary of~$\Surf$.)
\end{itemize}
Each arc $\gamma$ is considered up to isotopy inside the class of such
curves. 
\end{definition} 

The set $\AASM$ of all arcs in $\SM$ is typically infinite; 
the cases where it is finite can be easily classified. 

\begin{proposition}[{e.g., \cite[Section~2]{fg-dual-teich}}] 
\label{pr:finite-number-of-arcs}
The set $\AASM$ of all arcs in~$\SM$ is finite if and only if $\SM$ is
an unpunctured or once-punctured polygon. 
\end{proposition}

An arc whose endpoints coincide is called a \emph{loop}. 

\begin{definition}[\emph{Compatibility of arcs}]
\label{def:compat-arcs}
Two arcs are called \emph{compatible} if they do not intersect in the
interior of~$\Surf$; 
more precisely, there are curves in their respective isotopy classes
which do not intersect in the interior of~$\Surf$. 
\end{definition} 

For instance, each arc is compatible with itself. 

\begin{proposition}\cite{fhs-geodesics}
Any collection of pairwise compatible arcs can be realized by 
curves in respective isotopy classes 
which do not intersect in the interior of~$\Surf$. 
\end{proposition}

\begin{definition}[{\emph{Ideal triangulations}; see, e.g.,
    \cite[Section~5.5]{ivanov}}] 
A maximal collection of distinct pairwise compatible arcs is called an 
\emph{ideal triangulation}. 
(Recall that we excluded all cases in which $\SM$ cannot be
triangulated.) 
The arcs of a triangulation cut the
surface~$\Surf$ into \emph{ideal triangles}. 
The three sides of an ideal triangle do not have to be distinct, 
i.e., we allow \emph{self-folded} triangles. 
We also allow for a possibility that two triangles share more than one
side. 
\end{definition} 

\begin{figure}[htbp] 
\begin{center} 
\setlength{\unitlength}{1.5pt} 
\begin{picture}(40,27)(-20,3) 
\thicklines 

\qbezier(0,0)(-30,30)(0,30)
\qbezier(0,0)(30,30)(0,30)
\put(0,0){\line(0,1){20}}
 
\put(3,13){\makebox(0,0){$i$}}
\multiput(0,0)(0,20){2}{\circle*{2}} 
\end{picture} 

\end{center} 
\caption{Self-folded ideal triangle} 
\label{fig:self-folded} 
\end{figure}
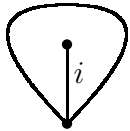

Figure~\ref{fig:self-folded} shows a triangle folded along the arc~$i$. 
Each side of an ideal triangle is either an arc 
or a segment of a boundary component between two marked points. 

\begin{figure}[htbp] 
\begin{center} 
\setlength{\unitlength}{3pt} 
\begin{picture}(20,19)(0,0) 
\thicklines 
\put(0,0){\line(1,0){20}}
\put(0,0){\line(100,173){10}}
\put(10,17.3){\line(100,-173){10}}
 
\put(0,0){\circle*{1}} 
\put(20,0){\circle*{1}} 
\put(10,17.3){\circle*{1}} 
\put(10,5.8){\circle*{1}} 

\put(0,0){\line(100,58){10}}
\put(20,0){\line(-100,58){10}}
\put(10,17.3){\line(0,-1){11.5}}

\end{picture} 
\qquad\qquad
\begin{picture}(20,19)(0,0) 
\thicklines 
\put(0,0){\line(1,0){20}}
\put(0,0){\line(100,173){10}}
\put(10,17.3){\line(100,-173){10}}
 
\put(0,0){\circle*{1}} 
\put(20,0){\circle*{1}} 
\put(10,17.3){\circle*{1}} 
\put(10,5.8){\circle*{1}} 

\put(0,0){\line(100,58){10}}
\put(20,0){\line(-100,58){10}}
\qbezier(0,0)(10,17.3)(20,0)

\end{picture} 
\qquad\qquad
\begin{picture}(20,19)(0,0) 
\thicklines 
\put(0,0){\line(1,0){20}}
\put(0,0){\line(100,173){10}}
\put(10,17.3){\line(100,-173){10}}
 
\put(0,0){\circle*{1}} 
\put(20,0){\circle*{1}} 
\put(10,17.3){\circle*{1}} 
\put(10,5.8){\circle*{1}} 

\put(0,0){\line(100,58){10}}
\qbezier(0,0)(10,17.3)(20,0)
\qbezier(0,0)(10,10)(12,7)
\qbezier(0,0)(13.5,3.8)(12,7)

\end{picture} 
\qquad\qquad
\begin{picture}(20,19)(0,0) 
\thicklines 
\put(0,0){\line(1,0){20}}
\put(0,0){\line(100,173){10}}
\put(10,17.3){\line(100,-173){10}}
 
\put(0,0){\circle*{1}} 
\put(20,0){\circle*{1}} 
\put(10,17.3){\circle*{1}} 
\put(10,5.8){\circle*{1}} 

\put(0,0){\line(100,58){10}}
\qbezier(0,0)(20,0)(10,17.3)
\qbezier(0,0)(10,10)(12,7)
\qbezier(0,0)(13.5,3.8)(12,7)

\end{picture} 

\end{center} 
\caption{Ideal triangulations of a once-punctured triangle} 
\label{fig:triangul-D3} 
\end{figure}
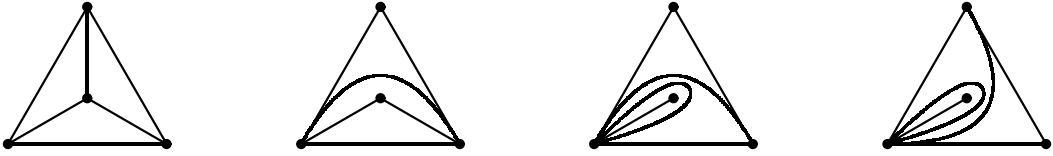

\begin{example}
Figure~\ref{fig:triangul-D3} shows four different kinds of ideal
triangulations of a once-punctured triangle. 
In this case, there are $10$ ideal triangulations altogether,
obtained from these four by $120^\circ$ and $240^\circ$ rotations. 
\end{example}

\begin{remark}
An alternative approach to working with triangulated surfaces 
is to start with a finite collection of disjoint oriented triangles,
take a partial matching of their sides (i.e., a set of disjoint pairs
of sides), then glue the triangles along the matched sides,
making sure that orientations match. 
\end{remark}

\begin{remark}
Yet another model is that of \emph{trivalent fat graphs} 
dual to (ideal) triangulations; see, e.g., 
  \cite{chekhov-penner,fock-goncharov1,fock-goncharov2}.
Note that by allowing self-folded triangles, 
we are allowing fat graphs with loops.
\end{remark}

The number $n$ of arcs in an ideal triangulation is an invariant of~$\SM$
called the \emph{rank} of~$\SM$. 
(The terminology comes from the fact that $n$ is the rank of any
cluster algebra in a class, to be introduced later, associated with
the surface~$\SM$.)
An easy count utilizing the Euler characteristic
(see, e.g., \cite[Section~2]{fg-dual-teich})
gives the following formula for the rank. 

\begin{proposition} 
\label{pr:n=6g+3b+3p+c-6}
Each ideal triangulation consists of 
\[
n=6g+3b+3p+c-6
\] 
arcs, where 
$g$ is the genus of~$\Surf$,
$b$ is the number of boundary components,
$p$ is the number of punctures, and 
$c$ is the number of marked points on the boundary. 
\end{proposition}

\begin{corollary}
For each positive integer~$n$, there are finitely many (homeomorphism classes of)
bordered surfaces $\SM$ with marked points whose rank is equal to~$n$. 
\end{corollary}

\begin{example}
By combining Proposition~\ref{pr:n=6g+3b+3p+c-6} with the restrictions 
in Definition~\ref{def:ciliated}, we can list all choices of~$\SM$
where the rank~$n$ is small (the nomenclature of ``types'' used below will be
explained later):
\begin{center}
\ \hspace{-.2in}
\begin{tabular}{ll}
$n=1$: & unpunctured square (type~$A_1$); \\
$n=2$: & unpunctured pentagon (type~$A_2$); \\
       & once-punctured digon (type~$A_1\times A_1$); \\
       & annulus with one marked point on each boundary component \\
       &\quad (type~$\widetilde A(1,1)$);\\
$n=3$: & unpunctured hexagon (type~$A_3$); \\
       & once-punctured triangle (type~$A_3$); \\
       & annulus with one marked point on one boundary component and two \\
       & \quad marked points on another (type~$\widetilde A(2,1)$); \\
       & once-punctured torus. 
\end{tabular}
\end{center}
\end{example}

Table~\ref{table:counts} shows the values of $g$, $b$, $p$, $m$,
and~$n$ for several examples of bordered surfaces with marked points. 

\begin{table}[htbp]
\begin{tabular}{lcccccc}\toprule
$\SM$ & $g$ & $b$ & $p$ & $c$ & $n$ & Type\\ \midrule
$(n+3)$-gon
        & $0$ & $1$ & $0$ & $n+3$ & $n$ & $A_n$\\
$n$-gon, $1$ puncture
        & $0$ & $1$ & $1$ & $n$ & $n$ & $D_n$\\
annulus,
$n_1+n_2$ marked points
        & $0$ & $2$ & $0$ & $n_1+n_2$ & $n_1+n_2$ & $\widetilde{A}(n_1,n_2)$ \\
$(n-3)$-gon, $2$ punctures
        & $0$ & $1$ & $2$ & $n-3$ & $n$ & $\widetilde{D}_{n-1}$\\
torus, $1$ puncture
        & $1$ & $0$ & $1$ & $0$ & $3$\\ \bottomrule
\end{tabular}
\caption{Examples of bordered surfaces with marked points}
\label{table:counts}
\end{table}

\begin{lemma}
\label{lem:honest-triangulation}
A bordered surface with marked points satisfying the
restrictions in Definition~\ref{def:ciliated}
has a triangulation without self-folded triangles.
\end{lemma}

\begin{proof}
Induction on~$n$. 
Let $\Surf'$ be the close oriented surface obtained by gluing a disk
to each boundary component of~$\Surf$. 
If the genus of $\Surf'$ is positive, then there is a loop $\gamma$ in
$\SM$ that goes around a handle in~$\Surf'$. 
Cut $\SM$ open along~$\gamma$, and proceed by induction.

It remains to treat the case of a sphere with holes and punctures.
In the absence of boundary, connect the punctures cyclically by
non-intersecting arcs; the polygons on both sides of the resulting
closed curve can be triangulated without self-folded triangles. 
If there is boundary, say with components labeled $1,\dots,b$, 
then draw non-intersecting arcs connecting components $1$ and~$2$, $2$
and~$3$, etc., and cut them open to reduce the claim to the case of
one boundary component (a disk with punctures).
It follows from the restrictions in Definition~\ref{def:ciliated} that
there will be at least $3$ marked points altogether, including at
least one marked point on the boundary of the disk.
Direct inspection shows that such a punctured disk can be triangulated
without self-folded triangles. 
\end{proof}

\section{Arc complexes}
\label{sec:arc-complexes}

\begin{definition}[\emph{Arc complex}; cf.\ \cite{harer-1985, harer-1986}]
\label{def:arc-complex} 
The \emph{arc complex} $\DDSM$ is the \emph{clique complex}
for the compatibility relation. 
That is, $\DDSM$ is the simplicial complex on the ground set 
$\AASM$ of all arcs in $\SM$ 
whose simplices are collections of distinct mutually compatible arcs,
and whose maximal simplices are the ideal triangulations. 
As they are all of the same cardinality 
(see Proposition~\ref{pr:n=6g+3b+3p+c-6}), 
the arc complex is \emph{pure}. 
\end{definition}

Our terminology differs from that of 
\cite{penner-probing, penner-arc-complexes} where the term ``arc complex''
is used for the quotient of $\DDSM$ by the \emph{pure mapping class
  group} of~$\SM$.

\begin{figure}[htbp]
\centerline{
        \includegraphics{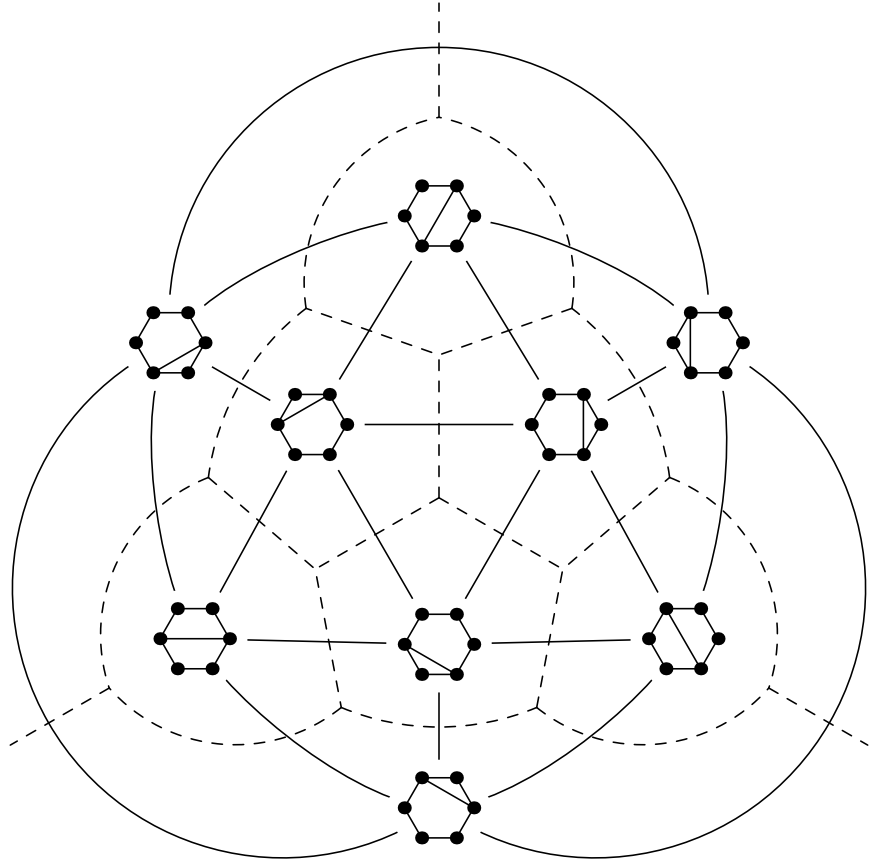}
}
\caption{The arc complex of a hexagon (solid lines)}
\label{fig:A3assoc_dual}
\end{figure}

\begin{example}
Figure~\ref{fig:A3assoc_dual} shows the arc complex of a hexagon,
i.e., an unpunctured disk with $6$~marked points on the boundary.
The exterior region is a triangle, which should be filled in.
It is then apparent that this arc
complex is topologically a $2$-dimensional sphere. 
Moreover, this complex can be realized as a boundary of a convex
polytope, the (polar) dual for a $3$-dimensional \emph{associahedron}. 
See, e.g., \cite[Section~3.1]{fomin-reading-pcmi} and references
therein. 
\end{example}

\begin{figure}
  \centering
  \includegraphics{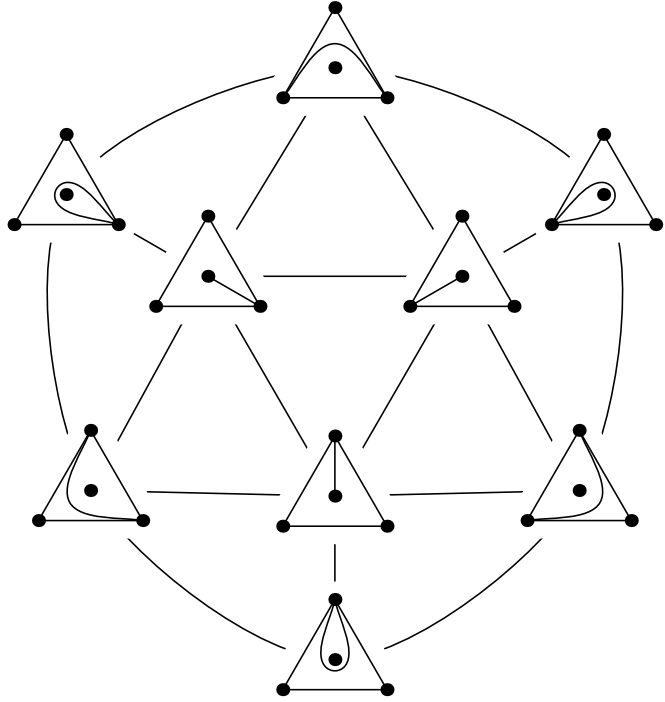}
  \caption{The arc complex of a once-punctured triangle}
\label{fig:arc-complex-D3}
\end{figure}

\begin{example}
  Figure~\ref{fig:arc-complex-D3} shows the arc complex of a
  once-punctured triangle, i.e., a disk with one puncture and 3 marked
  points on the boundary.  
\end{example}

\begin{example}
Figure~\ref{fig:Arc-complex-A(2,1)} shows the arc complex 
for an annulus with two marked points on one boundary component,
and one on another. Each maximal simplex is labeled by the
corresponding ideal triangulation.
\end{example}


\newsavebox{\ATwoOne}
\setlength{\unitlength}{1pt} 
\savebox{\ATwoOne}(40,60)[bl]{
\thicklines 
\put(20,30){\circle{20}} 
\put(20,30){\circle{60}} 
 
\put(20,0){\circle*{4}} 
\put(10,30){\circle*{4}} 
\put(30,30){\circle*{4}} 

\multiput(16,22)(2,0){5}{\circle*{0.5}} 
\multiput(14,24)(2,0){7}{\circle*{0.5}} 
\multiput(12,26)(2,0){9}{\circle*{0.5}} 
\multiput(12,28)(2,0){9}{\circle*{0.5}} 
\multiput(12,30)(2,0){9}{\circle*{0.5}} 
\multiput(12,32)(2,0){9}{\circle*{0.5}} 
\multiput(12,34)(2,0){9}{\circle*{0.5}} 
\multiput(14,36)(2,0){7}{\circle*{0.5}} 
\multiput(16,38)(2,0){5}{\circle*{0.5}} 
}

\newsavebox{\leftloopTwoOne}
\setlength{\unitlength}{1pt} 
\savebox{\leftloopTwoOne}(40,60)[bl]{
\thicklines 
\put(23,30){\circle{26}}
}

\newsavebox{\rightloopTwoOne}
\setlength{\unitlength}{1pt} 
\savebox{\rightloopTwoOne}(40,60)[bl]{
\thicklines 
\put(17,30){\circle{26}}
}

\newsavebox{\leftarcTwoOne}
\setlength{\unitlength}{1pt} 
\savebox{\leftarcTwoOne}(40,60)[bl]{
\thicklines 
\qbezier(20,0)(10,15)(10,30)
}

\newsavebox{\leftleftarcTwoOne}
\setlength{\unitlength}{1pt} 
\savebox{\leftleftarcTwoOne}(40,60)[bl]{
\thicklines 
\qbezier(15,45)(30,45)(30,30)
\qbezier(15,45)(0,45)(0,30)
\qbezier(20,0)(0,15)(0,30)
}

\newsavebox{\leftleftleftarcTwoOne}
\setlength{\unitlength}{1pt} 
\savebox{\leftleftleftarcTwoOne}(40,60)[bl]{
\thicklines 
\qbezier(25,15)(10,15)(10,30)
\qbezier(25,15)(40,15)(40,30)
\qbezier(15,50)(40,50)(40,30)
\qbezier(15,50)(-5,50)(-5,30)
\qbezier(20,0)(-5,15)(-5,30)
}

\newsavebox{\rightarcTwoOne}
\setlength{\unitlength}{1pt} 
\savebox{\rightarcTwoOne}(40,60)[bl]{
\thicklines 
\qbezier(20,0)(30,15)(30,30)
}

\newsavebox{\rightrightarcTwoOne}
\setlength{\unitlength}{1pt} 
\savebox{\rightrightarcTwoOne}(40,60)[bl]{
\thicklines 
\qbezier(25,45)(10,45)(10,30)
\qbezier(25,45)(40,45)(40,30)
\qbezier(20,0)(40,15)(40,30)
}

\newsavebox{\rightrightrightarcTwoOne}
\setlength{\unitlength}{1pt} 
\savebox{\rightrightrightarcTwoOne}(40,60)[bl]{
\thicklines 
\qbezier(15,15)(30,15)(30,30)
\qbezier(15,15)(0,15)(0,30)
\qbezier(25,50)(0,50)(0,30)
\qbezier(25,50)(45,50)(45,30)
\qbezier(20,0)(45,15)(45,30)
}


\begin{figure}[htbp] 

\begin{center} 
\setlength{\unitlength}{1.2pt} 
\begin{picture}(270,400)(0,-75) 

\thicklines


\put(30,90){\makebox(0,0){\usebox{\ATwoOne}}} 
\put(30,90){\makebox(0,0){\usebox{\leftloopTwoOne}}} 
\put(30,90){\makebox(0,0){\usebox{\leftarcTwoOne}}} 
\put(30,90){\makebox(0,0){\usebox{\rightrightarcTwoOne}}} 

\put(30,210){\makebox(0,0){\usebox{\ATwoOne}}} 
\put(30,210){\makebox(0,0){\usebox{\leftloopTwoOne}}} 
\put(30,210){\makebox(0,0){\usebox{\leftarcTwoOne}}} 
\put(30,210){\makebox(0,0){\usebox{\leftleftleftarcTwoOne}}} 

\put(110,90){\makebox(0,0){\usebox{\ATwoOne}}} 
\put(110,90){\makebox(0,0){\usebox{\rightarcTwoOne}}} 
\put(110,90){\makebox(0,0){\usebox{\leftarcTwoOne}}} 
\put(110,90){\makebox(0,0){\usebox{\rightrightarcTwoOne}}} 

\put(110,210){\makebox(0,0){\usebox{\ATwoOne}}} 
\put(110,210){\makebox(0,0){\usebox{\leftleftarcTwoOne}}} 
\put(110,210){\makebox(0,0){\usebox{\leftarcTwoOne}}} 
\put(110,210){\makebox(0,0){\usebox{\leftleftleftarcTwoOne}}} 

\put(150,30){\makebox(0,0){\usebox{\ATwoOne}}} 
\put(150,30){\makebox(0,0){\usebox{\rightarcTwoOne}}} 
\put(150,30){\makebox(0,0){\usebox{\rightrightrightarcTwoOne}}} 
\put(150,30){\makebox(0,0){\usebox{\rightrightarcTwoOne}}} 

\put(150,150){\makebox(0,0){\usebox{\ATwoOne}}} 
\put(150,150){\makebox(0,0){\usebox{\rightarcTwoOne}}} 
\put(150,150){\makebox(0,0){\usebox{\leftarcTwoOne}}} 
\put(150,150){\makebox(0,0){\usebox{\leftleftarcTwoOne}}} 

\put(230,30){\makebox(0,0){\usebox{\ATwoOne}}} 
\put(230,30){\makebox(0,0){\usebox{\rightarcTwoOne}}} 
\put(230,30){\makebox(0,0){\usebox{\rightrightrightarcTwoOne}}} 
\put(230,30){\makebox(0,0){\usebox{\rightloopTwoOne}}} 

\put(230,150){\makebox(0,0){\usebox{\ATwoOne}}} 
\put(230,150){\makebox(0,0){\usebox{\rightarcTwoOne}}} 
\put(230,150){\makebox(0,0){\usebox{\leftleftarcTwoOne}}} 
\put(230,150){\makebox(0,0){\usebox{\rightloopTwoOne}}}

\put(-40,150){\makebox(0,0){\usebox{\ATwoOne}}} 
\put(-40,150){\makebox(0,0){\usebox{\leftloopTwoOne}}} 

\put(70,30){\makebox(0,0){\usebox{\ATwoOne}}} 
\put(70,30){\makebox(0,0){\usebox{\rightrightarcTwoOne}}} 

\put(70,150){\makebox(0,0){\usebox{\ATwoOne}}} 
\put(70,150){\makebox(0,0){\usebox{\leftarcTwoOne}}} 

\put(70,270){\makebox(0,0){\usebox{\ATwoOne}}} 
\put(70,270){\makebox(0,0){\usebox{\leftleftleftarcTwoOne}}} 

\put(190,-30){\makebox(0,0){\usebox{\ATwoOne}}} 
\put(190,-30){\makebox(0,0){\usebox{\rightrightrightarcTwoOne}}} 

\put(190,90){\makebox(0,0){\usebox{\ATwoOne}}} 
\put(190,90){\makebox(0,0){\usebox{\rightarcTwoOne}}} 

\put(190,210){\makebox(0,0){\usebox{\ATwoOne}}} 
\put(190,210){\makebox(0,0){\usebox{\leftleftarcTwoOne}}} 

\put(300,90){\makebox(0,0){\usebox{\ATwoOne}}} 
\put(300,90){\makebox(0,0){\usebox{\rightloopTwoOne}}} 

\multiput(100,45)(0,120){2}{\line(2,1){60}}
\multiput(100,15)(0,120){3}{\line(2,-1){60}}
\multiput(70,64)(0,120){2}{\line(0,1){52}}
\multiput(190,4)(0,120){2}{\line(0,1){52}}

\put(70,304){\line(0,1){16}}
\put(190,244){\line(0,1){26}}

\put(190,-64){\line(0,-1){16}}
\put(70,-4){\line(0,-1){26}}

\put(100,285){\line(2,1){26}}
\put(160,-45){\line(-2,-1){26}}

\put(-6,150){\line(1,0){42}}
\put(224,90){\line(1,0){42}}

\qbezier(-30,182)(-10,240)(36,270)
\qbezier(-34,183)(-19,230)(2,260)
\qbezier(-37,184)(-30,210)(-16,240)
\qbezier(-39,184)(-35,200)(-28,220)

\qbezier(-30,118)(-10,60)(36,30)
\qbezier(-34,117)(-19,70)(2,40)
\qbezier(-37,116)(-30,90)(-16,60)
\qbezier(-39,116)(-35,100)(-28,80)

\qbezier(290,122)(270,180)(224,210)
\qbezier(294,123)(279,170)(258,200)
\qbezier(297,124)(290,150)(276,180)
\qbezier(299,124)(295,140)(288,160)

\qbezier(290,58)(270,0)(224,-30)
\qbezier(294,57)(279,10)(258,-20)
\qbezier(297,56)(290,30)(276,0)
\qbezier(299,56)(295,40)(288,20)

\end{picture} 

\end{center} 
\caption{Arc complex $\DDSM$ for an annulus of type~$\widetilde A(2,1)$} 
\label{fig:Arc-complex-A(2,1)} 
\end{figure}

\begin{definition}
\label{def:flip}
A \emph{flip} is a transformation of an ideal triangulation $T$ that
removes an arc $\gamma$ and replaces it with a (unique) different arc $\gamma'$ that, 
together with the remaining arcs, forms 
a new ideal triangulation~$T'$. 
\end{definition}

To illustrate, adjacent triangulations in Figure~\ref{fig:triangul-D3}
are related to each other by flips. 

There is at most one way to
flip an arc $\gamma$ in an ideal triangulation~$T$.
If $\gamma$ is a ``folded'' side of a self-folded ideal triangle 
(the interior edge in Figure~\ref{fig:self-folded}),
then $\gamma$ cannot be flipped. 
Otherwise, removing $\gamma$ creates a 
quadrilateral face on~$\SM$, and we let $\gamma'$ be the other
``diagonal'' of that quadrilateral. 
Another way to look at flips is the following. 

\begin{lemma}
\label{lem:arc-cplx-pseudo}
The arc complex is a \emph{pseudomanifold with boundary}, i.e., each
maximal simplex is of the same dimension and each simplex of
codimension~$1$ is contained in at most two maximal simplices.
\end{lemma}

The boundary of $\DDSM$ consists of collections of arcs which contain
a loop bounding a punctured monogon but do not contain the arc
enclosed.

\begin{theorem}[{\cite{harer-1985,hatcher}}]
\label{th:arc-cplx-is-contractible}
The arc complex is contractible
except when $\SM$ is a polygon, i.e., a disk with no punctures.
\end{theorem}

The \emph{dual graph} of a pseudomanifold has maximal simplices as
its vertices, with edges connecting maximal simplices sharing a
codimension-$1$ face. 

Let $\EESM$ denote the dual graph of~$\DDSM$. 
Thus, the vertices of $\EESM$ are (labeled by) the ideal
triangulations of~$\SM$, and the edges correspond to the flips.  
Examples are shown in Figure~\ref{fig:A3assoc_dual} (in dashed lines)
and in Figures~\ref{fig:A2assoc_basic}, \ref{fig:eesm-d3} and~\ref{fig:A(2,1)}.

\begin{figure}[phtb]
\centerline{
        \includegraphics{draws/arcs.30}
}
\caption{Graph $\EESM$ for an unpunctured pentagon}
\label{fig:A2assoc_basic}
\end{figure}


In Figure~\ref{fig:A(2,1)}, some readers may recognize the
``two-layer brick wall'' of \cite[Figure 4]{ca1}.  This coincidence
will be explained later; see
Example~\ref{example:cluster-type-A-affine}.

\begin{figure}
  \[
  \includegraphics{draws/arcs.10}
  \]
  \caption{Graph $\EESM$ for a once-punctured triangle}
  \label{fig:eesm-d3}
\end{figure}

\begin{figure}[phtb] 
\begin{center} 
\setlength{\unitlength}{1pt} 
\begin{picture}(270,400)(0,-75) 

\thicklines

\multiput(30,-70)(100,0){3}{\line(0,1){390}}
\multiput(30,0)(0,100){4}{\line(1,0){100}}
\multiput(130,-50)(0,100){4}{\line(1,0){100}}

\multiput(30,0)(0,100){4}{\circle*{5}} 
\multiput(130,-50)(0,50){8}{\circle*{5}} 
\multiput(230,-50)(0,100){4}{\circle*{5}} 

\put(-20,150){\makebox(0,0){\usebox{\ATwoOne}}} 
\put(-20,150){\makebox(0,0){\usebox{\leftloopTwoOne}}} 

\put(80,50){\makebox(0,0){\usebox{\ATwoOne}}} 
\put(80,50){\makebox(0,0){\usebox{\rightrightarcTwoOne}}} 

\put(80,150){\makebox(0,0){\usebox{\ATwoOne}}} 
\put(80,150){\makebox(0,0){\usebox{\leftarcTwoOne}}} 

\put(80,250){\makebox(0,0){\usebox{\ATwoOne}}} 
\put(80,250){\makebox(0,0){\usebox{\leftleftleftarcTwoOne}}} 

\put(180,0){\makebox(0,0){\usebox{\ATwoOne}}} 
\put(180,0){\makebox(0,0){\usebox{\rightrightrightarcTwoOne}}} 

\put(180,100){\makebox(0,0){\usebox{\ATwoOne}}} 
\put(180,100){\makebox(0,0){\usebox{\rightarcTwoOne}}} 

\put(180,200){\makebox(0,0){\usebox{\ATwoOne}}} 
\put(180,200){\makebox(0,0){\usebox{\leftleftarcTwoOne}}} 

\put(280,100){\makebox(0,0){\usebox{\ATwoOne}}} 
\put(280,100){\makebox(0,0){\usebox{\rightloopTwoOne}}} 

\end{picture} 

\end{center} 
\caption{Graph $\EESM$ for an annulus of type~$\widetilde A(2,1)$.  Each vertex
  corresponds to the triangulation made up of the three arcs labelling
  the adjacent regions.}
\label{fig:A(2,1)} 
\end{figure}

\begin{proposition}[{\cite{harer-1986,hatcher,mosher}}]
\label{pr:flips-connected}
The graph  $\EESM$ is connected. 
\end{proposition}

In other words, any two ideal triangulations are related by a sequence
of flips. Combining this statement with
Lemma~\ref{lem:honest-triangulation},
we obtain the following corollary, which can also be proved
directly. 

\begin{corollary}
\label{cor:honest-triangulation}
Any ideal triangulation can be transformed into a triangulation
without self-folded triangles by a sequence of flips. 
\end{corollary}

A more difficult result is 
Theorem~\ref{th:flip-cycles} below. 
It appeared explicitly in \cite[Theorem~1.1]{chekhov-penner};
see also \cite[Appendix~B]{teschner} and
\cite[Theorem~2.7]{chekhov-penner-quantizing}.
This result can also be extracted from 
\cite[Theorems 1.3 and~2.1]{harer-1986};  
(cf.\ the comments in \cite[pp.~190--191]{hatcher}. 

\pagebreak[3] 

\begin{theorem} 
\label{th:flip-cycles}
The fundamental group of  $\EESM$ 
is generated by cycles of length $4$ and~$5$, pinned down to a basepoint. 
More specifically, take the cycles of length~$4$ corresponding 
to pairs of commuting flips, and the cycles of length~$5$ obtained by 
alternately flipping  two arcs 
whose removal would create a pentagonal face 
(see Figure~\ref{fig:A2assoc_basic}). 
\end{theorem}

To clarify, by a ``pentagonal face'' in Theorem~\ref{th:flip-cycles}
we mean any face that looks like a pentagon to an observer located
inside it. One example of such a face is obtained by taking a
once-punctured triangle and connecting one vertex to the puncture. 

The graph $\EESM$ (or the complex~$\DDSM$) is finite, i.e., there are
finitely many ideal triangulations of~$\SM$, 
if and only if the total number of arcs in $\SM$ is finite, 
that is, in the cases listed in 
Proposition~\ref{pr:finite-number-of-arcs}. 

\section{Signed adjacency matrices and their mutations}
\label{sec:adj-mat}

In this section, we introduce and discuss the notion of a signed
adjacency matrix, a skew-symmetric matrix derived from an ideal
triangulation of a bordered surface. 
This construction is 
an extension of the ones proposed in the papers~\cite[Section~2]{fg-dual-teich}
and~\cite[Section~3]{gsv2}, which use different vocabulary and 
more restrictive setup. 

\begin{definition}
\label{def:signed-adj-matrix}
We associate to each ideal triangulation~$T$ the 
(generalized) \emph{signed adjacency matrix} $B=B(T)$ 
that reflects the combinatorics of~$T$. 
The rows and columns of $B(T)$ are naturally labeled by the arcs
in~$T$. For notational convenience, we arbitrarily 
label these arcs by the numbers $1,\dots,n$, so that the rows and
columns of~$B(T)$ are numbered from $1$ to~$n$ as customary,
with the understanding that this numbering of rows and columns is
temporary rather than intrinsic. 
For an arc (labeled)~$i$, let $\pi_T(i)$ denote (the label of) the arc
defined as follows: 
if there is a self-folded ideal triangle in~$T$ folded
along~$i$ (see Figure~\ref{fig:self-folded}), 
then $\pi_T(i)$ is its remaining side (the enclosing loop); 
if there is no such triangle, set $\pi_T(i)=i$. 

For each ideal triangle $\Delta$ in~$T$ which is not self-folded, 
define the $n \times n$
  integer matrix $B^\Delta=(b^\Delta_{ij})$ by setting 
\begin{equation}
\label{eq:B-proper}
b^\Delta_{ij}=
\begin{cases}
1 & \text{if $\Delta$ has sides labeled $\pi_T(i)$ and~$\pi_T(j)$,}\\
  & \text{\ \ with $\pi_T(j)$ following~$\pi_T(i)$ in the clockwise order;}\\
-1 & \text{if the same holds, with the counter-clockwise order;}\\
0 & \text{otherwise.}
\end{cases}
\end{equation} 
The matrix $B=B(T)=(b_{ij})$ is then defined by
\[
B=\sum_\Delta B^\Delta\,,
\]
the sum over all ideal triangles $\Delta$ in~$T$ which are not self-folded. 
The $n\times n$ matrix $B$ is skew-symmetric, and 
all its entries $b_{ij}$ are equal to $0$, $1$, $-1$, $2$, or~$-2$. 
\end{definition}

For triangulations without self-folded triangles, there is no need to
introduce the map~$\pi_T$, and the above definition of $B(T)$ 
coincides with the ones given in~\cite{fg-dual-teich,gsv2}. 

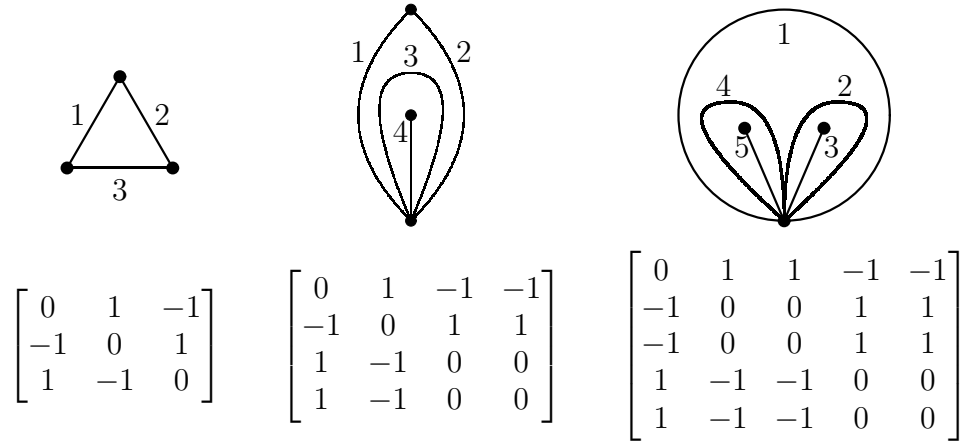
\begin{figure}[htbp] 
\begin{center} 
\begin{tabular}{ccc}
\setlength{\unitlength}{2pt} 
\begin{picture}(20,19)(0,-10) 
\thicklines 
\put(0,0){\line(1,0){20}}
\put(0,0){\line(100,173){10}}
\put(10,17.3){\line(100,-173){10}}
 
\put(0,0){\circle*{2}} 
\put(20,0){\circle*{2}} 
\put(10,17.3){\circle*{2}} 

\put(2,10){\makebox(0,0){$1$}} 
\put(18,10){\makebox(0,0){$2$}} 
\put(10,-4){\makebox(0,0){$3$}} 

\end{picture} 
&
\setlength{\unitlength}{8pt} 
\begin{picture}(10,10)(0,0) 

\qbezier(5,0)(0,5)(5,10)
\qbezier(5,0)(10,5)(5,10)
\put(5,0){\circle*{0.5}} 
\put(5,5){\circle*{0.5}} 
\put(5,10){\circle*{0.5}} 

\qbezier(5,0)(2,7)(5,7)
\qbezier(5,0)(8,7)(5,7)

\put(5,0){\line(0,1){5}}

\put(2.5,8){\makebox(0,0){$1$}} 
\put(7.5,8){\makebox(0,0){$2$}} 
\put(5,7.8){\makebox(0,0){$3$}} 
\put(4.5,4.2){\makebox(0,0){$4$}} 
\end{picture} 
&
\setlength{\unitlength}{1pt} 
\begin{picture}(80,80)(-20,0) 

\thicklines 
\put(20,40){\circle{80}} 
 
\put(20,0){\line(15,35){15}}
\put(20,0){\line(-15,35){15}}
 
\put(20,0){\circle*{4}} 
\put(5,35){\circle*{4}} 
\put(35,35){\circle*{4}} 

\qbezier(20,0)(-30,45)(0,45)
\qbezier(20,0)(70,45)(40,45)
\qbezier(20,0)(20,45)(0,45)
\qbezier(20,0)(20,45)(40,45)

\put(20,71){\makebox(0,0){$1$}}

\put(38,28){\makebox(0,0){$3$}} 
\put(4,28){\makebox(0,0){$5$}}

\put(43,51){\makebox(0,0){$2$}} 
\put(-3,51){\makebox(0,0){$4$}} 

\end{picture} 
\\[.1in]
$\begin{bmatrix}
0 & 1 & -1\\
-1& 0 & 1\\
1 & -1 &0
\end{bmatrix}$
&
\quad
$\begin{bmatrix}
0 & 1 & -1 &-1\\
-1& 0 & 1 & 1\\
1 & -1 &0 & 0\\
1 & -1 & 0& 0
\end{bmatrix}$
&
\quad
$\begin{bmatrix}
 0 & 1 & 1 &-1 &-1\\
-1 & 0 & 0 & 1 & 1\\
-1 & 0 & 0 & 1 & 1\\
 1 &-1 &-1 & 0 & 0\\
 1 &-1 &-1 & 0 & 0
\end{bmatrix}$
\end{tabular}
\end{center} 
\caption{The three puzzle pieces and signed adjacencies within each}
\label{fig:puzzle-pieces} 
\end{figure}

\begin{remark}
\label{rem:puzzle-pieces}
A more concrete version of Definition~\ref{def:signed-adj-matrix} 
that avoids the use of the map~$\pi_T$ is
based on representing an arbitrary ideal triangulation as a result of
gluing together a number of ``puzzle pieces'' which ``conceal''
self-folded triangles inside them. 
Specifically, any ideal triangulation, with a single exception
discussed below, can be obtained as follows:
\begin{itemize}
\item
Take several puzzle pieces of the form shown in
Figure~\ref{fig:puzzle-pieces} 
(a triangle, a triangulated once-punctured digon, or 
a triangulated twice-punctured monogon). 
\item
Take a partial matching of the exposed (outer) sides of these puzzle
pieces, never matching two sides of the same puzzle piece. In order to
obtain a connected surface, we need to ensure that
any two puzzle pieces in the collection are connected via matched
pairs. 
\item
Glue the puzzle pieces along the matched sides, making sure the
orientations match. 
\end{itemize}
The only exception to this construction is 
a triangulation of a $4$-punctured sphere obtained by gluing $3$
self-folded triangles to respective sides of an ordinary triangle.
(Recall that we do not allow a once-punctured monogon or a
thrice-punctured sphere, or else there would be more exceptions.) 

Each entry $b_{ij}$ of the matrix $B(T)$ is then obtained by 
adding up the $(0,\pm1)$-values describing signed adjacencies of 
the arcs $i$ and~$j$ inside each puzzle piece.
These are in turn determined by the matrices in
Figure~\ref{fig:puzzle-pieces}. 
\end{remark} 

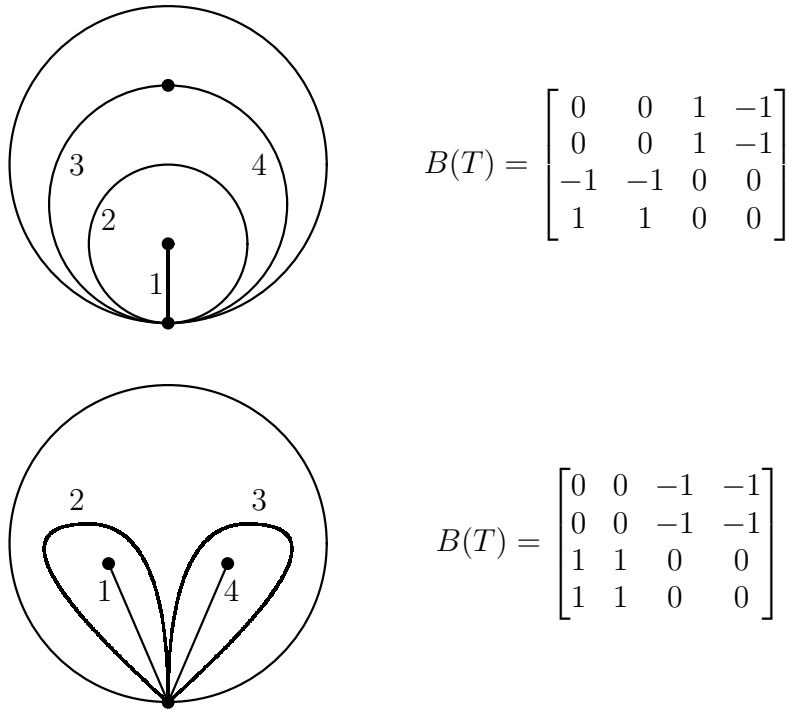
\begin{figure}[htbp] 
\begin{center} 
\setlength{\unitlength}{1.5pt} 
\begin{picture}(160,80)(0,0) 

\thicklines 
\put(20,20){\circle{40}} 
\put(20,30){\circle{60}} 
\put(20,40){\circle{80}} 
 
\put(20,0){\line(0,1){20}}
 
\multiput(20,0)(0,20){2}{\circle*{2.7}} 
\put(20,60){\circle*{2.7}} 

\put(17,10){\makebox(0,0){$1$}} 
\put(5,26){\makebox(0,0){$2$}} 
\put(-3,40){\makebox(0,0){$3$}} 
\put(43,40){\makebox(0,0){$4$}} 

\put(130,40){\makebox(0,0){
$\displaystyle
B(T)=\begin{bmatrix}
0 & 0 & 1 &-1\\
0 & 0 & 1 &-1\\
-1&-1 & 0 & 0\\
1 & 1 & 0 & 0
\end{bmatrix}
$}}
\end{picture} 

\setlength{\unitlength}{1.5pt} 
\begin{picture}(160,95)(0,0) 

\thicklines 
\put(20,40){\circle{80}} 
 
\put(20,0){\line(15,35){15}}
\put(20,0){\line(-15,35){15}}
 
\put(20,0){\circle*{2.7}} 
\put(5,35){\circle*{2.7}} 
\put(35,35){\circle*{2.7}} 

\qbezier(20,0)(-30,45)(0,45)
\qbezier(20,0)(70,45)(40,45)
\qbezier(20,0)(20,45)(0,45)
\qbezier(20,0)(20,45)(40,45)

\put(4,28){\makebox(0,0){$1$}} 
\put(36,28){\makebox(0,0){$4$}} 

\put(-3,51){\makebox(0,0){$2$}} 
\put(43,51){\makebox(0,0){$3$}} 

\put(130,40){\makebox(0,0){
$\displaystyle
B(T)=\begin{bmatrix}
0 & 0 & -1 & -1\\
0 & 0 & -1 &-1\\
1 & 1 & 0 & 0\\
1 & 1 & 0 & 0
\end{bmatrix}
$}}
\end{picture}

\end{center} 
\caption{Ideal triangulations of a twice-punctured monogon} 
\label{fig:twice-punctured-monogon} 
\end{figure}

\begin{example}[\emph{Twice-punctured monogon;
annulus of type~$\widetilde A(2,2)$}]
\label{example:monogon2}
Figure~\ref{fig:twice-punctured-monogon} shows two  
(labeled) ideal triangulations of a twice-punctured monogon,
and the corresponding signed adjacency matrices~$B(T)$. 

Figure~\ref{fig:A(2,2)} shows two 
ideal triangulations of an annulus with two marked points on each 
boundary component. The corresponding matrices $B(T)$ are exactly the
same as the ones in Figure~\ref{fig:twice-punctured-monogon}. 
\end{example}

\begin{figure}[htbp] 
\begin{center} 
\setlength{\unitlength}{2pt} 
\begin{picture}(60,63)(-10,0) 

\thicklines 
\put(20,30){\circle{20}} 
\put(20,30){\circle{60}} 
 
\put(20,0){\circle*{2}} 
\put(20,60){\circle*{2}} 
\put(10,30){\circle*{2}} 
\put(30,30){\circle*{2}} 

\multiput(16,22)(2,0){5}{\circle*{0.5}} 
\multiput(14,24)(2,0){7}{\circle*{0.5}} 
\multiput(12,26)(2,0){9}{\circle*{0.5}} 
\multiput(12,28)(2,0){9}{\circle*{0.5}} 
\multiput(12,30)(2,0){9}{\circle*{0.5}} 
\multiput(12,32)(2,0){9}{\circle*{0.5}} 
\multiput(12,34)(2,0){9}{\circle*{0.5}} 
\multiput(14,36)(2,0){7}{\circle*{0.5}} 
\multiput(16,38)(2,0){5}{\circle*{0.5}} 

\qbezier(20,0)(10,15)(10,30)
\qbezier(20,0)(30,15)(30,30)
\qbezier(20,60)(30,45)(30,30)

\qbezier(15,45)(30,45)(30,30)
\qbezier(15,45)(0,45)(0,30)
\qbezier(20,0)(0,15)(0,30)

\put(16,16){\makebox(0,0){$1$}} 
\put(30,11){\makebox(0,0){$4$}} 
\put(30,49){\makebox(0,0){$2$}} 
\put(10,49){\makebox(0,0){$3$}} 

\end{picture} 
\qquad\qquad
\begin{picture}(60,63)(-10,0) 

\thicklines 
\put(20,30){\circle{20}} 
\put(20,30){\circle{60}} 
 
\put(20,0){\circle*{2}} 
\put(20,60){\circle*{2}} 
\put(10,30){\circle*{2}} 
\put(30,30){\circle*{2}} 

\multiput(16,22)(2,0){5}{\circle*{0.5}} 
\multiput(14,24)(2,0){7}{\circle*{0.5}} 
\multiput(12,26)(2,0){9}{\circle*{0.5}} 
\multiput(12,28)(2,0){9}{\circle*{0.5}} 
\multiput(12,30)(2,0){9}{\circle*{0.5}} 
\multiput(12,32)(2,0){9}{\circle*{0.5}} 
\multiput(12,34)(2,0){9}{\circle*{0.5}} 
\multiput(14,36)(2,0){7}{\circle*{0.5}} 
\multiput(16,38)(2,0){5}{\circle*{0.5}} 

\qbezier(20,0)(10,15)(10,30)
\qbezier(20,0)(30,15)(30,30)
\qbezier(20,60)(30,45)(30,30)
\qbezier(20,60)(10,45)(10,30)

\put(30,11){\makebox(0,0){$4$}} 
\put(30,49){\makebox(0,0){$2$}} 
\put(10,11){\makebox(0,0){$1$}} 
\put(10,49){\makebox(0,0){$3$}} 

\end{picture} 

\end{center} 
\caption{Ideal triangulations of an annulus of type~$\widetilde A(2,2)$} 
\label{fig:A(2,2)} 
\end{figure}
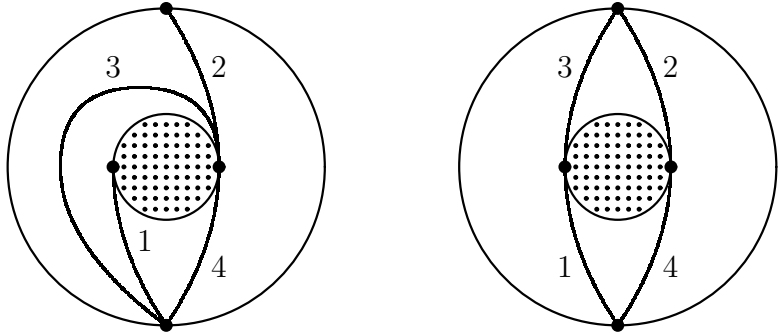 

\begin{example}[\emph{Annulus of type~$\widetilde A(1,1)$}]
Let $\SM$ be an annulus with one marked point on each boundary
component and no punctures. 
Then, for any (ideal) triangulation~$T$, the matrix~$B(T)$ has the form
\[
B(T)=\begin{bmatrix}
0 & 2\\
-2& 0
\end{bmatrix}\,,  
\] 
for an appropriate labeling of the arcs. 
See Figure~\ref{fig:A(1,1)}. 
\end{example}

\begin{figure}[htbp] 
\begin{center} 
\setlength{\unitlength}{2pt} 
\begin{picture}(60,60)(-10,0) 

\thicklines 
\put(20,30){\circle{20}} 
\put(20,30){\circle{60}} 
 
\put(20,0){\circle*{2}} 
\put(20,40){\circle*{2}} 

\multiput(16,22)(2,0){5}{\circle*{0.5}} 
\multiput(14,24)(2,0){7}{\circle*{0.5}} 
\multiput(12,26)(2,0){9}{\circle*{0.5}} 
\multiput(12,28)(2,0){9}{\circle*{0.5}} 
\multiput(12,30)(2,0){9}{\circle*{0.5}} 
\multiput(12,32)(2,0){9}{\circle*{0.5}} 
\multiput(12,34)(2,0){9}{\circle*{0.5}} 
\multiput(14,36)(2,0){7}{\circle*{0.5}} 
\multiput(16,38)(2,0){5}{\circle*{0.5}} 

\qbezier(20,0)(-15,40)(20,40)
\qbezier(20,0)(55,40)(20,40)

\put(37,40){\makebox(0,0){$1$}} 
\put(3,40){\makebox(0,0){$2$}} 

\end{picture} 

\end{center} 
\caption{Ideal triangulation of an annulus of type~$\widetilde A(1,1)$} 
\label{fig:A(1,1)} 
\end{figure}
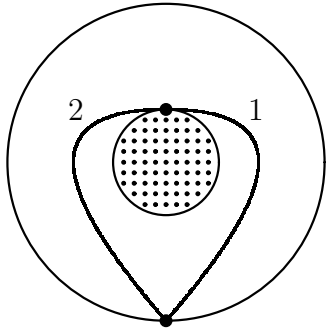 

\begin{example}[\emph{Once-punctured triangle; unpunctured hexagon}]
\label{example:triangle-hexagon}
Let $\SM$ be a disk with $3$~marked points on the boundary
and one puncture. 
Then each matrix $B=B(T)$ associated to an ideal triangulation $T$ of $\SM$ 
has one of the following four forms: 
\[
\begin{bmatrix}
0 & 1 & -1\\
-1& 0 & 1\\
 1&-1 & 0
\end{bmatrix}
,\quad
\begin{bmatrix}
0 & 1 & 0\\
-1& 0 & 1\\
 0&-1 & 0
\end{bmatrix}
,\quad
\begin{bmatrix}
0 & 1 & 0\\
-1& 0 & -1\\
 0& 1 & 0
\end{bmatrix}
,\quad
\begin{bmatrix}
0 & -1 & 0\\
1& 0 & 1\\
 0& -1 & 0
\end{bmatrix}
,
\] 
up to simultaneous relabeling of rows and columns. 
Cf.~Figure~\ref{fig:triangul-D3}. 

For an unpunctured hexagon,
we obtain exactly the same matrices. 
\end{example}

\begin{example}[\emph{Once-punctured torus}]
Let $\SM$ be a torus (no boundary) with one puncture. 
Then for any (ideal) triangulation~$T$, the matrix~$B(T)$ has the form
\[
B(T)=\begin{bmatrix}
0 & 2 & -2\\
-2& 0 & 2\\
 2&-2 & 0
\end{bmatrix}, 
\] 
for a suitable labeling of the arcs. 
\end{example}

In Section~\ref{sec:block-decomp}, we discuss the ``inverse problem''
of recovering information about the triangulated surface 
from the matrix~$B(T)$,
and the related problem of direct combinatorial description of the
class of matrices $B(T)$ that arise through this general construction. 


\begin{definition}[{\cite[Definition~4.2]{ca1}}]
\label{def:matrix-mutation} 
Let $B=(b_{ij})$ be an $n\times n$ integer matrix. 
We say that an $n\times n$ matrix $B'=(b'_{ij})$ is obtained from~$B$ by
\emph{matrix mutation} in direction~$k$,
and write $B' = \mu_k (B)$, 
if the entries of $B'$ are given by
\begin{equation}
\label{eq:matrix-mutation}
b'_{ij} =
\begin{cases}
-b_{ij} & \text{if $i=k$ or $j=k$;} \\[.05in]
b_{ij} + \displaystyle\frac{|b_{ik}| b_{kj} +
b_{ik} |b_{kj}|}{2} & \text{otherwise.}
\end{cases}
\end{equation}
\end{definition}

Matrix mutations are involutive, i.e., $\mu_k(\mu_k(B))=B$. 

If $B$ is a skew-symmetric matrix, then so is~$\mu_k(B)$. 

Proposition~\ref{pr:mut-tri} below shows that flips correspond to mutations
of the matrices~$B(T)$. 
This proposition can be verified by direct inspection. 
In restricted generality, it appeared
in~\cite{fock-goncharov1,fock-goncharov2,gsv2}. 

\begin{proposition}
\label{pr:mut-tri}
Suppose that an ideal triangulation $T'$ is obtained from $T$ by
a flip replacing an arc labeled~$k$. 
(The labeling of all other arcs remains unchanged.)
Then $B(T')= \mu_k(B(T))$.
\end{proposition}

\begin{proof}
Let us use the puzzle-piece decomposition of~$T'$ described in 
Remark~\ref{rem:puzzle-pieces}. 
Each flip occurs either inside a puzzle piece, or involves an arc 
shared by two puzzle pieces. 
In either case, the only entries of the
signed adjacency matrix affected by the flip 
are the ones located in the rows and columns 
associated with the arcs in the puzzle piece(s) involved. 
The proof therefore reduces to a routine case-by-case verification. 
\end{proof}

\begin{definition}\label{def:mutation-equivalent}
Two matrices are called \emph{mutation equivalent} if they can be
transformed into each other by a sequence of mutations. 
Let $\Bcal\SM$ denote the mutation equivalence class of a matrix
$B(T)$ associated with a triangulation~$T$ of~$\SM$.
\end{definition}

\begin{proposition}
  \label{pr:BSM-well-defined} The mutation equivalence
  class~$\Bcal\SM$ depends only on~$\SM$ but not on the choice
  of an ideal triangulation~$T$.
\end{proposition}

\begin{proof}
  This follows from Propositions~\ref{pr:flips-connected}
  and~\ref{pr:mut-tri}.  
\end{proof}

\section{Cluster algebras associated with triangulated surfaces}
\label{sec:cluster-triangul}

We begin this section by a swift review of cluster algebras, closely
following~\cite{ca2}.  This is not meant as an introduction to the
subject, but rather to fix notation and specify the generality we will
be working in.  For entry points to the theory, see,
for instance,~\cite{cdm,fomin-reading-pcmi}.

Let $(\PP,\oplus, \cdot)$ be a \emph{semifield}, i.e., 
$(\PP,\cdot)$ is an abelian group, 
$(\PP,\oplus)$ is a commutative semigroup, and the \emph{auxiliary
  addition}~$\oplus$ is distributive with respect to the multiplication. 

Let $\Fcal$ be a field isomorphic to the field of rational functions
in $n$ independent variables, 
with coefficients in~$\ZZ \PP$, the group ring of the multiplicative
group of~$\PP$.
 
\begin{definition}[\emph{Seeds and their mutations}] 
\label{def:seeds-and-mutations} 

A \emph{seed} in $\Fcal$ is a triple $\Sigma = (\xx, \pp, B)$, where 
 
\begin{itemize} 
\item 
the \emph{cluster} $\xx$ is a set of $n$ algebraically independent elements
  in~$\Fcal$ (called \emph{cluster variables}) 
that generate~$\Fcal$ over the field of fractions of $\ZZ \PP$. 
 
\item 
the \emph{coefficient tuple} 
$\pp = (p_x^\pm)_{x \in \xx}$ is a $2n$-tuple of elements of $\PP$ satisfying the 
\emph{normalization condition} 
$p_x^+ \oplus p_x^-  = 1$ for all $x \in \xx$. 
 
\item 
the \emph{exchange matrix} 
$B = (b_{xy})_{x,y \in \xx}$ is a skew-symmetrizable\footnote{A square
matrix is \emph{skew-symmetrizable} if it becomes skew-symmetric upon
rescaling of its columns by appropriate positive factors.}
  $n\times n$ integer matrix with rows and columns indexed by~$\xx$. 
\end{itemize} 

For such a seed $\Sigma = (\xx, \pp, B)$, and for any cluster variable
$z \in \xx$, the \emph{seed mutation} $\mu_z$ transforms $\Sigma$ into
a new seed 
$\overline \Sigma = (\overline \xx, \overline \pp, \overline B)$
defined as follows: 
\begin{itemize} 
 
\item 
$\overline \xx = \xx - \{z\} \cup \{\overline z\}$, where $\overline z \in \Fcal$ is determined 
by the \emph{exchange relation} 
\begin{equation} 
\label{eq:exchange-rel-xx} 
z\overline z
=p_z^+ \, \prod_{\substack{x\in\xx \\ b_{xz}>0}} x^{b_{xz}} 
+p_z^- \, \prod_{\substack{x\in\xx \\ b_{xz}<0}} x^{-b_{xz}} 
\end{equation} 
 
\item 
the $2n$-tuple $\overline \pp=(\overline p ^\pm_x)_{x\in\overline 
  \xx}$ 
is uniquely determined by the  normalization conditions 
$\overline p_x^+ \oplus \overline p_x^-  = 1$ together with 
 \begin{equation} 
\label{eq:p-mutation} 
{\overline p ^+_x}/{\overline p ^-_x} = 
\begin{cases} 
{p^-_z}/{p^+_z} & \text{if $x=\overline z$};\\[.05in] 
(p_z^+)^{b_{zx}}{p ^+_x}/{p ^-_x}   & \text{if $b_{zx}\geq 0$};\\[.05in] 
(p_z^-)^{b_{zx}}{p ^+_x}/{p ^-_x}   & \text{if $b_{zx}\leq 0$}. 
\end{cases} 
\end{equation}

\item 
the matrix $\overline B$ is obtained from $B$ by applying the matrix 
mutation in direction~$z$ and then relabelling one row and one column 
by replacing $z$ by~$\overline z$. 
 
\end{itemize} 
\end{definition}

It is easy to check that $\mu_{\overline z}(\overline\Sigma)=\Sigma$.

\begin{definition}[\emph{Normalized skew-symmetrizable cluster algebra}] 

\label{def:cluster-algebra} 
Let $\Scal$ be a mutation equivalence class of seeds in~$\Fcal$. 
Let $\Xcal$ denote the set of all cluster variables appearing in the
seeds of~$\Scal$.  
Let $\Rcal$ be any subring with unit in $\ZZ\PP$ 
containing all 
elements $p^\pm_x\in \pp$, for all seeds $\Sigma=(\xx, \pp, B)\in\Scal$. 
%
The (normalized) \emph{cluster algebra} 
$\Acal=\Acal(\Scal)$ is the $\Rcal$-subalgebra of $\Fcal$ generated 
by~$\mathcal{X}$. 

The exchange matrices $B$ appearing in the seeds of~$\Scal$ form a
mutation equivalence class, denoted by $\Bcal(\Acal)$. 
\end{definition} 

\begin{conjecture}[{\cite[p.~70]{ca2},
      \cite[Conjecture~4.14(2)]{cdm}}]
\label{conj:seed-by-cluster}
Each seed in a cluster algebra 
(i.e., each seed in a given mutation equivalence class)
is uniquely determined by its cluster. 
Two such seeds are related by a mutation 
if and only if the corresponding clusters share all elements but one. 
\end{conjecture}

\begin{definition}[\emph{Exchange graph, cluster
      complex}] 
\label{def:exchange-graph}
Let $\Acal=\Acal(\Scal)$ be the cluster algebra defined by~$\Scal$, 
a mutation equivalence class of seeds. 
The \emph{exchange graph} of $\Acal(\Scal)$
is the $n$-regular graph whose vertices are labeled by the seeds 
in~$\Scal$, and whose edges correspond to mutations. 

Assuming that Conjecture~\ref{conj:seed-by-cluster} holds for
a cluster algebra~$\Acal$, its underlying combinatorics 
is governed by its \emph{cluster complex}. 
This is the (possibly infinite) simplicial
complex on the ground set $\Xcal$ 
whose maximal simplices are the clusters. 
The exchange graph is then the dual graph of the cluster
complex: 
the vertices of the exchange graph can be identified with clusters,
with edges corresponding to pairs of clusters whose intersection has
cardinality~$n-1$. 
\end{definition}

Two cluster variables are \emph{compatible} if they appear together in 
some cluster. 

\begin{conjecture} 
\label{conj:flag-complex}
The cluster complex is always a ``flag complex.''
That is, the cluster complex is the clique complex 
for the compatibility relation on the set of all cluster variables.
\end{conjecture}

In other words, for any collection of pairwise compatible cluster
variables, there is a cluster containing all of them. 

We note that Conjecture~\ref{conj:flag-complex} would follow from 
\cite[Conjecture~7.4]{ca4} and \cite[Conjecture~4.14(3)]{cdm}. 

Our main object of study is the class of cluster algebras whose 
exchange matrices can be described by triangulated surfaces. 
More precisely, a cluster algebra~$\Acal$ belongs to this class 
if and only if $\Bcal(\Acal)=\Bcal\SM$,
for some bordered surface $\SM$ with marked points.  (See
Definition~\ref{def:mutation-equivalent} and
Proposition~\ref{pr:BSM-well-defined}.)
As we will see in Section~\ref{sec:types-diagrams},
this class in particular includes all cluster algebras of finite and
affine types $A$ and~$D$. 

The following result in particular establishes \cite[Conjecture~4.14, parts
  (1)--(3)]{cdm} and Conjecture~\ref{conj:flag-complex}
for the class of cluster algebras under consideration. 

\begin{theorem}
\label{th:cluster-complex-top}
Let $\SM$ be a bordered surface with marked points.
Assume that $\SM$ is \underbar{not} a closed surface with exactly $2$~punctures. 
Let $\Acal$ be a cluster algebra with $\Bcal(\Acal)=\Bcal\SM$. 
Then:
\begin{itemize}
\item
Each seed in $\Acal$ is uniquely determined by its cluster. 
\item
The cluster complex of $\Acal$ 
is uniquely determined by~$\SM$
(i.e., is independent of the choice of coefficients~$p^\pm_z$). 
\item
The seeds containing a particular cluster variable form a connected
subgraph of the exchange graph. 
\item
The cluster complex of $\Acal$ is the clique complex for the
compatibility relation on the set of cluster variables. 
\end{itemize}
\end{theorem} 

The proof of Theorem~\ref{th:cluster-complex-top} is given in
Section~\ref{sec:proofs}.  The cluster complex of any cluster algebra
with $\Bcal(\Acal) = \Bcal\SM$ is denoted by $\Delta\SM$.
A concrete description of $\Delta\SM$ in terms of
combinatorial topology of~$\SM$ is provided in
Theorem~\ref{th:cluster-cpx=arc-cpx}. 

We expect Theorem~\ref{th:cluster-complex-top} to hold in full
generality, including the case of closed surfaces with two punctures.
The proof presented in this paper does not however work in this
special case. A different proof covering all cases should appear
in~\cite{fst-hyper}. 

\pagebreak[3] 

\section{Types and diagrams}
\label{sec:types-diagrams}

We begin this section by reviewing the encoding of exchange matrices by
oriented graphs (quivers), the finite type classification
theorem~\cite{ca2}, and the related nomenclature of cluster algebras
of finite and affine type. 
We then discuss four examples (finite and affine types $A$ and~$D$)
that naturally fit into this paper's framework, i.e., 
they can be associated with bordered surfaces with marked points. 
(These examples correspond to the first four line entries in
Table~\ref{table:counts}.) 
In each case, we briefly mention appearances of cluster algebras of
the corresponding type as (homogeneous) coordinate rings of various
classical algebraic varieties such as Grassmannians, double Bruhat cells,
etc.; cf.\ \cite{ca3, ca2, scott}. 

\begin{definition}[\emph{Type of a cluster algebra}] 
\label{def:cluster-type}
We say that two cluster algebras $\Acal$ and~$\Acal'$ are of the same \emph{type}
if the corresponding mutation equivalence classes $\Bcal(\Acal)$ and
$\Bcal(\Acal')$ coincide, up to simultaneous relabeling of rows and columns. 
\end{definition}

A cluster algebra is of \emph{finite type} if it has finitely many
seeds. 
This terminology is consistent with Definition~\ref{def:cluster-type}
due to the following result.

\begin{theorem}[{\cite{ca2}}]
\label{th:finite-type-dep-on-B}
Whether a cluster algebra $\Acal$ is of finite or infinite type 
depends only on the mutation equivalence class~$\Bcal(\Acal)$. 
\end{theorem}

All cluster algebras of finite type have been classified~\cite{ca2}; 
this classification is completely parallel to the
Cartan-Killing classification of, say, finite crystallographic root
systems. 
As shown in~\cite{ca2}, Conjecture~\ref{conj:seed-by-cluster} holds
for any cluster algebra of finite type; its cluster complex is the
\emph{generalized associahedron}~\cite{yga} of the corresponding type. 

In this paper, we restrict our attention to cluster algebras whose
exchange matrices are \emph{skew-symmetric}. 
To state a more concrete version of
Theorem~\ref{th:finite-type-dep-on-B} in this restricted generality, 
we will need a construction that associates a mutation
equivalence class to an oriented graph (or \emph{quiver}).  

\begin{definition}[\emph{Encoding matrices by graphs}] 
\label{def:matrices-via-graphs}
Let $\Gamma$ be a finite directed graph with no loops
and no cycles of length~$2$; multiple edges are allowed. 
Let $B(\Gamma)=(b_{ij})$ denote the skew-symmetric
matrix whose rows and columns are labeled
by the vertices of~$\Gamma$, and whose entry $b_{ij}$ 
is equal to the number of edges going from $i$ to~$j$ 
minus the number of edges going from $j$ to~$i$. 
\end{definition}

\begin{lemma}[{\cite{ca2}}]
\label{lem:forest}
Let $\mathbb{T}$ be an unoriented graph obtained from a finite forest
(a disjoint union of trees)  
by replacing some of its edges by multiple ones. 
If $\Gamma_1$ and $\Gamma_2$ are two orientations of~$\mathbb{T}$, 
then $B(\Gamma_1)$ and $B(\Gamma_2)$ are mutation equivalent. 
\end{lemma}

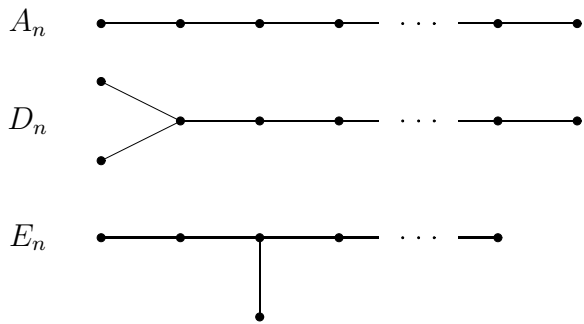
\begin{figure}[htbp] 
\vspace{-.2in} 
\[ 
\begin{array}{ccl} 
A_n && 
\setlength{\unitlength}{1.5pt} 
\begin{picture}(140,17)(0,-2) 
\put(0,0){\line(1,0){70}} 
\put(90,0){\line(1,0){30}} 
\multiput(0,0)(20,0){4}{\circle*{2}} 
\multiput(100,0)(20,0){2}{\circle*{2}} 
\multiput(76,0)(4,0){3}{\circle*{0.5}} 
\end{picture}\\[.1in] 
D_n 
&& 
\setlength{\unitlength}{1.5pt} 
\begin{picture}(120,17)(0,-2) 
\put(20,0){\line(1,0){50}} 
\put(90,0){\line(1,0){30}} 
\put(0,10){\line(2,-1){20}} 
\put(0,-10){\line(2,1){20}} 
\multiput(20,0)(20,0){3}{\circle*{2}} 
\multiput(100,0)(20,0){2}{\circle*{2}} 
\put(0,10){\circle*{2}} 
\put(0,-10){\circle*{2}} 
\multiput(76,0)(4,0){3}{\circle*{0.5}} 
\end{picture} 
\\[.2in] 
E_n 
&& 
\setlength{\unitlength}{1.5pt} 
\begin{picture}(140,17)(0,-2) 
\put(0,0){\line(1,0){70}} 
\put(90,0){\line(1,0){10}} 
\put(40,0){\line(0,-1){20}} 
\put(40,-20){\circle*{2}} 
\multiput(0,0)(20,0){4}{\circle*{2}} 
\multiput(100,0)(20,0){1}{\circle*{2}} 
\multiput(76,0)(4,0){3}{\circle*{0.5}} 
\end{picture} 
\\[.25in] 
\end{array} 
\] 
\vspace{-.1in} 
\caption{Simply-laced Coxeter-Dynkin diagrams} 
\label{fig:ADE} 
\end{figure}

It follows that there is a well-defined mutation equivalence class
$\Bcal(\mathbb{T})$ (associated with any such graph~$\mathbb{T}$ 
(i.e., with any finite forest with multiple edges). 
This in particular applies in a situation where~$\mathbb{T}$ is a 
simply-laced Coxeter-Dynkin diagram of type $A_n$ ($n\geq 1$), 
$D_n$ ($n\geq 4$), or~$E_n$ ($6\leq n\leq 8$),
or a disjoint union of such diagrams. 
See Figure~\ref{fig:ADE}.

\begin{theorem}[\cite{ca2}] 
\label{th:ADE-classif} 
A cluster algebra $\Acal$ of rank~$n$ with skew-symmetric exchange
matrices is of finite type if and only if the mutation equivalence
class $\Bcal(\Acal)$ is of the form $\Bcal(\mathbb{T})$,
where $\mathbb{T}$ is a disjoint union of simply-laced Coxeter-Dynkin
diagrams. 
\end{theorem}

In the situation described in Theorem~\ref{th:ADE-classif}, 
the type of~$\mathbb{T}$ (in the Cartan-Killing nomenclature) is uniquely
determined by the cluster algebra~$\Acal$, 
and is called \emph{the (cluster) type} of~$\Acal$
(or of the mutation equivalence class $\Bcal(\Acal)$,
or of any matrix in that class).
This terminology has already been used in Section~\ref{sec:surfaces}. 

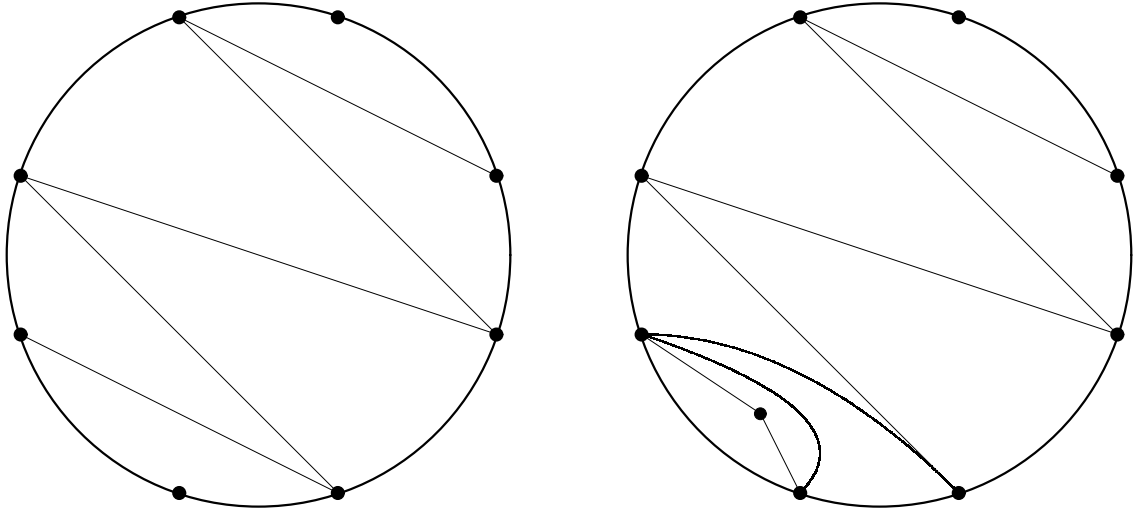
\begin{figure}[htbp] 
\begin{center} 
\setlength{\unitlength}{3pt} 
\begin{picture}(60,63)(0,-1.5) 
\thicklines 
\put(30,30){\circle{63.5}} 
 
  \multiput(20,0)(20,0){2}{\circle*{1.5}} 
  \multiput(20,60)(20,0){2}{\circle*{1.5}} 
  \multiput(0,20)(0,20){2}{\circle*{1.5}} 
  \multiput(60,20)(0,20){2}{\circle*{1.5}} 
 
\thinlines 
  \multiput(0,40)(20,20){2}{\line(1,-1){40}} 
\put(0,20){\line(2,-1){40}} \put(0,40){\line(3,-1){60}} 
\put(20,60){\line(2,-1){40}}

\end{picture} 
\qquad\qquad
\setlength{\unitlength}{3pt} 
\begin{picture}(60,63)(0,-1.5) 
\thicklines 
\put(30,30){\circle{63.5}} 
 
  \multiput(20,0)(20,0){2}{\circle*{1.5}} 
  \multiput(20,60)(20,0){2}{\circle*{1.5}} 
  \multiput(0,20)(0,20){2}{\circle*{1.5}} 
  \multiput(60,20)(0,20){2}{\circle*{1.5}} 
 
\thinlines 
  \multiput(0,40)(20,20){2}{\line(1,-1){40}} 
\qbezier(0,20)(20,20)(40,0)
\put(0,40){\line(3,-1){60}} 
\put(20,60){\line(2,-1){40}} 
 
\put(20,0){\line(-1,2){5}} 
\qbezier(0,20)(30,10)(20,0)
\put(0,20){\line(3,-2){15}} 
 
\put(15,10){\circle*{1.5}}

\end{picture} 
\end{center} 
\caption{Triangulations for types $A$ and $D$} 
\label{fig:octagon-snake} 
\end{figure}

\begin{example}[\emph{Type $A_n$: A polygon}]
Let $\SM$ be an unpunctured $(n+3)$-gon ($n\geq 1$). 
Then $\Bcal\SM$ is a mutation equivalence class of type~$A_n$. 
This can be verified by looking at a triangulation of the form shown
in Figure~\ref{fig:octagon-snake} on the left. 

As shown in~\cite[Section~12.2]{ca2}, a cluster algebra of this type
can be constructed by taking (the $\ZZ$-form of) the homogeneous
coordinate ring of the Grassmannian ${\rm Gr}_{2,n+3}$ of
$2$-dimensional subspaces of an $(n+3)$-dimensional complex vector
space, with respect to its Pl\"ucker embedding. 
Another example is the coordinate ring of the affine base space $\SL_4(\CC)/N$
(type~$A_3$). (Here and below, $N$ denotes the subgroup of unipotent
upper-triangular matrices.) 
\end{example}

\begin{example}[\emph{Type $D_n$: Once-punctured polygon}]
Let $\SM$ be a once-punctured $n$-gon ($n\geq 4$). 
Then $\Bcal\SM$ is a mutation equivalence class of type~$D_n$. 
This can be verified by looking at a triangulation of the form shown
in Figure~\ref{fig:octagon-snake} on the right. 

Examples of cluster algebras of type~$D$ include the (homogeneous) 
coordinate rings of the Grassmannian ${\rm Gr}_{3,6}$ (type~$D_4$),
the special linear group $\SL_3(\CC)$ viewed as a surface in~$\CC^9$
(type~$D_4$), the affine base space $\SL_5(\CC)/N$ (type~$D_6$),
and (the affine cone over) the Schubert divisor in the Grassmannian
${\rm Gr}_{2,n+2}$ (type $D_n$)~\cite[Example~12.15]{ca2}.
\end{example}

\begin{figure}[htbp] 
\vspace{-.2in} 
\[ 
\begin{array}{ccl} 
\widetilde A_1 && 
\setlength{\unitlength}{1.5pt} 
\begin{picture}(20,17)(0,-2) 
\put(0,1){\line(1,0){20}} 
\put(0,-1){\line(1,0){20}} 
\multiput(0,0)(20,0){2}{\circle*{2}} 
\end{picture}\\[.1in] 
\begin{array}{c}
\widetilde A_{n-1}\\
{\scriptstyle (n\geq 3)}
\end{array} && 
\setlength{\unitlength}{1.5pt} 
\begin{picture}(120,22)(0,18) 
\put(0,30){\line(1,0){50}} 
\put(70,30){\line(1,0){50}} 
\put(0,30){\line(3,-1){60}} 
\put(60,10){\line(3,1){60}} 
\put(60,10){\circle*{2}} 
\multiput(0,30)(20,0){3}{\circle*{2}} 
\multiput(80,30)(20,0){3}{\circle*{2}} 
\multiput(56,30)(4,0){3}{\circle*{0.5}} 
\end{picture}\\[.2in] 
\widetilde D_ {n-1} && 
\setlength{\unitlength}{1.5pt} 
\begin{picture}(100,17)(0,-2) 
\put(20,0){\line(1,0){30}} 
\put(70,0){\line(1,0){30}} 
\put(0,10){\line(2,-1){20}} 
\put(0,-10){\line(2,1){20}} 
\put(100,0){\line(2,-1){20}} 
\put(100,0){\line(2,1){20}} 
\multiput(20,0)(20,0){2}{\circle*{2}} 
\multiput(80,0)(20,0){2}{\circle*{2}} 
\put(0,10){\circle*{2}} 
\put(0,-10){\circle*{2}} 
\put(120,10){\circle*{2}} 
\put(120,-10){\circle*{2}} 
\multiput(56,0)(4,0){3}{\circle*{0.5}} 
\end{picture} 
\\[.25in] 
\end{array} 
\] 
\vspace{-.1in} 
\caption{Some affine Coxeter-Dynkin diagrams} 
\label{fig:affine-dynkin} 
\end{figure}
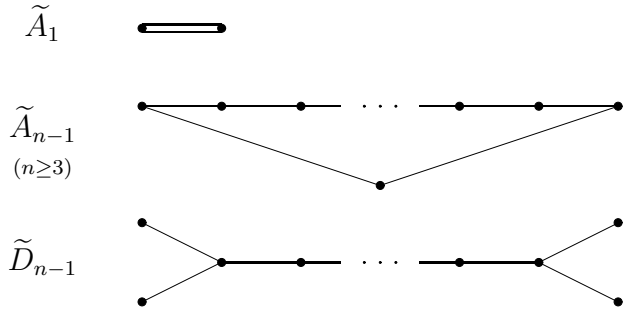

By analogy with finite type, one can use (orientations of) affine Coxeter-Dynkin
diagrams to define mutation equivalence classes of ``affine types.''
Some examples are shown in Figure~\ref{fig:affine-dynkin}. 
In the case of affine type $\widetilde A_{n-1}$ ($n\geq 3$),
the underlying graph (an $n$-cycle) 
is not a forest, so Lemma~\ref{lem:forest} does
not apply. In fact, the mutation equivalence class in this case 
does depend on the orientation of the edges in the diagram; see
Lemma~\ref{lem:affine-a} below. 

\begin{lemma}
\label{lem:affine-a}
Let $\Gamma$ and $\Gamma'$ be two $n$-cycles ($n\geq 3$) whose edges have
been oriented, so that in~$\Gamma$ (resp.,~$\Gamma'$),  
there are $n_1$ (resp.,~$n_1'$) edges of one direction 
and $n_2=n-n_1$ (resp., $n_2'=n-n_1'$) edges of the opposite direction. 
Then the matrices $B(\Gamma)$ and $B(\Gamma')$ 
are mutation equivalent 
(up to simultaneous relabeling of rows and columns) 
if and only if the unordered pairs $\{n_1,n_2\}$ and $\{n_1',n_2'\}$
coincide. 
If $n_1=0$ or $n_1=n$, then $B(\Gamma)$ has type~$D_n$. 
\end{lemma}

The ``only if'' part of the lemma is proved in
Section~\ref{sec:recover-topology}. 

\begin{proof}[Proof (the ``if'' direction)] 
Suppose that $\{n_1,n_2\}=\{n_1',n_2'\}$. It is then easy to check
that $\Gamma$ can be transformed into $\Gamma'$ by
\emph{shape-preserving mutations} (see \cite[Section~9.2]{ca2}),
i.e., mutations performed at sinks or sources of the diagram. 
\end{proof} 

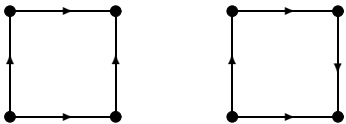
\begin{figure}[htbp] 
\begin{center}
\setlength{\unitlength}{2pt} 
\begin{picture}(20,23)(0,-2) 
\multiput(0,0)(0,20){2}{\line(1,0){20}} 
\multiput(0,0)(20,0){2}{\line(0,1){20}} 
\multiput(0,0)(20,0){2}{\circle*{2}} 
\multiput(0,20)(20,0){2}{\circle*{2}} 
\put(0,7){\vector(0,1){5}} 
\put(20,7){\vector(0,1){5}} 
\put(7,0){\vector(1,0){5}} 
\put(7,20){\vector(1,0){5}} 
\end{picture} 
\hspace{.5in}
\begin{picture}(20,23)(0,-2) 
\multiput(0,0)(0,20){2}{\line(1,0){20}} 
\multiput(0,0)(20,0){2}{\line(0,1){20}} 
\multiput(0,0)(20,0){2}{\circle*{2}} 
\multiput(0,20)(20,0){2}{\circle*{2}} 
\put(0,7){\vector(0,1){5}} 
\put(20,13){\vector(0,-1){5}} 
\put(7,0){\vector(1,0){5}} 
\put(7,20){\vector(1,0){5}} 
\end{picture} 
\end{center}
\vspace{-.1in} 
\caption{Diagrams of types $\widetilde A(2,2)$ and $\widetilde A(3,1)$} 
\label{fig:tilde-a22-a13} 
\end{figure}

In the situation described in Lemma~\ref{lem:affine-a}, 
with $n_1\geq n_2>0$, we say that the
mutation equivalence class of~$B(\Gamma)$ (or a cluster algebra
with this exchange matrix) is of \emph{type~$\widetilde A(n_1,n_2)$}.
See Figure~\ref{fig:tilde-a22-a13}. 

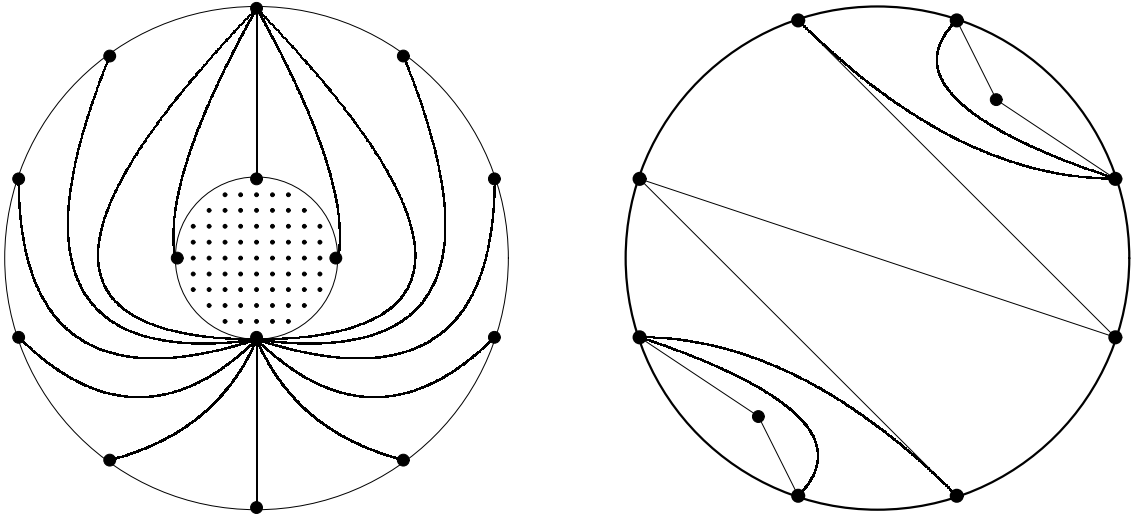
\begin{figure}[htbp] 
\begin{center} 
\setlength{\unitlength}{3pt} 
\begin{picture}(60,63)(0,0) 
\thinlines 

\put(30,30){\circle{63.5}} 

\put(30,30){\circle{20.5}}

  \put(30,-1.5){\circle*{1.5}} 
  \put(11.5,4.5){\circle*{1.5}} 
  \put(48.5,4.5){\circle*{1.5}} 
  \put(30,61.5){\circle*{1.5}} 
  \put(11.5,55.5){\circle*{1.5}} 
  \put(48.5,55.5){\circle*{1.5}} 
  \multiput(0,20)(0,20){2}{\circle*{1.5}} 
  \multiput(60,20)(0,20){2}{\circle*{1.5}} 
 
\put(30,20){\circle*{1.5}} 
\put(30,40){\circle*{1.5}} 
\put(20,30){\circle*{1.5}} 
\put(40,30){\circle*{1.5}} 

\multiput(26,22)(2,0){5}{\circle*{0.5}} 
\multiput(24,24)(2,0){7}{\circle*{0.5}} 
\multiput(22,26)(2,0){9}{\circle*{0.5}} 
\multiput(22,28)(2,0){9}{\circle*{0.5}} 
\multiput(22,30)(2,0){9}{\circle*{0.5}} 
\multiput(22,32)(2,0){9}{\circle*{0.5}} 
\multiput(22,34)(2,0){9}{\circle*{0.5}} 
\multiput(24,36)(2,0){7}{\circle*{0.5}} 
\multiput(26,38)(2,0){5}{\circle*{0.5}} 

\put(30,-1.5){\line(0,1){21.5}}
\put(30,61.5){\line(0,-1){21.5}}
\qbezier(30,61.5)(18,40)(19.7,30)
\qbezier(30,61.5)(42,40)(40.3,30)
\qbezier(30,61.5)(-10,20)(30,19.7)
\qbezier(30,61.5)(70,20)(30,19.7)

\qbezier(11.5,4.5)(25,8)(30,19.7)
\qbezier(48.5,4.5)(35,8)(30,19.7)

\qbezier(0,20)(15,5)(30,19.7)
\qbezier(60,20)(45,5)(30,19.7)

\qbezier(0,40)(0,10)(30,19.7)
\qbezier(60,40)(60,10)(30,19.7)

\qbezier(11.5,55.5)(-5,15)(30,19.7)
\qbezier(48.5,55.5)(65,15)(30,19.7)

\end{picture} 
\qquad\qquad
\setlength{\unitlength}{3pt} 
\begin{picture}(60,63)(0,0) 
\thicklines 
\put(30,30){\circle{63.5}} 
 
  \multiput(20,0)(20,0){2}{\circle*{1.5}} 
  \multiput(20,60)(20,0){2}{\circle*{1.5}} 
  \multiput(0,20)(0,20){2}{\circle*{1.5}} 
  \multiput(60,20)(0,20){2}{\circle*{1.5}} 
 
\thinlines 
  \multiput(0,40)(20,20){2}{\line(1,-1){40}} 
\qbezier(0,20)(20,20)(40,0)
\put(0,40){\line(3,-1){60}} 
 
\put(20,0){\line(-1,2){5}} 
\qbezier(0,20)(30,10)(20,0)
\put(0,20){\line(3,-2){15}} 
 
\put(15,10){\circle*{1.5}} 
 
\put(40,60){\line(1,-2){5}} 
\qbezier(60,40)(30,50)(40,60)
\put(60,40){\line(-3,2){15}} 
\qbezier(60,40)(40,40)(20,60)

\put(45,50){\circle*{1.5}} 
\end{picture} 

\end{center} 
\caption{Triangulations for types $\widetilde A(10,4)$ 
and $\widetilde D_{10}$}
\label{fig:a-d-affine} 
\end{figure}

\begin{example}[\emph{Type $\widetilde A(n_1,n_2)$: An annulus}]
\label{example:cluster-type-A-affine}
Let $\SM$ be an unpunctured annulus with $n_1$ marked points on one
boundary component and $n_2$ on another. 
Then $\Bcal\SM$ is a mutation equivalence class of type~$\widetilde
A(n_1,n_2)$. 
This can be verified by looking at a triangulation of the form shown
in Figure~\ref{fig:a-d-affine} on the left.

In the special case of type $\widetilde A(2,1)$, we recover the example
illustrated in Figures~\ref{fig:Arc-complex-A(2,1)} and~\ref{fig:A(2,1)}.

Type $\widetilde A(2,2)$ can also be obtained by taking
a twice-punctured monogon; see Example~\ref{example:monogon2}.

Coordinate rings of certain double Bruhat cells and 
Schubert varieties/cells in (affine cones over) Grassmannians have natural
cluster algebra structures of type $\widetilde A(n_1,n_2)$.
For instance, cluster type  $\widetilde A(6,3)$ can be realized by
coordinate rings of
hypersurfaces in $\operatorname{Gr}_{4,8}$ and/or  
$\SL_4$ (embedding~$\SL_4$ as the hypersurface 
$\det=1$ in~$\CC^{16}$).
\end{example}

\begin{example}[\emph{Type $\widetilde D_{n-1}$: Twice-punctured polygon}]
Let $\SM$ be a twice-punctured $(n-3)$-gon, with $n\geq 5$. 
Then $\Bcal\SM$ is a mutation equivalence class of type~$\widetilde D_{n-1}$. 
This can be verified by looking at a triangulation of the form shown
in Figure~\ref{fig:a-d-affine} on the right. 

Examples of cluster algebras of this type include 
coordinate rings of double Bruhat cells, Schubert varieties, etc. 
For instance, cluster type~$\widetilde D_8$ can be realized by the 
coordinate ring of a hypersurface $P=1$ in $\operatorname{Gr}_{3,9}$,
where $P$ is a Pl\"ucker coordinate. On the combinatorial level, the
natural cluster structure on $\operatorname{Gr}_{3,9}$ can be represented
by a $2\times 5$ grid graph on $10$ vertices, with all $4$-cycles properly 
oriented (see~\cite{scott}); removing (or "freezing") a vertex in the 
middle yields the aforementioned hypersurface.
\end{example}

\pagebreak[3]

\section{Tagged arc complexes} 
\label{sec:tagged-arcs} 

In view of Proposition~\ref{pr:BSM-well-defined}, associated with each
marked surface~$\SM$ there is a class of cluster algebras $\Acal$
whose exchange matrices include the sign adjacency matrices of
triangulations of $\SM$.  Proposition~\ref{pr:mut-tri} suggests (and
we will later prove) that the arc
complex~$\DDSM$ is a subcomplex of the cluster complex of~$\Acal$, and
its dual graph~$\EESM$ is a subgraph of the exchange graph of~$\Acal$.
(In general, one does not get the full picture because some arcs in some
triangulations are not flippable.)  We will now present a
combinatorial construction for the entire cluster complex and exchange
graph.

\begin{definition}[\emph{Tagged arcs}]
\label{def:tagged-arcs}
Each arc $\gamma$ in $\SM$ has two \emph{ends} obtained by arbitrarily 
cutting~$\gamma$ into three pieces, then throwing out the middle one. 
We think of the two ends as locations near the
endpoints to be used for labeling (``tagging'') an arc. 
A \emph{tagged arc} is an arc in which each end has been tagged in
one of two ways, 
\emph{plain} or \emph{notched}, so that the following conditions are
satisfied: 
\begin{itemize}
\item
the arc does not cut out a once-punctured monogon; 
\item
an endpoint lying on the boundary is tagged plain; and
\item
both ends of a loop are tagged in the same way.
\end{itemize}
\end{definition} 

In the figures, the plain tags are omitted while the notched tags are
represented by the $\notch$ symbol.

\begin{figure}[htbp] 
\begin{center} 
\setlength{\unitlength}{1.5pt} 
\begin{picture}(40,38)(0,-5) 
\thicklines 


\qbezier(20,0)(8,32)(20,32)
\qbezier(20,0)(32,32)(20,32)

\multiput(20,0)(0,20){2}{\circle*{2}} 
\put(10,10){\makebox(0,0){$\gamma$}} 
\put(20,-5){\makebox(0,0){$a$}} 
\put(20,25){\makebox(0,0){$b$}} 
\end{picture} 
\qquad\qquad
\begin{picture}(40,38)(0,-5) 
\thicklines 

\put(20,0){\line(0,1){20}}
\put(20,16){\makebox(0,0){$\notch$}}
\put(10,8){\makebox(0,0){$\tau(\gamma)$}} 

\multiput(20,0)(0,20){2}{\circle*{2}} 
\put(20,-5){\makebox(0,0){$a$}} 
\put(20,25){\makebox(0,0){$b$}} 
\end{picture} 

\end{center} 
\caption{Representing an arc bounding a punctured monogon by a tagged arc} 
\label{fig:arc-as-tagged-arc} 
\end{figure}
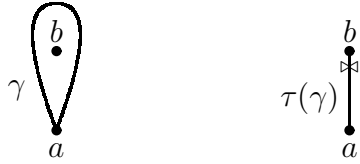

\begin{definition}
\label{def:plain-arc-as-tagged-arc}
Each ordinary (``plain'') arc~$\gamma$ can be represented by a 
tagged arc~$\tau(\gamma)$, as follows. If $\gamma$ does not cut out
a once-punctured monogon, then~$\tau(\gamma)$ is $\gamma$ with
both ends tagged plain.
Otherwise, let $\gamma$ be a loop, based at a marked point~$a$, 
cutting out a punctured monogon with the sole puncture~$b$ inside it;  
see Figure~\ref{fig:arc-as-tagged-arc}.  Note that there is
only one possible monogon since $\SM$ is not the 3-punctured sphere
(by assumption)
so $\gamma$ cannot bound a punctured monogon on both sides.
Let $\beta$ be the unique arc connecting $a$ and $b$ and compatible
with~$\gamma$. Then $\tau(\gamma)$ is
obtained by tagging $\beta$ plain at~$a$ and notched at~$b$. 
\end{definition}

The set of all tagged arcs in $\SM$ is denoted by~$\ATSM$. 
As we just explained, there is a canonical map $\tau$ from the set of
plain arcs~$\AASM$ into $\ATSM$. 

\begin{remark}
Each tagged arc $\gamma\in\ATSM$ belongs to one of the following
classes: 
\begin{itemize}
\item
$\gamma$ is a plain arc that connects different marked points on the
  boundary; 
\item
$\gamma$ is a plain arc that connects a marked point on the boundary to itself and 
does not enclose, on either side, a once-punctured
monogon; 
\item
$\gamma$ connects a puncture to itself and does not enclose,
on either side, a once-punctured monogon;
$\gamma$ may be either plain or with both ends notched; 
\item
$\gamma$ connects a marked point on the boundary to a puncture
(two tagged arcs per untagged isotopy class); 
\item
$\gamma$ connects two punctures at different locations
(four tagged arcs per untagged isotopy class).
\end{itemize}
Table~\ref{table:tagged-arcs} shows, for each of several examples of
bordered surfaces with marked points, 
the cardinalities of these five classes of tagged arcs. 
\end{remark}

\begin{table}[ht]
{\Small
\begin{tabular}{lccccc}
\toprule
$\SM$ &
\begin{tabular}{@{}c@{}}Boundary\\to boundary\\(not a loop)\end{tabular} &
\begin{tabular}{@{}c@{}}Loop at\\boundary\end{tabular} &
\begin{tabular}{@{}c@{}}Loop at\\puncture\end{tabular} &
\begin{tabular}{@{}c@{}}Boundary\\to puncture\end{tabular} &
\begin{tabular}{@{}c@{}}Puncture\\to puncture\\(not a loop)\end{tabular}
\\
\midrule
$(n+3)$-gon
& $\dfrac{n(n+3)}{2}$ & $0$ & $0$ & $0$ & $0$  \\
$n$-gon, $1$ puncture
& $n^2-2n$ & $0$ & $0$ & $2n$ & $0$ \\
annulus $\tilde A(n_1,n_2)$;  
$n_1,n_2 \!\ge\!2$
& $\infty$ & $n_1+n_2$ & $0$ & $0$ & $0$  \\
annulus $\tilde A(n_1,1)$;  
$n_1\!\ge\!2$
& $\infty$ & $n_1$ & $0$ & $0$ & $0$  \\
annulus $\tilde A(1,1)$
& $\infty$ & $0$ & $0$ & $0$ & $0$  \\
$(n\!-\!3)$-gon, $2$ punctures; $n\!\ge\!5$
& $\infty$ & $n-3$ & $0$ & $\infty$ & $4$\\
monogon, $2$ punctures  
& $0$ & $0$ & $0$ & $\infty$ & $4$\\
torus, $1$ puncture & $0$ & $0$ & $\infty$ & $0$ & $0$\\
\bottomrule
\end{tabular}
}
\caption{Tagged arcs of different type}
\label{table:tagged-arcs}
\end{table}

We next adapt the concept of compatibility of arcs 
(see Definition~\ref{def:compat-arcs}) 
to the tagged setup. 

\begin{definition}[\emph{Compatibility of tagged arcs}]
\label{def:compat-tagged-arcs}
Two tagged arcs $\alpha,\beta\in\ATSM$ are called
\emph{compatible} if and only if the following conditions are satisfied:
\begin{itemize}
\item
the untagged versions of $\alpha$ and $\beta$ are compatible; 
\item
if the untagged versions of $\alpha$ and $\beta$ are different, 
and $\alpha$ and $\beta$ share an endpoint~$a$, 
then the ends of $\alpha$ and $\beta$ connecting to~$a$ 
must be tagged in the same way;
\item 
if the untagged versions of $\alpha$ and $\beta$ coincide, 
then at least one end of $\alpha$ must be tagged in the same way 
as the corresponding end of~$\beta$. 
\end{itemize}
\end{definition} 
\begin{remark}
The following is a complete list of possible compatible pairs~$\{\alpha,\beta\}$ 
of tagged arcs:
\begin{itemize}
\item
$\alpha$ and $\beta$ have endpoints at four different locations, and do
not intersect; 
\item
$\alpha$ and $\beta$ connect from different locations to the same
boundary point, and do not intersect in the interior of~$\Surf$; 
\item
$\alpha$ and $\beta$ connect from different locations to the same
puncture, are tagged in the same way at it, and do not intersect
 in the interior of~$\Surf$; 
\item
$\alpha$ and $\beta$ are not loops, have the same endpoints,
are tagged in the same way at their respective ends, and do not intersect
 in the interior of~$\Surf$; 
\item
$\alpha$ and $\beta$ are two different loops based at the same
marked point,
are tagged in the same way at all four ends, and do not intersect 
in the interior of~$\Surf$; 
\item
$\alpha$ and $\beta$ are not loops, their untagged versions coincide,
they are tagged in the same way at one end and
differently at another; 
\item
$\alpha$ and $\beta$ coincide. 
\end{itemize}
\end{remark}

\pagebreak[3]

\begin{remark}
\label{rem:compat-plain/tagged}
If two plain arcs~$\alpha$ and $\beta$ are compatible, then
the tagged arcs~$\tau(\alpha)$ and $\tau(\beta)$ representing
them (in the sense of Definition~\ref{def:plain-arc-as-tagged-arc})
are compatible.
The converse is \emph{false}: in a once-punctured digon, 
the two loops are not compatible even though the tagged arcs
representing them are compatible. 
On the other hand, the two notions coincide for arcs that do not cut
out a once-punctured monogon. 
\end{remark}

\begin{definition}[\emph{Tagged arc complex}]
\label{def:tagged-arc-complex}
As before, let $\SM$ be a bordered surface with marked points. 
The \emph{tagged arc complex} $\DTSM$ is 
the clique complex for the compatibility relation on 
the set of tagged arcs ~$\ATSM$. 
That is, $\DTSM$ is the simplicial
complex on the ground set $\ATSM$ whose
simplices are collections of mutually compatible tagged arcs. 
In the absence of punctures, the complex $\DTSM$ obviously
coincides with the ordinary arc complex~$\DDSM$ of
Definition~\ref{def:arc-complex}. 
\end{definition}


\newsavebox{\digon}
\setlength{\unitlength}{4pt} 
\savebox{\digon}(10,10)[bl]{
\qbezier(5,0)(0,5)(5,10)
\qbezier(5,0)(10,5)(5,10)
\put(5,0){\circle*{1}} 
\put(5,5){\circle*{1}} 
\put(5,10){\circle*{1}} 
}

\newsavebox{\hdigon}
\setlength{\unitlength}{4pt} 
\savebox{\hdigon}(10,10)[bl]{
\qbezier(5,0)(0,5)(5,10)
\qbezier(5,0)(10,5)(5,10)
\put(5,0){\circle*{1}} 
\put(5,5){\circle{1.5}} 
\put(5,10){\circle*{1}} 
}

\newsavebox{\lowbar}
\savebox{\lowbar}(10,10)[bl]{
\put(5,0){\line(0,1){5}}
}

\newsavebox{\highbar}
\setlength{\unitlength}{4pt} 
\savebox{\highbar}(10,10)[bl]{
\put(5,5){\line(0,1){5}}
}

\newsavebox{\lowloop}
\setlength{\unitlength}{4pt} 
\savebox{\lowloop}(10,10)[bl]{
\qbezier(5,0)(2,7)(5,7)
\qbezier(5,0)(8,7)(5,7)
}

\newsavebox{\highloop}
\setlength{\unitlength}{4pt} 
\savebox{\highloop}(10,10)[bl]{
\qbezier(5,10)(2,3)(5,3)
\qbezier(5,10)(8,3)(5,3)
}

\newsavebox{\lowtag}
\setlength{\unitlength}{4pt} 
\savebox{\lowtag}(10,10)[bl]{
\put(5,3.5){\makebox(0,0){$\notch$}}
}

\newsavebox{\hightag}
\setlength{\unitlength}{4pt} 
\savebox{\hightag}(10,10)[bl]{
\put(5,6.5){\makebox(0,0){$\notch$}}
}

\begin{figure}
  \centering
  \includegraphics{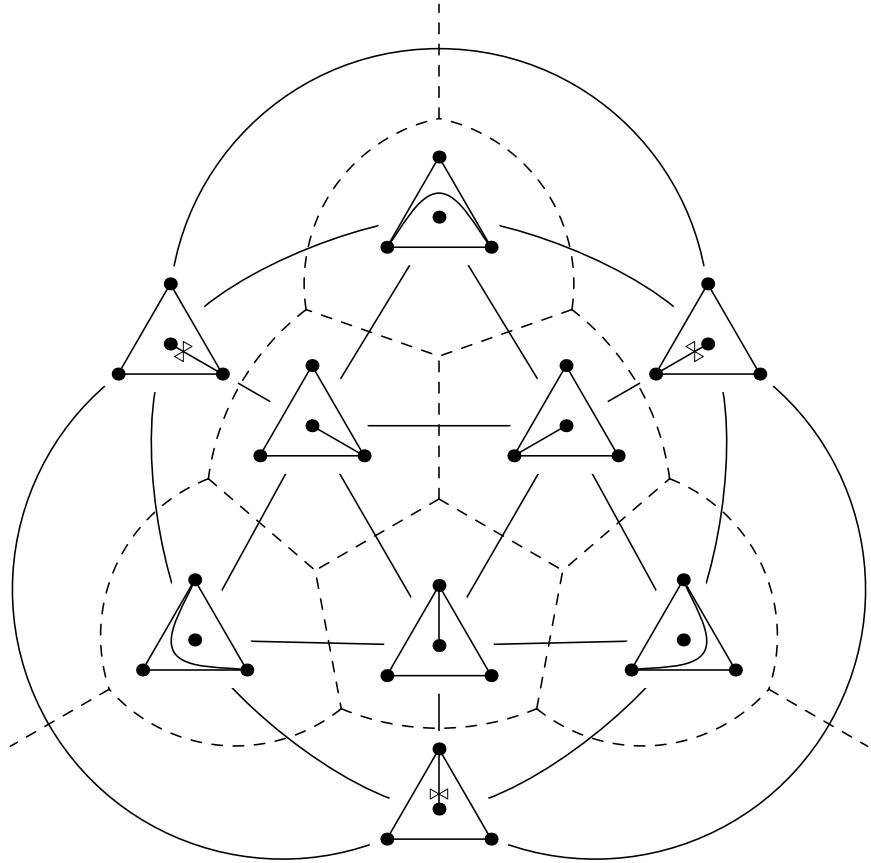}
  \caption{The tagged arc complex of a once-punctured triangle}
  \label{fig:tagged-arc-d3}
\end{figure}

\begin{example}
  Figure~\ref{fig:tagged-arc-d3} shows the tagged arc complex of a
  punctured triangle.  In view of
  Remark~\ref{rem:compat-plain/tagged}, the untagged arc
  complex~$\DDSM$ can be viewed as
  a subcomplex of the tagged arc complex~$\DTSM$, but $\DDSM$ is
  \emph{not} an induced subcomplex of $\DTSM$, as you can see by
  comparing Figures~\ref{fig:arc-complex-D3}
  and~\ref{fig:tagged-arc-d3}.
\end{example}

A maximal (by inclusion) collection of pairwise
compatible tagged arcs in~$\SM$ is called a \emph{tagged
  triangulation}. 

\begin{theorem}
\label{th:tagged-pseudo}
The tagged arc complex $\DTSM$ is a pseudomanifold.
That is, each maximal simplex in $\DTSM$ (a tagged triangulation) has the
same cardinality (equal to the rank of~$\SM$), 
and each simplex of codimension~$1$ is contained in
precisely two maximal simplices.  
\end{theorem}

Theorem~\ref{th:tagged-pseudo} is proved in
Section~\ref{sec:basic-properties}.

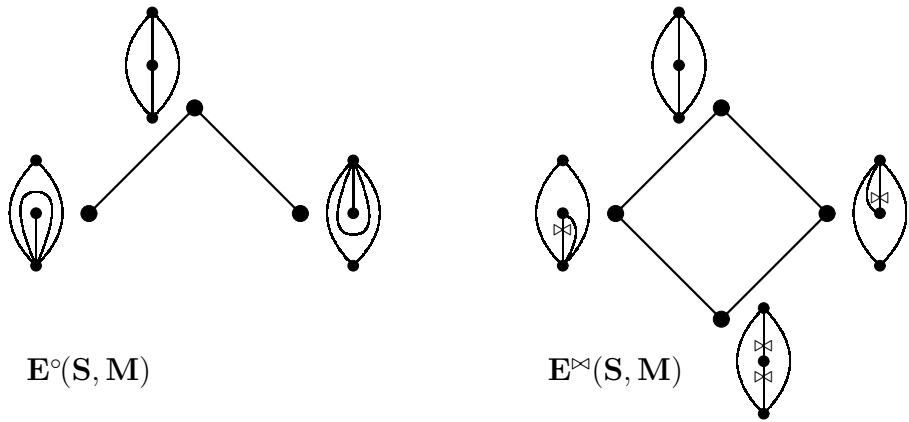
\begin{figure}[htbp] 
\begin{center} 

\setlength{\unitlength}{4pt} 
\begin{picture}(40,37)(5,-3) 

\thinlines 


\put(40,15){\makebox(0,0){\usebox{\digon}}} 
\put(40,15){\makebox(0,0){\usebox{\highbar}}} 
\put(40,15){\makebox(0,0){\usebox{\highloop}}} 

\put(10,15){\makebox(0,0){\usebox{\digon}}} 
\put(10,15){\makebox(0,0){\usebox{\lowbar}}} 
\put(10,15){\makebox(0,0){\usebox{\lowloop}}}

\put(21,29){\makebox(0,0){\usebox{\digon}}} 
\put(21,29){\makebox(0,0){\usebox{\lowbar}}} 
\put(21,29){\makebox(0,0){\usebox{\highbar}}}

\put(25,25){\circle*{1.5}} 
\put(15,15){\circle*{1.5}} 
\put(35,15){\circle*{1.5}} 

\thicklines

\put(15,15){\line(1,1){10}}
\put(25,25){\line(1,-1){10}}

\put(15,0){\makebox(0,0){$\EESM$}}

\end{picture} 
\qquad\quad
\begin{picture}(40,37)(5,-3) 

\thinlines 

\put(29,1){\makebox(0,0){\usebox{\digon}}} 
\put(29,1){\makebox(0,0){\usebox{\lowbar}}} 
\put(29,1){\makebox(0,0){\usebox{\lowtag}}} 
\put(29,1){\makebox(0,0){\usebox{\highbar}}} 
\put(29,1){\makebox(0,0){\usebox{\hightag}}} 

\put(40,15){\makebox(0,0){\usebox{\digon}}} 
\qbezier(40,15)(37.5,16)(40,20)
\put(40,15){\makebox(0,0){\usebox{\highbar}}} 
\put(40,15){\makebox(0,0){\usebox{\hightag}}} 

\put(10,15){\makebox(0,0){\usebox{\digon}}} 
\qbezier(10,15)(12.5,14)(10,10)
\put(10,15){\makebox(0,0){\usebox{\lowbar}}} 
\put(10,15){\makebox(0,0){\usebox{\lowtag}}}

\put(21,29){\makebox(0,0){\usebox{\digon}}} 
\put(21,29){\makebox(0,0){\usebox{\lowbar}}} 
\put(21,29){\makebox(0,0){\usebox{\highbar}}}

\put(25,5){\circle*{1.5}} 
\put(25,25){\circle*{1.5}} 
\put(15,15){\circle*{1.5}} 
\put(35,15){\circle*{1.5}} 

\thicklines

\put(25,5){\line(1,1){10}}
\put(25,5){\line(-1,1){10}}
\put(15,15){\line(1,1){10}}
\put(25,25){\line(1,-1){10}}

\put(15,0){\makebox(0,0){$\ETSM$}}

\end{picture} 

\end{center} 
\caption{Graphs $\EESM$ and $\ETSM$ for a once-punctured digon} 
\label{fig:exch-graph-A1A1} 
\end{figure}

We denote by $\ETSM$ the dual graph of~$\DTSM$.
We will soon see that $\ETSM$ almost always coincides with
the exchange graph~$\Exch\SM$ 
of a cluster algebra~$\Acal$ with $\Bcal(\Acal)=\Bcal\SM$
(see Definition~\ref{def:exchange-graph}).
Again, if there are no punctures, then $\ETSM$ is just the
dual graph $\EESM$ of the ordinary arc complex. 
By analogy with the ordinary case,
the edges of $\ETSM$ are called \emph{flips}.
Thus, each flip replaces a tagged arc in a tagged triangulation by a
(uniquely defined, different) new tagged arc. 
See Figure~\ref{fig:tagged-arc-d3} (dashed lines) and
Figure~\ref{fig:exch-graph-A1A1}.

\begin{proposition}
\label{pr:dual-tagged-connected}
If $\SM$ is not a closed surface with exactly one puncture
then the graph $\ETSM$ is connected, i.e., 
any two tagged triangulations are connected by a
sequence of flips. 
Consequently, $\DTSM$ is connected.

If $\SM$ is a closed surface with one puncture, then $\ETSM$ and
$\DTSM$ each have two isomorphic components: one in which all ends of all arcs
are plain and one in which they are all notched.
\end{proposition}

The following is our main result, proved in Section~\ref{sec:proofs}.

\begin{theorem}
\label{th:cluster-cpx=arc-cpx}
Assume that $\SM$ is not a closed surface with 
one or two punctures.
Let $\Acal$ be a cluster algebra with $\Bcal(\Acal)=\Bcal\SM$. 
Then the cluster complex~$\Delta\SM$ of~$\Acal$ is isomorphic to the tagged arc
complex~$\DTSM$, and the exchange graph of~$\Acal$ is isomorphic to~$\ETSM$.

If $\SM$ is a closed surface with exactly one puncture, $\Delta\SM$ is
isomorphic to a connected component of $\DTSM$, and the exchange graph
of~$\Acal$ is isomorphic to a connected component of~$\ETSM$.
\end{theorem}

The case of a two-punctured closed surface will be treated
later~\cite{fst-hyper}.

Thus, the cluster variables in~$\Acal$ can be labeled by the tagged
arcs in~$\SM$ so that two cluster variables are compatible if and only
if the corresponding tagged arcs are compatible.  If $\SM$ is a closed
surface with exactly one puncture, we only use the plain arcs.

By Theorem~\ref{th:cluster-cpx=arc-cpx}, the tagged arc complex
of~$\SM$ is identified with the cluster complex of a cluster
algebra~$\Acal$ whose seeds are labeled by the tagged triangulations
of~$\SM$. 
The exchange matrix at each seed can be directly described 
in terms of the corresponding tagged triangulation;
the exact combinatorial recipe is given in 
Section~\ref{sec:mutations-tagged}. 

\pagebreak[3]

\section{Denominators and intersection numbers}
\label{sec:intersection-pairing}

One of the most fundamental results of cluster algebra theory is the
following. 

\begin{theorem}[Laurent phenomenon for cluster algebras~\cite{ca1}]
\label{th:laurent-phenom}
Any cluster variable, when expressed as a rational function 
in the elements of a given cluster, is a Laurent polynomial. 
\end{theorem}

It is a natural idea---already used in many papers, starting
with~\cite{ca2}---to use the \emph{denominators} of these Laurent
polynomials to distinguish the cluster variables from each other. 
Let us recall this setup, largely following~\cite[Section~7]{ca4}. 

\begin{definition}[\emph{Denominator vectors}]
Let $\xx_0$ 
be a cluster
from some seed in a cluster algebra~$\Acal$.
By Theorem~\ref{th:laurent-phenom}, any nonzero element
$z \in \Acal$ can be uniquely written as
\begin{equation}
\label{eq:Laurent-normal-form}
z = \frac{P(\xx_0)}{\displaystyle\prod_{x\in\xx_0} x^{d(x|z)}} \, ,
\end{equation}
where $\dd_{\xx_0}(z)=(d(x|z):x\in\xx_0)$ is an $n$-tuple of integers, and 
$P(\xx_0)$ is a polynomial 
not divisible by any cluster variable~$x\in\xx_0$.
The vector $\dd_{\xx_0}(z)$ 
is called the \emph{denominator vector}
of~$z$ with respect to the cluster~$\xx_0$.
\end{definition}

According to \cite[Conjecture~7.4(3)]{ca4}, the number $d(x|z)$ depends only on
$x$ and~$z$ but not on the choice of cluster $\xx_0$ containing~$x$. 
Since we cannot assume this to hold \emph{a priori}, 
we will sometimes use the notation $d_{\xx_0}(x|z)$ to reflect the possible dependence
on~$\xx_0$. 


As noted in~\cite{ca4}, the denominator vectors of cluster variables
(and in fact their entire Newton polytopes) 
do not depend on the choice of coefficients. 
This is because the exchange relations \eqref{eq:exchange-rel-xx} 
imply the following explicit recurrence for the denominator vectors. 

\begin{lemma}[{\cite[Section~7]{ca4}}]
\label{lem:denom-recurrence}
The components $d(x|z)=d_{\xx_0}(x|z)$ of 
denominator vectors of cluster variables 
with respect to a fixed
initial cluster~$\xx_0$ are uniquely
determined by the initial conditions
\begin{equation}
\label{eq:denom-cluster-variable-initial}
\begin{array}{l}
d(x|x) = -1\,,\\
d(x|z)=0 \quad \text{for} \quad z\in\xx_0-\{x\}\,, 
\end{array}
\end{equation}
together with the recurrence relations
\begin{equation}
\label{eq:exchange-denominator}
d(x|\overline z) =
- d(x|z) +
\max\Bigl(\displaystyle\sum_{u\in\xx} \,[b_{uz}]_+\, d(x|u),
\sum_{u\in\xx}\, [-b_{uz}]_+ \,d(x|u) \Bigr) \,, 
\end{equation}
for each pair of adjacent seeds $\Sigma \!=\! (\xx, \pp, B)$ and 
$\mu_z(\Sigma)\!=\! (\overline \xx, \overline \pp,\overline B)$,
with $B\! = \!(b_{uv})_{u,v \in \xx}$ 
and $\overline \xx= \xx - \{z\} \cup \{\overline z\}$
(cf.\ \eqref{eq:exchange-rel-xx}). 
Here we use the notation $[a]_+=\max(a,0)$. 
\end{lemma} 

The relation~\eqref{eq:exchange-denominator} can be viewed as a
\emph{tropical} version of the exchange
relation~\eqref{eq:exchange-rel-xx}.  That is,
\eqref{eq:exchange-denominator} can be obtained from
\eqref{eq:exchange-rel-xx} by replacing ordinary multiplication and
addition by their tropical counterparts, addition and maximum,
respectively. 

We will show that for the class of cluster algebras under our 
consideration, the denominators $d(x|z)$ can be interpreted as 
certain modified intersection numbers.
The existence of such an interpretation was also
suggested in~\cite{gsv2}. 

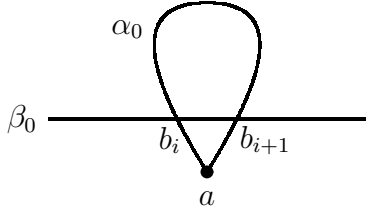
\begin{figure}[htbp] 
\begin{center} 
\setlength{\unitlength}{2pt} 
\begin{picture}(40,35)(0,-2) 
\thicklines 


\qbezier(20,0)(0,32)(20,32)
\qbezier(20,0)(40,32)(20,32)

\put(20,0){\circle*{2}} 
\put(5,27){\makebox(0,0){$\alpha_0$}} 
\put(20,-5){\makebox(0,0){$a$}} 

\put(-10,10){\line(1,0){60}}
\put(-15,10){\makebox(0,0){$\beta_0$}} 

\put(13,6){\makebox(0,0){$b_i$}} 
\put(31,6){\makebox(0,0){$b_{i+1}$}} 

\end{picture} 

\end{center} 
\caption{Computing $B$ in Definition~\ref{def:inters-pairing}} 
\label{fig:loop-cutout}
\end{figure}

\begin{definition}[\emph{Intersection pairing}]
\label{def:inters-pairing}
Let $\alpha$ and $\beta$ be two tagged arcs in~$\SM$.
The \emph{intersection number} $(\alpha|\beta)$ is defined as
follows. 
Let $\alpha_0$ and $\beta_0$ be two non-self-intersecting curves
homotopic rel~$\Mark$ to the untagged versions of $\alpha$
and~$\beta$, respectively, and intersecting transversally in a minimum
number of points in $\Surf\setminus \Mark$. 
Then $(\alpha|\beta)\stackrel{\rm def}{=}A\!+\!B\!+\!C\!+\!D$, where 
\begin{itemize}
\item
$A$ is the number of intersection points of $\alpha_0$ and $\beta_0$
in $\Surf \setminus \Mark$;
\item
$B=0$ unless $\alpha_0$ is a loop based at a marked point~$a$,
in which case $B$ is computed as follows:  
assume that $\beta_0$ intersects~$\alpha_0$ (in the interior of~$\Surf\setminus\Mark$)
at the points $b_1,\dots,b_m$ (numbered along $\beta_0$ in this order); 
then $B$ is the \emph{negative} of the number of segments
$\gamma_i=[b_i,b_{i+1}]\subset\beta_0$ 
such that $\gamma_i$ together with the segments 
$[a,b_i]\subset\alpha_0$ and $[a,b_{i+1}]\subset\alpha_0$
form three sides of a contractible triangle disjoint from the punctures (see
Figure~\ref{fig:loop-cutout});  
\item
$C=0$ unless $\alpha_0=\beta_0$, in which case $C=-1$; 
\item
$D$ is the number of ends of $\beta$ that are incident to 
an endpoint of~$\alpha$ and, at that endpoint, carry a tag different
from the one $\alpha$ does. 
\end{itemize}
\end{definition}

\begin{figure}[htbp]
\begin{center} 
\setlength{\unitlength}{1.5pt} 
\begin{picture}(40,60)(0,-5) 
\thicklines 
\qbezier(20,20)(8,52)(20,52)
\qbezier(20,20)(32,52)(20,52)

\multiput(20,0)(0,20){2}{\circle*{2}} 
\put(10,40){\makebox(0,0){$\alpha$}} 
\put(20,-5){\makebox(0,0){$b$}} 
\put(25,22){\makebox(0,0){$a$}} 

\put(20,0){\line(0,1){20}}
\put(20,16){\makebox(0,0){$\notch$}}
\put(15,8){\makebox(0,0){$\beta$}} 

\multiput(20,0)(0,20){2}{\circle*{2}} 

\end{picture} 

\end{center} 
\caption{\hbox{Example with $(\alpha|\beta) \ne (\beta|\alpha)$}} 
\label{fig:denom-not-symmetric}
\end{figure}
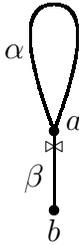

\begin{example}
\label{example:inters-num-nonsymm}
To illustrate Definition~\ref{def:inters-pairing},
let $\alpha$ be a plain loop based at a puncture~$a$, 
and let $\beta$ be an arc connecting $a$ with $b\neq a$,
notched at~$a$ and not having other common points with~$\alpha$.  (See
Figure~\ref{fig:denom-not-symmetric}).
Then $(\alpha|\beta)=1$ while $(\beta|\alpha)=2$. 
\end{example}

{}From now on, we assume that $\SM$ is not a closed surface with 
exactly two punctures,
and $\Acal$ be a cluster algebra with $\Bcal(\Acal)=\Bcal\SM$. 
By Theorem~\ref{th:cluster-cpx=arc-cpx}, the cluster variables
in~$\Acal$ can be labeled by the tagged arcs in~$\SM$.
Let $x[\beta]$ denote the cluster variable labeled by~$\beta$. 

\begin{theorem}
\label{th:denoms-as-inters-numbers}
For any tagged arcs $\alpha$ and $\beta$, the denominator component 
$d(x[\alpha]|x[\beta])$ is equal to the intersection
number~$(\alpha|\beta)$. In particular, $d(x[\alpha]|x[\beta])$
does not depend on the choice of a cluster containing~$x[\alpha]$. 
\end{theorem}

In the special case of ordinary (rather than tagged) arcs, 
Theorem~\ref{th:denoms-as-inters-numbers}
rectifies a similar albeit incomplete statement made in 
\cite[Theorem~3.4]{gsv2}. 

Theorem~\ref{th:denoms-as-inters-numbers} establishes 
\cite[Conjecture~7.4]{ca4} for the class of
  cluster algebras considered in this paper. 

Theorem~\ref{th:denoms-as-inters-numbers} can be restated as saying
that for any tagged arcs $\alpha_1,\dots,\alpha_n$ forming a tagged
triangulation of~$\SM$, and any tagged arc~$\beta$, we have 
\begin{equation}
x[\beta]=\dfrac{P(x[\alpha_1],\dots,x[\alpha_n])}{\displaystyle\prod_{i=1}^n
  x[\alpha]^{(\alpha|\beta)}}\,, 
\end{equation}
where $P$ is a polynomial not divisible by any variable~$x[\alpha_i]$. 

Theorem~\ref{th:denoms-as-inters-numbers} is proved in 
Section~\ref{sec:proofs} by verifying that the intersection numbers
$(\alpha|\beta)$ satisfy the
recurrence~\eqref{eq:exchange-denominator}. 

\pagebreak[3] 

\section{Proofs of main results}
\label{sec:proofs}

In this section, we prove Theorems~\ref{th:cluster-complex-top},
\ref{th:tagged-pseudo},
\ref{th:cluster-cpx=arc-cpx}, \ref{th:denoms-as-inters-numbers}, 
and Proposition~\ref{pr:dual-tagged-connected}. 



\subsection{Basic properties of tagged arc complexes}
\label{sec:basic-properties}

\begin{definition}[\emph{Signature of a partial tagged triangulation}] 
\label{def:sign-partial}
Let $C$ be a simplex in $\DTSM$, i.e., a collection of pairwise compatible
tagged arcs. For a puncture~$a$, let us temporarily denote by $S_C(a) \subset
\{\text{plain},\text{notched}\}$ the set of tags at $a$ that appear
in~$C$.  The \emph{signature} $\delta_{C}$ of~$C$ is the
$(0,\pm1)$-valued function on the set of all punctures in~$\SM$
defined by
\begin{equation}
\delta_{C}(a)=
\begin{cases}
1 & \text{if $S_C(a) = \{\text{plain}\}$;} \\
-1 & \text{if $S_C(a) = \{\text{notched}\}$;} \\
0 & \text{if $S_C(a) = \{\text{plain}, \text{notched}\}$ or $S_C(a) = \emptyset$.} \\
\end{cases}
\end{equation}
Note that if $S_C(a) = \{\text{plain}, \text{notched}\}$, then there are
precisely two arcs in $C$ incident to $a$: the untagged versions of
these two arcs coincide and they are tagged in the same way at the end
different from~$a$.
\end{definition}

\begin{definition}
\label{def:epsilon-tagging}
  Let $\varepsilon$ be a $(\pm1)$-valued function on the set of
  punctures in $\SM$.  For $\gamma$ an (ordinary) arc, define
  $\tau(\gamma,\varepsilon)$ to be the tagged arc obtained from
  $\tau(\gamma)$ (see Definition~\ref{def:plain-arc-as-tagged-arc}) by
  changing 
  the tags at all punctures $a$ with $\varepsilon(a) = -1$ 
  from plain to notched and vice versa.  For a collection~$C$ of compatible
  ordinary arcs, define $\tau(C,\varepsilon) = \{\tau(\gamma, \varepsilon) :
  \gamma \in C\}$. 
  It is clear that the tagged arcs in $\tau(C,\varepsilon)$ are pairwise compatible. 

  Conversely, for $C$ a collection of compatible tagged arcs, let
  $C^\circ$ be the collection of ordinary arcs obtained as follows:
  \begin{itemize}
  \item replace all (notched) tags at the punctures~$a$ with
    $\delta_C(a) = -1$ by plain ones; 
  \item for each puncture $b$ with $\delta_C(b) = 0$, replace the arc
    $\beta$ notched at~$b$ (if any) by a loop~$\gamma$ enclosing~$b$
    and~$\beta$.  That is, $\gamma$~is based at the endpoint $a\neq b$
    of~$\beta$, and is obtained by closely wrapping around~$\beta$, as
    shown in Figure~\ref{fig:arc-as-tagged-arc}.
  \end{itemize}
  It is straightforward to check that the arcs in $C^\circ$ are pairwise
  compatible and distinct.
\end{definition}

\begin{lemma}
\label{lem:detag-tag-inverse}
  Let $C$ be a collection of compatible tagged arcs.  Let
  $\varepsilon$ be any $(\pm1)$-valued function on the set of
  punctures in $\SM$ such that $\varepsilon(a)\delta_C(a) \ge 0$ for
  all punctures~$a$.  Then $\tau(C^\circ, \varepsilon) = C$.
\end{lemma}
\begin{proof}
  Clear from the definitions.
\end{proof}

\begin{proof}[Proof of Theorem~\ref{th:tagged-pseudo}]
We need to show that the tagged arc complex $\DTSM$ is a pseudomanifold.
First, let us prove that $\DTSM$ is pure
of dimension~$n-1$. 
(As before, $n$ denotes the rank of~$\SM$, i.e., 
the cardinality of each maximal simplex in~$\DDSM$.)
For any collection~$C$ of $k$ pairwise compatible tagged arcs,
$C^\circ$ is a collection of $k$ pairwise compatible ordinary arcs.
It follows that any simplex in $\DTSM$ has cardinality at
most~$n$.


Suppose that $C$ is a simplex in $\DTSM$ of cardinality less
than~$n$. We must show that $C$ is not maximal.  Since $C^\circ$ is
not maximal, there exists an arc $\gamma\notin C^\circ$ 
compatible with all arcs in~$C^\circ$.  Let $\varepsilon$ be any
$(\pm1)$-valued function with $\varepsilon(a)\delta_C(a) \ge 0$ for all
punctures $a$.  Then by Lemma~\ref{lem:detag-tag-inverse},
$\tau(\gamma,\varepsilon)$ is compatible with~$C$.

It remains to demonstrate that each codimension~$1$ simplex~$C$ in
$\DTSM$ is contained in exactly two maximal simplices
(tagged triangulations).  Note that if a loop starting at a
puncture~$a$ enclosing a monogon with a single puncture~$b$ appears
in~$C^\circ$, the enclosed arc from $a$ to $b$ also appears
in~$C^\circ$.  Thus for any codimension~$1$ simplex~$C$, the ordinary
collection~$C^\circ$ is not on the boundary of $\DDSM$ and so can be
completed by two different arcs~$\gamma$ and~$\gamma'$.  The tagged arcs
$\tau(\gamma,\varepsilon)$ and $\tau(\gamma',\varepsilon)$ give two
different ways to complete~$C$ to a maximal simplex.
To show that these are the only two completions,
reverse the construction: 
Let $\tilde C\supset C$ be a maximal simplex in $\DTSM$. 
Then $\tilde C^\circ\setminus C^\circ$ consists of a single ordinary
arc compatible with~$C^\circ$, so it must be one of the arcs $\gamma$
and~$\gamma'$. 
Furthermore, we have $\varepsilon(d)\delta_{\tilde C}(d)\ge 0$ 
for all punctures~$d$, so by Lemma~\ref{lem:detag-tag-inverse},
the sole arc in $\tilde C\setminus C$ is either
$\tau(\gamma,\varepsilon)$ or $\tau(\gamma',\varepsilon)$. 
\end{proof}



\subsection{Strata of tagged triangulations} 

We next discuss the relationship between the graphs $\EESM$ 
and $\ETSM$ of ordinary and tagged flips, making use of the
following natural concepts. 

\begin{definition}[\emph{Open and closed strata of tagged triangulations}] 
\label{def:sign-strat}
For a $(0,\pm1)$-function $\delta$ on the set of punctures,   
the \emph{open stratum}~$\Omega_\delta$ is the set of all 
tagged triangulations $T$ with signature $\delta_{T}=\delta$. 
For a $(\pm1)$-function  $\varepsilon$ on the set of punctures, 
the \emph{closed stratum}~$\overline\Omega_\varepsilon$ consists of all 
tagged triangulations $T$ such that $\delta_{T}(a)=\varepsilon(a)$ 
whenever $\delta_{T}(a)=\pm1$.
In other words, $\overline\Omega_\varepsilon$ is the union of all open
strata~$\Omega_\delta$ satisfying $\delta(a)\varepsilon(a)\geq 0$ for 
all  punctures~$a$. 

The ``positive'' closed stratum~$\overline\Omega_{\mathbf{1}}$
associated with the function $\varepsilon=\mathbf{1}$ 
defined by setting $\varepsilon(a)=1$ for all punctures~$a$ is identified with
the set of all ideal triangulations of~$\SM$ by the map $T \mapsto
\tau(T)=\tau(T,\mathbf{1})$.  Accordingly, it makes sense to
talk about the \emph{signature} of an ideal triangulation, 
defined as a $(0,1)$-function on the set of punctures which takes
values $0$ and~$1$ depending on whether a puncture is located inside
a self-folded triangle or not.
\end{definition}

As before, let $p$ denote the number of punctures in~$\SM$.  For each
$(\pm1)$-function~$\varepsilon$ on the set of punctures, the induced
subgraph of $\ETSM$ whose vertex set is the closed
stratum~$\overline\Omega_\varepsilon$ is canonically isomorphic to the
graph~$\EESM$ via the construction of
Definition~\ref{def:epsilon-tagging}.  
Thus, $\ETSM$~can be obtained by gluing together $2^p$
copies of~$\EESM$.

\begin{proof}[Proof of Proposition~\ref{pr:dual-tagged-connected}
(Connectedness of the dual graph $\ETSM$)] 
If $\SM$ is a closed surface with one puncture, then the tagged arcs
forming a simplex in
$\DTSM$ have either all plain ends or all notched ends (by
Definition~\ref{def:tagged-arcs}), so there are two connected
components. Switching between plain and notched tags
gives an isomorphism between the two components.

Otherwise, in view of Proposition~\ref{pr:flips-connected}, it
suffices to show that any tagged triangulation can be transformed by a
sequence of flips into an ideal triangulation---more precisely, into a
tagged triangulation with no notches.  First, we pick some closed
stratum~$\overline\Omega_\varepsilon$ and use
Corollary~\ref{cor:honest-triangulation} within that stratum to deduce
that any tagged triangulation can be transformed by flips into a
tagged triangulation~$T$ without punctures~$a$ with $\delta_T(a)=0$.
We then transform the punctures~$a$ with $\delta_T(a) = -1$, one by
one, into punctures with $\delta_T(a)=1$ by first placing a puncture
inside a once-punctured digon and then applying (tagged) flips shown
in Figure~\ref{fig:exch-graph-A1A1} (moving from the bottom to the
top). Note that we have excluded all cases in which placing a puncture inside a
digon is not possible, either in Definition~\ref{def:ciliated} or in
the case $\SM$ is a closed surface with one puncture.
\end{proof}

\begin{lemma}
\label{lem:strata-nonempty}
If $\Surf$ is not closed, then each open stratum~$\Omega_\delta$ is
non-empty.  If $\SM$ is closed, then $\Omega_\delta$
is non-empty unless $\delta(a)=0$ for all~$a$, 
in which case $\Omega_\delta$ is empty.
\end{lemma}

\begin{proof}
It is enough to prove the lemma for the open strata contained in the
distinguished ``positive'' closed stratum~$\overline\Omega_{\mathbf{1}}$.
That is, it suffices to show that for any $(0,1)$-vector~$\delta$,
except for the zero vector in the closed case, 
there exists an ideal triangulation with signature~$\delta$. 

If there is some boundary, let $a$ be a marked point on the boundary.
Otherwise, let $a$ be a puncture with $\delta(a) \ne 0$.  In either
case, connect $a$ with each puncture $b$ such that $\delta(b)=0$,
so that these arcs do not intersect pairwise. Wrap a loop based at $a$
around each of these arcs, forming a self-folded
triangle. Triangulate the rest without such
triangles if necessary (see Lemma~\ref{lem:honest-triangulation}).
Each self-folded triangle we cut out reduces the number of punctures~$p$
by one, possibly increases the number of boundary components, and adds
$1$ to the number of marked points on the boundary; this never goes
from a surface allowing triangulations to one with no triangulations.

To see that the stratum with $\delta(a) = 0$ for all punctures~$a$ is
empty in the closed case, note that for a self-folded triangle in an
ideal triangulation, the exterior vertex~$a$ must have $\delta(a) =
\pm 1$.
\end{proof}

\subsection{Mutations and flips in the tagged case}
\label{sec:mutations-tagged} 


We next define adjacency matrices of tagged triangulations,
and verify that, in analogy with Proposition~\ref{pr:mut-tri}, 
flips (i.e., edges in~$\ETSM$) correspond to
mutations of these matrices. 

\begin{definition}[\emph{Signed adjacency matrix of a tagged
      triangulation}] 
\label{def:B-tagged}
Let $T$ be a tagged triangulation of~$\SM$, with its tagged arcs
labeled $1$ through~$n$.
We define the (generalized) signed adjacency matrix $B(T)=(b_{ij})$
as follows. Let $T^\circ$ be the ideal triangulation associated
with~$T$; for each tagged arc
in~$T$, the corresponding ordinary arc in~$T^\circ$ retains the same
label. We then set $B(T)=B(T^\circ)$. 
\end{definition}

Definition~\ref{def:B-tagged} can also be restated in the spirit of 
Remark~\ref{rem:puzzle-pieces}. 

\begin{lemma}
\label{lem:flips-tagged}
Let $T$ and $\overline T$ be labeled tagged triangulations related by a
flip of a tagged arc labeled~$k$.
(The rest of the labeling stays put.) 
Then $B(\overline T)=\mu_k(B(T))$. 
\end{lemma}

\begin{proof}
There is always a closed stratum~$\overline\Omega_\varepsilon$
containing both $T$ and~$\overline T$.  Within that stratum
Proposition~\ref{pr:mut-tri} applies.
\end{proof}

\subsection{Cycles in the graph of flips} 

Our next goal is to prove, using Theorem~\ref{th:flip-cycles}, 
that all cycles in the graph of ``tagged flips'' 
are generated by cycles of length~$4$ and~$5$ (barring the exceptional
cases of closed surfaces with two punctures).  
Our proof is based on a connectivity result
(Lemma~\ref{lem:lozenge-connected} below) 
which in turn relies on the following more general version of the notion
of arc complex.

\begin{definition}
  \label{def:gen-surface}
  Let $\SM$ be a surface as in Definition~\ref{def:ciliated}, without
  the condition that there be at least one marked point on each
  boundary component.  An \emph{arc} in $\SM$ is as in
  Definition~\ref{def:arcs}, and two arcs are \emph{compatible} if
  they do not intersect in the interior of~$\Surf$ 
  (as in Definition~\ref{def:compat-arcs}).  The
  \emph{arc complex}~$\DDSM$ is as in Definition~\ref{def:arc-complex}
  the clique complex for the compatibility relation, while a
  \emph{(generalized) ideal triangulation} of~$\SM$ is a maximal
  collection of disjoint compatible arcs.  The triangulation is
  generalized because if there there are boundary components without marked
  points (such a component is called a \emph{hole}), 
  then some of the complementary regions will not be triangles but
  rather annuli with one marked point on one boundary component.
\end{definition}

\begin{figure}[htbp] 
\begin{center} 
\setlength{\unitlength}{4pt} 
\begin{picture}(40,15)(5,-7) 
\thinlines 
\put(40,0){\makebox(0,0){\usebox{\hdigon}}} 
\put(40,0){\makebox(0,0){\usebox{\highloop}}} 

\put(10,0){\makebox(0,0){\usebox{\hdigon}}} 
\put(10,0){\makebox(0,0){\usebox{\lowloop}}}

\put(15,0){\circle*{1.5}} 
\put(35,0){\circle*{1.5}} 

\thicklines

\put(15,0){\line(1,0){20}}

\end{picture} 
\end{center} 
\caption{A flip inside a digon with a hole} 
\label{fig:digon-flip} 
\end{figure}
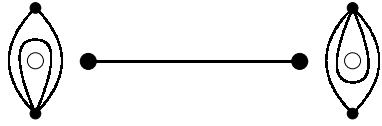

Inside these generalized ideal triangulations, there are two kinds of
(generalized) flips: the flips inside a quadrilateral we saw before, and flips
inside a digon with a hole removed, as in
Figure~\ref{fig:digon-flip}.  The graph~$\EESM$ is defined as before
with this new notion of flips.

\begin{proposition}\cite{hatcher}
\label{pr:digon-connected}
  For $\SM$ as in Definition~\ref{def:gen-surface}, the graph $\EESM$
  is connected.
\end{proposition}

In fact Hatcher proves a much more general result.

\begin{theorem}\cite{hatcher}
\label{th:gen-arc-cplx-is-contractible}
For $\SM$ as in Definition~\ref{def:gen-surface}, the arc complex
$\DDSM$ is contractible unless $\SM$ is an unpunctured polygon or a
polygon with one hole removed; in these exceptional cases $\DDSM$ is a
sphere~\cite{harer-1986}.
\end{theorem}
Since $\DDSM$ is still a
pseudomanifold, this implies the proposition.

\begin{definition}[\emph{$\lozenge$-flips}]
\label{def:lozenge-flips}
A \emph{$\lozenge$-flip} 
is a transformation of an ideal triangulation that
is either an ordinary flip, or a combination of two flips occurring
inside a once-punctured digon (see Figure~\ref{fig:lozenge-flip}). 
This notion is naturally extended to tagged triangulations by
replicating the above definition inside each closed
stratum. 
That is, two tagged triangulations $T_1$ and~$T_2$
are related by a $\lozenge$-flip if, first, they belong to the same closed
stratum~$\overline\Omega_\varepsilon$ and, second, the ideal
triangulations $T_1^\circ$ and~$T_2^\circ$ associated to them are
related by a $\lozenge$-flip. 
\end{definition}

\begin{figure}[htbp] 
\begin{center} 
\setlength{\unitlength}{4pt} 
\begin{picture}(40,10)(5,-5) 

\thinlines

\put(40,0){\makebox(0,0){\usebox{\digon}}} 
\put(40,0){\makebox(0,0){\usebox{\highbar}}} 
\put(40,0){\makebox(0,0){\usebox{\highloop}}} 

\put(10,0){\makebox(0,0){\usebox{\digon}}} 
\put(10,0){\makebox(0,0){\usebox{\lowbar}}} 
\put(10,0){\makebox(0,0){\usebox{\lowloop}}}

\put(15,0){\circle*{1.5}} 
\put(35,0){\circle*{1.5}} 

\thicklines

\put(15,0){\line(1,0){20}}

\end{picture} 
\end{center} 
\caption{A $\lozenge$-flip} 
\label{fig:lozenge-flip} 
\end{figure} 

Comparing Figures~\ref{fig:digon-flip} and~\ref{fig:lozenge-flip},
we see that a $\lozenge$-flip corresponds to a generalized flip 
for a surface in which the puncture inside the digon has been replaced
by a hole. 

\begin{definition}[\emph{$\lozenge$-paths}]
\label{def:lozenge-paths}
A \emph{$\lozenge$-path} is a sequence of (ideal or tagged) triangulations in
which any two consecutive triangulations are related by a
$\lozenge$-flip. 
A path $\Pi$ in $\ETSM$ is a \emph{realization} of a
$\lozenge$-path $\Pi^\diamond$ if $\Pi$ can be obtained from $\Pi^\diamond$ 
by refining each $\lozenge$-flip which is not an ordinary flip 
by a pair of ordinary flips. 
\end{definition}

\begin{lemma}
\label{lem:lozenge-connected}
Assume that $\SM$ has more than one puncture.
Let $b$ be a puncture. 
Then any two ideal triangulations~$T$ with $\delta_T(b)=0$
(i.e., $b$~lies inside a self-folded triangle) are $\lozenge$-connected. 
Consequently, the intersection of each closed stratum of tagged
triangulations with the set $\{T:\delta_{T}(b)=0\}$ is
$\lozenge$-connected. 
\end{lemma}

\begin{proof}
  This follows from Proposition~\ref{pr:digon-connected}.
\end{proof}

\begin{definition}[\emph{R2-cycles}] 
\label{def:R2-cycles}
Let $T$ be a tagged triangulation, and let $B(T)=(b_{ij})$ be the
corresponding (generalized) signed adjacency matrix. 
Consider two arcs in~$T$ labeled $i$ and~$j$ 
such that $b_{ij}=0$ or $b_{ij}=\pm1$. 
It is easy to check that alternately flipping the two arcs labeled $i$
and~$j$ (while retaining the labeling of the remaining arcs) 
brings us back to the original tagged triangulation~$T$ after
$4$~flips (if $b_{ij}=0$) or after $5$~flips (if $b_{ij}=\pm1$). 
The corresponding cycle in the exchange graph $\ETSM$ is called an 
\emph{R2-cycle}. 
\end{definition}

\begin{remark}
In \cite[Section~2.1]{ca2}, the R2-cycles were called 
\emph{geodesic loops} (with respect to the canonical discrete
connection on the exchange graph). We avoid using this terminology here 
because of the misleading connotations it brings as we speak of curves
on surfaces. 

The terminology in Definition~\ref{def:R2-cycles}
reflects the fact that, in effect, we are performing
flips in a (possibly disconnected) bordered surface of rank~$n=2$, 
specifically a pentagon, or a disjoint union of two quadrilaterals, or
a once-punctured digon. 
\end{remark} 

\begin{definition}[\emph{R2-homotopy}] 
Two paths in $\ETSM$ are called \emph{R2-homotopic} to each other
if they can be obtained from each other by a sequence of 
transformations of the following two kinds:
\begin{itemize}
\item
removing two consecutive edges on a path which trace the same edge
in $\ETSM$ in opposite directions; 
\item
replacing a fragment of a path (i.e., several consecutive edges of it)
contained in an R2-cycle by the complement of that fragment 
inside this R2-cycle. 
\end{itemize}
We note that any two realizations of the same $\lozenge$-path 
(see Definition~\ref{def:lozenge-paths}) are R2-homotopic. 
Accordingly, two $\lozenge$-paths are called R2-homotopic if (any of) their 
respective realizations in $\ETSM$  are R2-homotopic to each
other. 
\end{definition}

\begin{theorem}
\label{th:R2-contractible}
Assume that $\SM$ is not a closed surface with exactly two punctures. 
Then any cycle in $\ETSM$ is R2-homotopic to a point. 
\end{theorem}

\begin{proof}
By Theorem~\ref{th:flip-cycles}, the claim holds for any cycle
contained in any closed stratum~$\overline\Omega_\varepsilon$,
since we can identify it with~$\EESM$ via the map
$\tau(\cdot,\varepsilon)$ (see
Definition~\ref{def:sign-strat}). 
To rephrase, any two paths in $\ETSM$ that have common endpoints
and are entirely contained in~$\overline\Omega_\varepsilon$ are
R2-homotopic to each other. 
The idea of proof is to successively use such R2-homotopies to
contract a given cycle~$\Pi$. 

First, we fix a linear order on the set of closed strata.
This will determine the order in which the portions of~$\Pi$
contained in these strata will be deformed. 
Specifically, let $\varepsilon_1,\dots,\varepsilon_{2^p}$ 
be any linear ordering on the $2^p$-element
set of $(\pm1)$-functions on the set of punctures that satisfies 
\[
\sum_p \varepsilon_i(p)< \sum_p \varepsilon_j(p)\Longrightarrow
i<j.
\]
The only purpose of this condition is to ensure that, for any~$i$ 
with $1\leq i<2^p$, we have 
\begin{equation}
\label{eq:strata-ordering}
\overline\Omega_{\varepsilon_i}\cap
\Bigl(\bigcup_{j>i}\overline\Omega_{\varepsilon_j} \Bigr)
=\bigcup_{j>i}\Bigl(\overline\Omega_{\varepsilon_i}
\cap\overline\Omega_{\varepsilon_j} \Bigr) 
=\bigcup_{\substack{j>i\\ \varepsilon_j\gtrdot\varepsilon_i}}
\Bigl(\overline\Omega_{\varepsilon_i}
\cap\overline\Omega_{\varepsilon_j} \Bigr) \,,
\end{equation}
where the notation $\varepsilon_j\gtrdot\varepsilon_i$ means that the
$(\pm1)$-functions $\varepsilon_i$ and $\varepsilon_j$ coincide at all
punctures except for a single puncture~$a$ for which
$\varepsilon_j(a)>\varepsilon_i(a)$.

Starting with the given cycle~$\Pi$,
we are going to successively perform R2-homo\-topies 
of the fragments of a cycle contained inside the closed strata
$\overline\Omega_{\varepsilon_1}$,
$\overline\Omega_{\varepsilon_2}$,\dots 
(in this order). 
For this plan to work, we need to make sure that, 
after each deformation occurring
inside a closed stratum~$\overline\Omega_{\varepsilon_i}$, 
the deformed fragment is contained in the union of closed strata
$\displaystyle\cup_{j>i}\overline\Omega_{\varepsilon_j}\,$.

Let $\Theta$ be a fragment of interest to us; that is, $\Theta$~is a path in
$\ETSM$ contained entirely in~$\overline\Omega_{\varepsilon_i}$ 
which begins and ends at some tagged triangulations $T_1$ and $T_2$ 
located at the ``boundary'' of this closed stratum. 
More precisely, we may assume, taking~\eqref{eq:strata-ordering} into
account, that 
\begin{itemize}
\item 
$T_1\in \overline\Omega_{\varepsilon_j}$ for some $j>i$ such
    that $\varepsilon_j\gtrdot\varepsilon_i$; 
    say, $\varepsilon_j(a)>\varepsilon_i(a)$;
\item 
$T_2\in \overline\Omega_{\varepsilon_k}$ for some $k>i$ such
    that $\varepsilon_k\gtrdot\varepsilon_i$;  
    say, $\varepsilon_k(b)>\varepsilon_i(b)$.
\end{itemize}
In particular, $\delta_{T_1}(a)=\delta_{T_2}(b)=0$. 

Suppose that $j=k$ (so $a=b$).
Then both $T_1$ and $T_2$ lie in 
$\overline\Omega_{\varepsilon_i}\cap
\overline\Omega_{\varepsilon_j}$,
which is precisely the set of tagged triangulations $T$ in
$\overline\Omega_{\varepsilon_i}$ for which $\delta_T(a)=0$. 
By Lemma~\ref{lem:lozenge-connected}, 
there exists a $\lozenge$-path  inside 
$\overline\Omega_{\varepsilon_i}\cap
\overline\Omega_{\varepsilon_j}$ that connects $T_1$ and~$T_2$. 
By Theorem~\ref{th:flip-cycles}, this $\lozenge$-path is R2-homotopic
to~$\Theta$. The latter is therefore R2-homotopic to a path in
$\ETSM$ contained entirely in~$\overline\Omega_{\varepsilon_j}$,  
and we are done. 

Now suppose that $j\neq k$ (so $a\neq b$). 
It follows from Lemma~\ref{lem:strata-nonempty}
(here we use the condition in the theorem) that there exists a tagged
triangulation~$T_0\in \overline\Omega_{\varepsilon_i} \cap \overline\Omega_{\varepsilon_j}\cap
\overline\Omega_{\varepsilon_k}$. 
By Lemma~\ref{lem:lozenge-connected}, there exist $\lozenge$-paths  inside 
$\overline\Omega_{\varepsilon_i}\cap
\overline\Omega_{\varepsilon_j}$ and $\overline\Omega_{\varepsilon_i}\cap
\overline\Omega_{\varepsilon_k}$, respectively,
that connect $T_0$ with $T_1$ and~$T_2$, respectively. 
Again by Theorem~\ref{th:flip-cycles}, $\Theta$~is R2-homotopic to
the concatenation of these paths, which can in turn be realized as
a path in $\ETSM$ contained entirely in
$\overline\Omega_{\varepsilon_j}\cup\overline\Omega_{\varepsilon_k}$. 
The theorem is proved. 
\end{proof}

\begin{remark}
  The above proof can be viewed as a sequence of repeated applications of
  the Mayer-Vietoris Theorem. 
\end{remark}

\begin{remark}
The statement of Theorem~\ref{th:R2-contractible} 
is \emph{false} for closed surfaces with
  $2$ punctures, for the following reasons.
Let $\SM$ be such a surface. 
Then the set of all tagged triangulations is a disjoint union of the
  $8$ open strata determined by the signatures at the two punctures
(recall that the signature $(0,0)$ is impossible): 
\begin{equation}
\label{eq:8-strata}
\begin{array}{ccccc}
\Omega_{-+} & \longleftrightarrow & \Omega_{0+} & \longleftrightarrow & \Omega_{++}\\
\updownarrow &&&& \updownarrow \\
\Omega_{-0} &                     &             &                     & \Omega_{+0}\\
\updownarrow &&&& \updownarrow \\
\Omega_{--} & \longleftrightarrow & \Omega_{0-} & \longleftrightarrow & \Omega_{+-}
\end{array}
\end{equation}
All edges in $\ETSM$ connect tagged triangulations that either
belong to the same open stratum or else to a pair of strata connected
by an arrow in the diagram~\eqref{eq:8-strata}. 
Consider a cycle $\Pi$ that goes through these $8$ strata clockwise. 
Since each R2-cycle is of length $4$ or~$5$, it is impossible to
represent $\Pi$ as a product of R2-cycles. 
\end{remark} 

\subsection{Tropical exchange relations for intersection numbers} 

In this section, we return to the study of intersection
numbers~$(\alpha|\beta)$ introduced in
Section~\ref{sec:intersection-pairing}. 
We will demonstrate that these numbers satisfy the ``tropical
recurrence''~\eqref{eq:exchange-denominator}. 

\pagebreak[3] 

\begin{lemma}
\label{lem:trop-cluster-inters-num}
Let $T=(\beta_1,\dots,\beta_n)$ be a tagged triangulation of~$\SM$, 
with the signed adjacency matrix $B(T)=(b_{ij})$. 
Let $\overline T$ be the tagged
triangulation obtained from~$T$ by a flip replacing 
$\beta_k$ by a tagged arc~$\overline \beta_k$. 
Then, for any tagged arc~$\alpha$ in~$\SM$, we have
\begin{equation}
\label{eq:exchange-inters-num}
(\alpha|\overline \beta_k) =
- (\alpha|\beta_k) +
\max\Bigl(\displaystyle\sum_{i=1}^n \,[b_{ik}]_+\, (\alpha|\beta_i),
\sum_{i=1}^n\, [-b_{ik}]_+ \,(\alpha|\beta_i) \Bigr) \,, 
\end{equation}
\end{lemma}

\begin{proof}
Without loss of generality, we may assume that the tagged
triangulations~$T$ and~$\overline T$ belong to the positive closed stratum
(see Definition~\ref{def:sign-strat}). 
The proof then proceeds case by case, for each combinatorially
different type of a flip $T\to\overline T$. 
(Note that the roles of $T$ and $\overline T$ can be switched without
affecting the statement~\eqref{eq:exchange-inters-num} to be proved.) 
These types can be catalogued using the puzzle-piece decomposition
of (the ideal triangulation corresponding to)~$T$, 
similarly to the proof of Proposition~\ref{pr:mut-tri}.

Two most basic types of flips are (cf.\ Figure~\ref{fig:flips-F1-F2}): 
\begin{itemize}
\item[{(F1)}]
flipping a diagonal of a
quadrilateral~$Q$ formed by two triangles (none of them self-folded)
sharing precisely one edge;
\item[{(F2)}]
flipping a tagged arc inside a once-punctured digon. 
\end{itemize}
We will prove~\eqref{eq:exchange-inters-num} for the flips of type~F1,
sketch the proof for type~F2, and  
omit the tedious verification of the remaining cases.

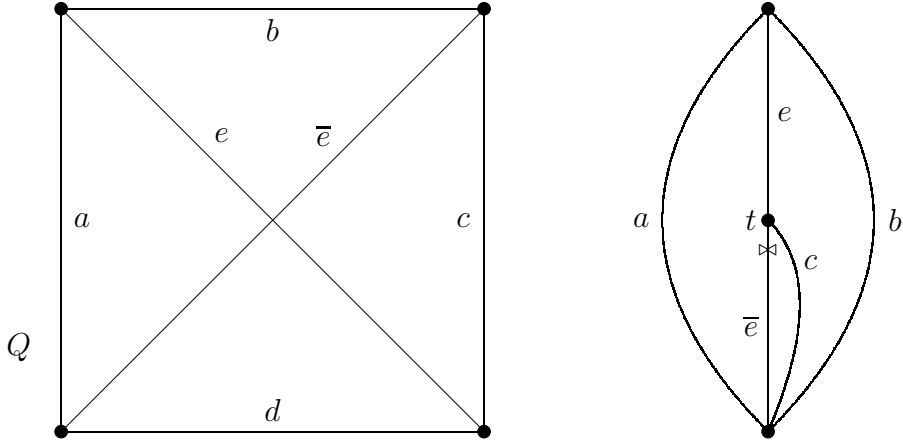
\begin{figure}[htbp] 
\begin{center} 
\setlength{\unitlength}{16pt} 
\begin{picture}(10,10)(0,0) 

\thinlines 

\put(0,0){\circle*{.3}} 
\put(0,10){\circle*{.3}} 
\put(10,0){\circle*{.3}} 
\put(10,10){\circle*{.3}} 
\put(0,0){\line(0,1){10}}
\put(10,0){\line(0,1){10}}
\put(0,0){\line(1,0){10}}
\put(0,10){\line(1,0){10}}
\put(0,0){\line(1,1){10}}
\put(0,10){\line(1,-1){10}}

\put(0.5,5){\makebox(0,0){$a$}} 
\put(5,9.5){\makebox(0,0){$b$}} 
\put(9.5,5){\makebox(0,0){$c$}} 
\put(5,0.5){\makebox(0,0){$d$}} 
\put(3.8,7){\makebox(0,0){$e$}} 
\put(6.2,7){\makebox(0,0){$\overline e$}} 

\put(-1,2){\makebox(0,0){$Q$}} 
\end{picture} 
\qquad
\begin{picture}(8,10)(0,0) 

\thinlines 

\qbezier(5,0)(0,5)(5,10)
\qbezier(5,0)(10,5)(5,10)
\qbezier(5,0)(6.5,3.5)(5,5)
\put(5,0){\circle*{.3}} 
\put(5,5){\circle*{.3}} 
\put(5,10){\circle*{.3}} 
\put(5,0){\line(0,1){5}}
\put(5,5){\line(0,1){5}}
\put(4.6,5){\makebox(0,0){$t$}} 
\put(2,5){\makebox(0,0){$a$}} 
\put(8,5){\makebox(0,0){$b$}} 
\put(6,4){\makebox(0,0){$c$}} 
\put(4.6,2.5){\makebox(0,0){$\overline e$}} 
\put(5.4,7.5){\makebox(0,0){$e$}} 

\put(5,4.3){\makebox(0,0){$\notch$}} 
\end{picture} 
\end{center} 
\caption{Flips of types F1 and F2} 
\label{fig:flips-F1-F2} 
\end{figure} 

The only arcs $\beta\in T$ contributing
to~\eqref{eq:exchange-inters-num} are $\beta_k$, $\overline \beta_k$,
and those $\beta_i$ for which $b_{ik}\neq 0$.  
In the case of a flip of type~F1, these arcs are precisely 
the sides and the diagonals of the quadrilateral~$Q$. 
Let $a,\dots,e,\overline e$ denote the corresponding intersection
numbers, as shown in Figure~\ref{fig:flips-F1-F2} on the left. 
Then \eqref{eq:exchange-inters-num} takes the form 
\begin{equation}
\label{eq:e+e=max}
e+e'=\max(a+c,b+d)\,. 
\end{equation}
To prove this, we first get rid of the tagging. 
If $\alpha$ is notched at some vertex of~$Q$, 
then removing the notch would subtract~$1$ from each of the sums
$a+c$, $b+d$, and $e+\overline e$. 
We may therefore assume that $\alpha$ is a plain arc. 
We may also assume that $\alpha$ is none of the $6$ arcs shown in
Figure~\ref{fig:flips-F1-F2}, since otherwise \eqref{eq:e+e=max} can be checked
directly. 
Thus, we need not worry about the summands $C$ and $D$ in the
definition of intersection numbers, and only take careful account of the
summands $A$ and~$B$. 

First let us assume that $\alpha$ is not a loop, so that $B=0$ throughout.   
Each segment $\sigma\subset\alpha$ cutting across~$Q$ makes contributions to
the 
the three sums $a+c$, $b+d$, and $e+\overline e$;
we denote these contributions by $q_{ac}$, $q_{bd}$,
and~$q_{e\overline e}$, respectively. 
Then the following cases may occur: 
\begin{itemize}
\item[{(F1-1})]
$\sigma$ cuts through two consecutive sides of~$Q$, 
with $q_{ac}=q_{bd}=q_{e\overline e}=1$; 
\item[{(F1-2})]
$\sigma$ cuts through two opposite sides of~$Q$,  
with $q_{e\overline e}=2$, $q_{ac}=2$, $q_{bd}=0$; 
\item[{(F1-3})]
$\sigma$ cuts through two opposite sides of~$Q$,  
with $q_{e\overline e}=2$, $q_{ac}=0$, $q_{bd}=2$; 
\item[{(F1-4})]
$\sigma$ starts at a vertex of~$Q$, 
with $q_{e\overline e}=1$, $q_{ac}=1$, $q_{bd}=0$; 
\item[{(F1-5})]
$\sigma$ starts at a vertex of~$Q$, 
with $q_{e\overline e}=1$, $q_{ac}=0$, $q_{bd}=1$.  
\end{itemize}
It is easy to check 
(keeping in mind that different segments $\sigma$ do not
intersect each other, and $\alpha$ is not a loop) 
that we are left with the following possibilities:
\begin{itemize}
\item 
all $\sigma$'s are of types F1-1, F1-2, or F1-4, so that 
$q_{e\overline e}=q_{ac}\geq q_{bd}$ in each case; 
\item 
all $\sigma$'s are of types F1-1, F1-3, or F1-5, so that 
$q_{e\overline e}=q_{bd}\geq q_{ac}$ in each case. 
\end{itemize}
We then conclude by additivity that in these cases, 
$e+\overline e=a+c\geq b+d$ or $e+\overline e=b+d\geq a+c$,
respectively, yielding~\eqref{eq:e+e=max}. 

If $\alpha$ is a loop, we need to make some adjustments to the above
proof. Namely, we need to account for the summands~$B$ in
Definition~\ref{def:inters-pairing}, and also allow for the 
additional possibility that 
\begin{itemize}
\item 
all $\sigma$'s but two are of type F1-1, one is of type F1-4, 
and one of type F1-5. 
\end{itemize}
In the latter case, the total $B$-contribution to~$e+\overline e$
is~$-1$ whereas the $B$-contributions to $a+c$ and $b+d$ are
both~$0$. Combined with the $A$-contributions of the two segments of
types F1-4 and F1-5, this results in equal quantities 
(namely~$-1$) for all three sums; so in this case we actually get 
$e+\overline e=a+c=b+d$. 
In all other cases (that is, unless F1-4 and F1-5 occur
simultaneously), the $B$-contributions associated with pairs of
segments of $\alpha$ should be considered together
with the corresponding $A$-contributions; 
each time, the same argument as above works, i.e., all inequalities
will go the same way, consistently favouring $a+c$ over $b+d$
(or vice versa).

We note that essentially the same proof works in the case where two
opposite sides of~$Q$ (or both pairs of opposite sides)
are identified; keep in mind that in that case, the corresponding
entries of the matrix~$B(T)$ are doubled, making up for the gluing of
intersection points. 

For a flip of type~F2, the formula \eqref{eq:exchange-inters-num} takes the form 
\begin{equation}
\label{eq:e+e=max(a,b)}
e+\overline e=\max(a,b)\,. 
\end{equation}
The argument is similar to the one given above in the F1 case,
with the additional complication due to the presence of a puncture~$t$ 
inside the domain in which the flip occurs---in this case, inside a digon. 
As a result, we should consider the possibility that $\alpha$ may have 
$t$ as an endpoint, or even be a loop based at~$t$. 
As above, the proof examines several possible ways in which segments of
$\alpha$ may cut through the digon.
Each time the contributions to $e+\overline e$ turn out to be equal to 
the maximum of their counterparts for $a$ and~$b$. 
As before, the maximum is determined consistently for all segments
(or pairs of segments, in the case of $B$-contributions, which are
pooled together with the $A$-contributions from the same segments);
that is, either all contributions to $a$ are larger than respective
contributions to~$b$, or the opposite happens. 
Either way, \eqref{eq:e+e=max(a,b)} follows. 

Other types of flips are examined in a similar fashion.
One of them is illustrated in Figure~\ref{fig:flip-F3}.
This flip (dubbed~F3) occurs at an arc shared by a
triangle and a punctured digon. 
The formula \eqref{eq:exchange-inters-num} takes the form 
\begin{equation}
\label{eq:e+e,type-F3}
e+\overline e=\max(a+c,b+d+f) \,. 
\end{equation}

\begin{figure}[htbp] 
\begin{center} 
\setlength{\unitlength}{8pt} 
\begin{picture}(20,17)(0,0) 
\thicklines 
\put(0,0){\line(1,0){20}}
\put(0,0){\line(100,173){10}}
\put(10,17.3){\line(100,-173){10}}
 
\put(3.5,8){\makebox(0,0){$c$}}
\put(16.5,8){\makebox(0,0){$a$}}
\put(10,-1){\makebox(0,0){$b$}}
\put(6,5){\makebox(0,0){$e'$}}
\put(14,5){\makebox(0,0){$e$}}
\put(8.7,10){\makebox(0,0){$d$}}
\put(11.3,10){\makebox(0,0){$f$}}

\put(0,0){\circle*{0.5}} 
\put(20,0){\circle*{0.5}} 
\put(10,17.3){\circle*{0.5}} 
\put(10,5.8){\circle*{0.5}} 

\qbezier(20,0)(0,0)(10,17.3)
\qbezier(0,0)(20,0)(10,17.3)

\qbezier(10,17.3)(8.5,11.7)(10,5.8)
\qbezier(10,17.3)(11.5,11.7)(10,5.8)

\put(10.43,7.5){\makebox(0,0){$\notch$}}

\end{picture} 
\end{center} 
\caption{Flip of type F3} 
\label{fig:flip-F3} 
\end{figure}
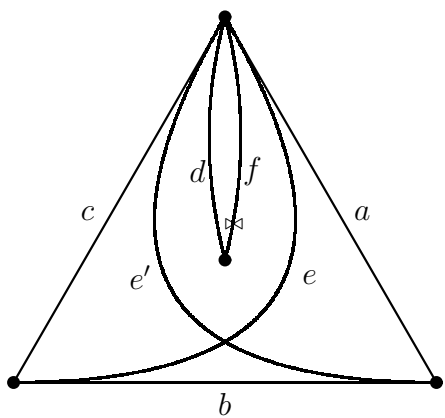 

We omit the consideration of the remaining types of flips,
which include flipping an arc shared by two
once-punctured digons, and flips involving a twice-punctured monogon. 
\end{proof} 

\subsection{Proof of Theorems~\ref{th:cluster-complex-top},
  \ref{th:cluster-cpx=arc-cpx}, and~\ref{th:denoms-as-inters-numbers}}
\label{sec:main-proofs}

The proof is based on the following result.

\begin{proposition}[{\cite[Proposition~2.3, Lemma~2.4]{ca2}}]
\label{pr:covering-exch-graph}
Let $\Delta$ be a simplicial complex on a (possibly infinite) ground
set~$\Psi$ satisfying the following conditions: 
\begin{itemize}
\item
$\Delta$ is an $(n-1)$-dimensional pseudomanifold with a connected dual graph~$\Exch$; 
\item
the dual graph of the link of every non-maximal simplex in $\Delta$ is
connected;
\item
the fundamental group of $\Exch$ is generated by cycles 
(pinned down to a fixed basepoint) that arise as links of
codimension~$2$ simplices.  
\end{itemize}
Suppose that $\{B(C)=(b_{\beta\gamma})_{\beta,\gamma\in C}\}$ is 
a collection of skew-symmetrizable matrices labeled by the maximal
simplices of~$\Delta$ such that: 
\begin{itemize}
\item
for each edge $(C,\overline C)$
in~$\Exch$ with $\overline C=C-\{\gamma\}\cup\{\overline\gamma\}$, 
the matrix $B(\overline C)$ is obtained from $B(C)$ by a matrix 
mutation in direction (of the label of)~$\gamma$ 
(we assume that the labels for the elements of $C\cap\overline C$
are the same in $C$ and in~$\overline C$);
\item
for each codimension~$2$ simplex $D\in\Delta$ whose link is a finite
cycle (say of length~$\ell$), 
and for any maximal simplex $C=D\cup\{\beta,\gamma\}$,
the number $|b_{\beta\gamma}b_{\gamma\beta}|$ is equal to $0$, $1$,
$2$, or~$3$, and furthermore $\ell$ is equal to $4$, $5$, $6$, or~$8$,
respectively. 
\end{itemize}
Finally, let $\Acal$ be a normalized cluster algebra with exchange matrix
$B=B(C_0)$ for some $C_0\in\Exch$. 
Then there is a (canonically defined) 
surjection $C\mapsto \Sigma(C)$ from $\Exch$ to the set of
seeds in~$\Acal$ and a (canonical) surjection $\beta\mapsto x[\beta]$ from $\Psi$
to the set of cluster variables of~$\Acal$ such that 
\begin{itemize}
\item
for each edge $(C,\overline C)$
in~$\Exch$ with $\overline C=C-\{\gamma\}\cup\{\overline\gamma\}$, 
the seeds $\Sigma(C)$ and $\Sigma(\overline C)$ are related by a
mutation in direction~$\gamma$; 
\item
the cluster variables in each seed $\Sigma(C)$ are precisely $\{x[\beta]:\beta\in
C\}$. 
\end{itemize}
\end{proposition}

We are going to use Proposition~\ref{pr:covering-exch-graph}
for $\Delta=\DTSM$, so that $\Psi$ is the set of tagged arcs
in~$\SM$, and $\Exch=\ETSM$ is the graph of tagged triangulations. 
The restrictions on~$\Delta$ stated in
Proposition~\ref{pr:covering-exch-graph} are satisfied by 
Theorem~\ref{th:tagged-pseudo},
Proposition~\ref{pr:dual-tagged-connected},
and Theorem~\ref{th:R2-contractible}, as long as $\SM$ is not a
surface with exactly two punctures.
The connectedness of links follows by noting that each link is in fact
a tagged arc complex of another (possibly disconnected) surface. 

The matrices $B(C)$ associated with the maximal simplices of~$\Delta$ 
are the signed adjacency matrices of tagged triangulations. 
The restrictions on them stated in
Proposition~\ref{pr:covering-exch-graph} are satisfied by
Lemma~\ref{lem:flips-tagged} and by direct examination of all possible
R2-cycles (which are links of codimension~$2$ simplices). 
The only values of $|b_{\beta\gamma}b_{\gamma\beta}|$ that actually
come up in this context are $0$ and~$1$, 
corresponding to $4$-cycles and $5$-cycles in~$\ETSM$. 

The conclusion of Proposition~\ref{pr:covering-exch-graph} then tells
us that in a cluster algebra $\Acal$ associated with~$\SM$, 
there is a well defined cluster variable $x[\beta]$ corresponding to
each tagged arc~$\beta$, and a well defined seed $\Sigma(T)$ corresponding to
each tagged triangulation~$T$;
the cluster variables participating in $\Sigma(T)$ are precisely the
$x[\beta]$'s for $\beta\in T$. 
In short, $\ETSM$ and $\DTSM$ \emph{cover} the exchange graph
and cluster complex of~$\Acal$. 

To complete the proof of the main claims in
Theorems~\ref{th:cluster-complex-top} and~\ref{th:cluster-cpx=arc-cpx}, 
we need to show that this covering is in fact an isomorphism,
that is, the aforementioned surjections are actually bijections. 
In practical terms, we need to show that all cluster variables
$x[\beta]$ are distinct. 

But first, let us prove Theorem~\ref{th:denoms-as-inters-numbers}.
Now that there is a well defined cluster variable $x[\beta]$ for each
tagged arc~$\beta$, the statement of
Theorem~\ref{th:denoms-as-inters-numbers}
follows directly from Lemma~\ref{lem:denom-recurrence} and 
Lemma~\ref{lem:trop-cluster-inters-num}. 

The fact that all cluster variables $x[\beta]$ are distinct is then
immediate from Theorem~\ref{th:denoms-as-inters-numbers}, 
as the latter implies that for any
$\gamma\neq\beta$, the expression for $x[\gamma]$ in terms of a seed
containing~$x[\beta]$ may not be equal to~$x[\beta]$
since $(\gamma|\beta)>-1=(\beta|\beta)$. 


It remains to verify that the seeds containing a given cluster
variable form a connected subgraph of~$\ETSM$. 
This statement is a special case of the ``connectedness of links''
property discussed earlier in the proof. 
\qed 

\begin{remark}
The assertion that all cluster variables $x[\beta]$ are distinct 
can be given a different proof that avoids the use of denominators and
intersection numbers but instead relies on hyperbolic geometry.
This argument will be given in~\cite{fst-hyper}. 
\end{remark}

\section{On the topology of cluster complexes}
\label{sec:cluster-cplx-conj}

In this section we determine (up to homotopy) the topology of 
the cluster complex associated with a general bordered surface 
with marked points---or,
equivalently (see Theorem~\ref{th:cluster-cpx=arc-cpx}), 
the topology of the corresponding tagged arc complex~$\DTSM$. 

As before, $\SM$ is a bordered surface with marked points
(as described in Definition~\ref{def:ciliated}),
and $\DTSM$ is its tagged arc complex
(see Definition~\ref{def:tagged-arc-complex}). 

\begin{theorem}
\label{th:cluster-complex-contractible-or-sphere}
The cluster complex of a cluster algebra $\Acal$
associated with $\SM$ is either contractible or 
homotopy equivalent to a sphere provided $\SM$ is not a closed surface with
exactly two punctures.
Specifically:
\begin{itemize}
\item
If $\Acal$ is of finite type, 
the cluster complex  is homeomorphic to an $(n-1)$-dimensional
sphere~$S^{n-1}$.  
\item
If $\SM$ is a closed surface with $p\ge 2$ punctures, 
the cluster complex  is homotopy equivalent to~$S^{p-1}$. 
\item
In all other cases, the cluster complex  is contractible. 
\end{itemize}
\end{theorem}

Theorem~\ref{th:cluster-complex-contractible-or-sphere} is a
consequence of the following version for tagged arc complexes.

\begin{theorem}
\label{th:arc-complex-conractible-or-sphere}
The tagged arc complex $\DTSM$ is either contractible or 
homotopy equivalent to a sphere. 
Specifically:
\begin{itemize}
\item
If 
$\SM$ is a polygon
  or a once-punctured polygon, 
then $\DTSM$ is homeomorphic to~$S^{n-1}$.  
\item
If $\SM$ is a closed surface with $p$ punctures, 
then $\DTSM$ is homotopy equivalent to~$S^{p-1}$. 
\item
In all other cases, $\DTSM$ is contractible. 
\end{itemize}
\end{theorem}

The proof of Theorem~\ref{th:arc-complex-conractible-or-sphere} will
use the following classical tool from combinatorial topology
(see, e.g., \cite[Theorem~10.7]{bjorner-encyclop} and references therein). 

\begin{definition}
The \emph{nerve} $\mathcal{N}(X)$ 
of a collection of sets $X=(X_i)_{i\in I}$ 
is the simplicial complex whose ground set is~$I$ and whose simplices
are finite subsets $\{i_1,\dots,i_k\}\subset I$ such that 
$X_{i_1}\cap\cdots\cap X_{i_k}\neq\emptyset$. 
\end{definition}

\begin{lemma}[The Nerve Lemma]
\label{lem:nerve-lemma}
Let $\Delta$ be a triangulable topological space
and $X=(X_i)_{i\in I}$ a finite family of closed
subsets of~$\Delta$ such that $\Delta=\bigcup_{i\in I} X_i\,$.
If each intersection $X_{i_1}\cap\cdots\cap X_{i_k}$ is either empty
or contractible, then $\Delta$ is homotopy equivalent to the
nerve~$\mathcal{N}(X)$. 
\end{lemma}

\begin{proof}[Proof of Theorem~\ref{th:arc-complex-conractible-or-sphere}]
  If $\SM$ has no punctures and is not a polygon, then the tagged arc complex $\DTSM$
  coincides with the ordinary untagged arc complex $\DDSM$, which is
  contractible by Theorem~\ref{th:arc-cplx-is-contractible}.  If
  $\SM$ is an unpunctured or a once-punctured polygon, 
  then $\DTSM$ is the cluster complex of
  a cluster algebra of finite type $A_n$ or~$D_n$, 
  which is homeomorphic to~$S^{n-1}$
  \cite[Theorem~1.13]{ca2}\cite[Corollary~1.11]{yga}.

  From now on assume that $p \ge 1$ and $\SM$ is not a once-punctured polygon.
  Consider the map $\varphi : \DTSM
  \rightarrow \RR^p$ which is linear on simplices and defined on
  vertices of $\DTSM$ (tagged arcs $\gamma$) by
  \begin{equation}
    \varphi(\gamma; i) =
      \begin{cases}
        0 & \text{if $\gamma$ does not meet puncture~$i$;}\\
        1 & \text{if $\gamma$ is tagged plain at~$i$;}\\
        -1 & \text{if $\gamma$ is notched at~$i$.}\\
      \end{cases}
  \end{equation}
  For a $(\pm1)$-valued function $\varepsilon$ on the set of punctures
  in~$\SM$, define 
\[
X_\varepsilon = X_\varepsilon\SM
\stackrel{\rm def}{=}
\{x \in
  \DTSM : \varphi(x; i) \,\varepsilon(i) \ge 0 \text{\ for all
    punctures $i$}\}.
\]
  The sets $X_\varepsilon$ are closed and cover $\DTSM$.  Similarly, for
  $\delta$ a $\{0, \pm 1\}$-valued function on the set of
  punctures, define 
\[
X_\delta
\stackrel{\rm def}{=}
\{ x \in \DTSM :
\text{$\delta(i) = 0\,\Rightarrow\, \varphi(x; i) = 0$
  and $\varphi(x;i)\, \delta(i) \ge 0$ for all $i$}\}.  
\]
The intersection of
  any collection of $X_\varepsilon$'s is an~$X_\delta$.

  For $\delta$ as above, let $\SM_\delta$ be the surface with marked
  points obtained from $\SM$ by replacing each puncture~$i$ with $\delta(i) = 0$
  by a \emph{hole}, a boundary component with no marked points
(cf.\ Definition~\ref{def:gen-surface}).

  \begin{lemma}
    \label{lem:x-delta-is-arc-complex}
    For $\delta$ a $\{0,\pm1\}$-valued function, $X_\delta$ is
    homeomorphic to the arc complex~$\DDSM_\delta$.
  \end{lemma}

  \begin{proof}
    By symmetry, we may assume that $\delta$ only has values of $0$
    and $+1$.  
    Define the map $\tilde\tau:\DDSM\to \DTSM$, 
    a variant of the map~$\tau$ of
    Definition~\ref{def:plain-arc-as-tagged-arc}, as follows.
    If $\gamma$ does not cut out
    a once-punctured monogon, then $\tilde\tau(\gamma)$ is $\gamma$ with both
    ends tagged plain, as before.  Otherwise, let $\beta$ be the plain
    arc enclosed by~$\gamma$, and let $\beta^{\notch}$ be
    the same arc with the encircled end notched and the other
    end plain; then $\tilde\tau(\gamma)=(\beta+\beta^{\notch})/2$.
    Extend $\tilde\tau$ by linearity over simplices in $\DDSM$.  Note
    that the vertices of any simplex in $\DDSM$ map inside a single
    simplex of $\DTSM$, so $\tilde\tau$ is well-defined.

    The map~$\tilde\tau$ maps $\DDSM_\delta\subset \DDSM$ inside
    $X_\delta \subset \DTSM$, since for every vertex $\gamma$ of
    $\DDSM_\delta$ and every puncture~$i$ with $\delta(i)=0$, 
    we have $\varphi(\tilde\tau(\gamma);i)=0$.  Conversely, we can define an
    inverse map from $X_\delta$ to $\DDSM_\delta$: the arcs
    incident to a puncture $i$ with $\delta(i)=0$ must come in equal weights
    plain and notched, and we map such a mixture to a loop enclosing
    the puncture.
  \end{proof}

   \begin{lemma}
    \label{lem:x-delta-contractible}
    $X_\delta$ is contractible if it is non-empty.  $X_\delta$ is
    empty only if $\Surf$ is a closed surface and $\delta$ is zero
    everywhere.
  \end{lemma}

  \begin{proof}
    The first statement follows from
    Lemma~\ref{lem:x-delta-is-arc-complex} and
    Theorem~\ref{th:gen-arc-cplx-is-contractible}, since we have
    already ruled out the non-contractible cases in that theorem.  The
    second statement is analogous to Lemma~\ref{lem:strata-nonempty}.
  \end{proof}

  To complete the proof of
  Theorem~\ref{th:arc-complex-conractible-or-sphere}, consider the
  covering of $\DTSM$ by the sets~$X_\varepsilon$.  By
  Lemma~\ref{lem:x-delta-contractible}, this covering satisfies
  the hypotheses of Lemma~\ref{lem:nerve-lemma}.  If~$\SM$ is not a
  closed surface with punctures then every intersection of a subset of the
  $X_\varepsilon$'s is non-empty.  
  Thus $\mathcal N(X)$ and $\DTSM$ are contractible.  Otherwise, we
  use the following lemma.

\pagebreak[3]

\begin{lemma}
\label{lem:cube-nerve}
Let $\mathcal{N}_p$ be  the simplicial complex defined as follows:
\begin{itemize}
\item
The ground set of $\mathcal{N}_p$ is the 
set $\{1,-1\}^p$ of all
$(\pm1)$-vectors $\varepsilon=(\varepsilon_1,\dots,\varepsilon_p)$. 
\item
A subset $J\subset \{1,-1\}^p$ is a simplex in $\mathcal{N}_p$ 
if and only if there exist an index $k\in\{1,\dots,p\}$ and a sign 
$\sigma\in\{1,-1\}$ such that $\varepsilon_k\cdot\sigma\ge 0$ for all 
$\varepsilon\in J$. 
\end{itemize}
Then $\mathcal{N}_p$ is homotopy equivalent to the $(p-1)$-dimensional
sphere~$S^{p-1}$. 
\end{lemma}

\begin{proof}
Let $\Delta$ be the boundary of the cube
$[-1,1]^p\subset\mathbb{R}^p$. 
For $\varepsilon=(\varepsilon_1,\dots,\varepsilon_p)\in\{1,-1\}^p$, 
let 
\[
X_\varepsilon=\Delta\cap\{x=(x_1,\dots,x_p): x_k\,\varepsilon_k\ge 0
\text{\ for all $k$}\}. 
\]
Then $\mathcal{N}(X)=\mathcal{N}_p$ and $\Delta\cong S^{p-1}$, so 
applying Lemma~\ref{lem:nerve-lemma} establishes the claim. 
\end{proof}

  If $\SM$ is a closed surface with $p$~punctures, then in view of
  Lemma~\ref{lem:x-delta-contractible} the nerve of the family
  $(X_\varepsilon\SM)$ coincides with $\mathcal N_p$. Therefore by
  Lemma~\ref{lem:cube-nerve},  $\mathcal N(X)$ and thus $\DTSM$ are homotopy
  equivalent to the sphere~$S^{p-1}$.
\end{proof}

\begin{remark}
The mapping class group of $\SM$ naturally acts on the cluster
complex, or the tagged arc complex $\DTSM$. 
The quotient is a finite simplicial complex that generalizes the
``arc complexes'' studied by R.~Penner in~\cite{penner-arc-complexes}.
\end{remark}

\section{Polynomial vs.~exponential growth} 
\label{sec:poly-growth}

A cluster algebra (or the corresponding exchange graph) 
has \emph{polynomial growth} if the number of
distinct seeds which can be obtained from
a fixed initial seed by at most $n$ mutations is bounded from above by
a polynomial function of~$n$.  A cluster algebra has \emph{exponential growth}
if the number of such seeds is bounded from below by an exponentially
growing function of~$n$.  (There is always an exponential upper bound.)



\begin{figure}
\setlength{\unitlength}{1.5pt} 
\begin{picture}(220,35)(-10,-15) 
\put(-15,-2){\makebox(0,0){$\Gamma(n_1,n_2)$}}
\put(20,-6){\makebox(0,0){$a_1$}}
\put(40,-6){\makebox(0,0){$a_2$}}
\put(60,-6){\makebox(0,0){$\cdots$}}
\put(80,-6){\makebox(0,0){$a_{n_1\!-\!1}$}}
\put(104,14){\makebox(0,0){$a_{n_1}$}}
\put(104,-16){\makebox(0,0){$b_{n_2\!+\!1}$}}

\put(226,-10){\makebox(0,0){$b_0$}}
\put(226,10){\makebox(0,0){$b'_0$}}
\put(200,-6){\makebox(0,0){$b_1$}}
\put(180,-6){\makebox(0,0){$b_2$}}
\put(155,-6){\makebox(0,0){$\cdots$}}
\put(125,-6){\makebox(0,0){$b_{n_2}$}}

\put(20,0){\line(1,0){60}} 
\put(80,0){\line(2,-1){20}} 
\put(80,0){\line(2,1){20}} 
\put(120,0){\line(-2,-1){20}} 
\put(120,0){\line(-2,1){20}} 
\multiput(20,0)(20,0){4}{\circle*{2}} 
\put(220,10){\circle*{2}} 
\put(220,-10){\circle*{2}} 
\put(100,10){\circle*{2}} 
\put(100,-10){\circle*{2}} 
\multiput(120,0)(20,0){5}{\circle*{2}} 
\put(120,0){\line(1,0){80}} 
\put(200,0){\line(2,-1){20}} 
\put(200,0){\line(2,1){20}} 
\put(99.3,-10){\line(0,1){20}} 
\put(100.7,-10){\line(0,1){20}} 
\end{picture} 

\setlength{\unitlength}{1.5pt} 
\begin{picture}(220,48)(-10,-18) 
\put(-15,-2){\makebox(0,0){$\Gamma(n_1,n_2,n_3)$}}
\put(20,-6){\makebox(0,0){$a_1$}}
\put(40,-6){\makebox(0,0){$\cdots$}}
\put(60,-6){\makebox(0,0){$a_{n_1\!-\!1}$}}
\put(84,14){\makebox(0,0){$a_{n_1}$}}
\put(84,-16){\makebox(0,0){$b_{n_2\!+\!2}$}}

\put(140,-6){\makebox(0,0){$b_1$}}
\put(121,-6){\makebox(0,0){$\cdots$}}
\put(105,-6){\makebox(0,0){$b_{n_2\!+\!1}$}}
\put(164,14){\makebox(0,0){$c_{n_3}$}}
\put(164,-16){\makebox(0,0){$b_0$}}

\put(20,0){\line(1,0){40}} 
\put(60,0){\line(2,-1){20}} 
\put(60,0){\line(2,1){20}} 
\put(100,0){\line(-2,-1){20}} 
\put(100,0){\line(-2,1){20}} 
\multiput(20,0)(20,0){3}{\circle*{2}} 
\put(160,10){\circle*{2}} 
\put(160,-10){\circle*{2}} 
\put(80,10){\circle*{2}} 
\put(80,-10){\circle*{2}} 
\multiput(100,0)(20,0){3}{\circle*{2}} 
\put(100,0){\line(1,0){40}} 
\put(140,0){\line(2,-1){20}} 
\put(140,0){\line(2,1){20}} 
\put(79.3,-10){\line(0,1){20}} 
\put(80.7,-10){\line(0,1){20}} 

\put(159.3,-10){\line(0,1){20}} 
\put(160.7,-10){\line(0,1){20}} 
\put(180,0){\line(-2,-1){20}} 
\put(180,0){\line(-2,1){20}} 
\multiput(180,0)(20,0){3}{\circle*{2}} 
\put(180,0){\line(1,0){40}} 
\put(185,-6){\makebox(0,0){$c_{n_3\!-\!1}$}}
\put(220,-6){\makebox(0,0){$c_1$}}
\put(201,-6){\makebox(0,0){$\cdots$}}
\end{picture} 

\caption{Diagrams for the cluster algebras of quadratic and cubic growth.
All triangles are oriented. Orientations of the remaining edges are of no importance.} 
\label{fig:poly-growth} 
\end{figure}
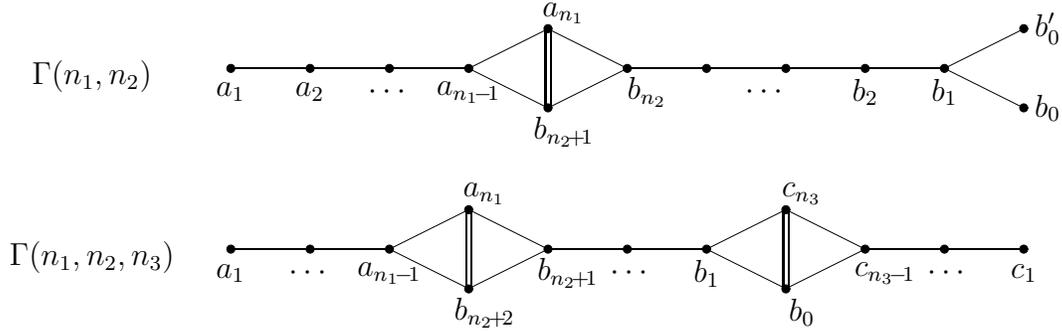

\begin{proposition} 
\label{pr:poly-growth-cluster}
Let $\Acal$ be a cluster algebra associated with 
a connected bordered surface with marked points.
Then $\Acal$ has polynomial growth if and only if 
it is defined by a diagram/quiver in the following list:
\begin{itemize}
\item finite type $A_n$ (finite);
\item finite type $D_n$ (finite);
\item affine type $\tilde A(n_1, n_2)$ (linear growth);
\item affine type $\widetilde D_{n-1}$ (linear growth);
\item diagram $\Gamma(n_1,n_2)$ ($n_1,n_2\in\ZZ_{>0}$) shown in Figure~\ref{fig:poly-growth}  (quadratic
  growth); 
\item diagram $\Gamma(n_1,n_2,n_3)$ ($n_1,n_2,n_3\in\ZZ_{>0}$) shown in Figure~\ref{fig:poly-growth}  (cubic growth).
\end{itemize}
Otherwise $\Acal$ has exponential growth.
\end{proposition}

For the class of cluster algebras considered in this paper, the growth
rate coincides with the growth rate of the graph~$\ETSM$ of tagged
flips (by Theorem~\ref{th:cluster-cpx=arc-cpx}).  
The counterpart of Proposition~\ref{pr:poly-growth-cluster} formulated in
terms of the surface is given below. 

\begin{proposition}
\label{pr:poly-growth-surface}
Let $\SM$ be a connected bordered surface with marked points,
which has $b$ boundary components and $p$ punctures. 
Then the graph of tagged flips $\ETSM$ has polynomial growth if and
only if $\Surf$ is a sphere, and $b+p\leq 3$. 
In all other cases, $\ETSM$ has exponential growth.
\end{proposition}

Proposition~\ref{pr:poly-growth-surface} implies
Proposition~\ref{pr:poly-growth-cluster}: the possibilities for the
numbers of punctures and boundary components are
\begin{itemize}
\item one boundary component (a disk) with $n+3$ marked points: type~$A_n$;
\item one boundary component with $n$ marked points, one puncture (a
  once-punc\-tured disk): type~$D_n$;
\item two boundary components (an annulus) with $n_1$ and $n_2$ marked points,
    respectively: type $\tilde A(n_1, n_2)$;
\item one boundary component with $n-3$ marked points and two punctures (a twice-punctured disk):
  type~$\widetilde D_{n-1}$;
\item two boundary components with $n_1$ and $n_2$ marked points and one puncture
  (a once-punctured annulus): type $\Gamma(n_1,n_2)$;
\item three boundary components with $n_1$, $n_2$ and $n_3$ marked points (a
    pair of pants): type $\Gamma(n_1,n_2,n_3)$. 
\end{itemize}

\begin{proof}[Proof of Proposition~\ref{pr:poly-growth-surface}]
Let a \emph{feature} of a surface be either a puncture or a boundary
component (with a positive number of marked points). 
  The graph~$\ETSM$ is covered by finitely many copies of the ordinary
  graph of ideal triangulations $\EESM$, which in turn is
  quasi-isometric to the mapping class group~$\MCG\SM$ (see,
  e.g.,~\cite{ivanov}). If $\SM$ is a
  sphere with at most $3$ features, then $\MCG\SM$ is generated
  by Dehn twists along the boundary components, which commute with
  each other; hence $\MCG\SM$ has polynomial growth. 
  Otherwise, we can find two simple closed curves~$C_1, C_2$ which
  intersect non-trivially on $\Surf$, as follows. If $\SM$ has more
  than 3~features, pick features $p_1$ through $p_4$ and pick a
  separating curve~$C_1$ which has $p_1$ and $p_2$ on one side and
  $p_3$ and $p_4$ on the other and a separating curve~$C_2$ which has
  $p_2$ and $p_3$ on one side and $p_4$ and $p_1$ on the other.
  If $\Surf$ has genus $g>0$,
  take $C_1$ and $C_2$ which have algebraic intersection number~1
  in~$H_1(\Surf)$.  In either case the Dehn twists along $C_1$ and $C_2$
  will   not commute.
  Therefore, by the work of Ivanov~\cite{ivanov-algebraic} and
  McCarthy~\cite{mccarthy} on the Tits alternative for mapping class
  groups, $\MCG\SM$ contains a free group on two generators and hence
  has exponential growth.
\end{proof}

It is straightforward to verify the precise growth rate.

\begin{remark}
Curiously, the mutation classes defined by the diagrams in
Figure~\ref{fig:poly-growth} 
(thus corresponding to a once-punctured annulus and a pair of pants,
respectively) 
can also be obtained by reorienting some of the edges in the ``minimal
$2$-infinite diagrams'' of types
$D_n^{(1)}(m,r)$ and $D_n^{(1)}(m,r,s)$, in A.~Seven's
nomenclature~\cite[Section~8]{seven-2-infinite}. 
\end{remark}

\section{Finite mutation classes}
\label{sec:finit-mut-classes}

It is an important problem of cluster combinatorics to describe all 
skew-sym\-metrizable $n\times n$ matrices~$B$ whose mutation equivalence
class is finite; see, e.g., \cite{buan-reiten} \cite[Section~4]{ca4}. 
This is obviously the case for $n\leq 2$, and also 
whenever $B$ defines a cluster algebra of finite type. 
But there are many more examples, as we will see, even if one
restricts (as we do in this paper) to the skew-symmetric case. 


A.~Buan and I.~Reiten used representation theory of quivers to obtain
the following nontrivial result.

\begin{theorem}[{\cite[Theorem~3.6]{buan-reiten}}]
\label{th:buan-reiten}
Let $\Gamma$ be a finite directed graph with 
$n\geq 3$ vertices and no oriented cycles.
Then the mutation equivalence class of the corresponding
$n\times n$ skew-symmetric matrix~$B=B(\Gamma)$
(see Definition~\ref{def:matrices-via-graphs}) 
is finite if and only if $\Gamma$ is an orientation of a simply-laced
Dynkin or extended Dynkin diagram. 
\end{theorem}

Our main theorems imply the following statement.

\begin{corollary}
\label{cor:mut-finite}
The mutation class of any matrix $B(T)$ associated to a triangulation
$T$ of a bordered surface with marked points is finite. 
\end{corollary}

Corollary~\ref{cor:mut-finite} is immediate from the fact that all matrices in
such a mutation class have entries~$\leq 2$. 
This is of course a very crude argument that hardly tells us anything
about the mutation class in question. 
The following simple observations shed some light on this issue. 

\begin{proposition}
\label{pr:mut-class-ideal}
The mutation class $\Bcal\SM$ consists of the
signed adjacency matrices of all ideal triangulations of~$\SM$
(up to simultaneous permutations of rows and columns).
\end{proposition}

Proposition~\ref{pr:mut-class-ideal} is immediate from 
Definition~\ref{def:B-tagged}. 
It can be used to give an alternative proof of the mutation-finiteness
property: Even though the number of ideal triangulations is typically 
infinite
(cf.\ Proposition~\ref{pr:finite-number-of-arcs}), 
the number of orbits of the mapping class group action is always
finite. Since the matrix $B(T)$ is invariant under this action, 
the claim follows. 

Yet another, completely combinatorial proof of mutation-finiteness 
can be given using block decompositions discussed in
Section~\ref{sec:block-decomp}. 

A skew-symmetric integer matrix~$B$ is \emph{acyclic} if the rows and
columns of~$B$ can be simultaneously permuted so that all entries
above the diagonal are nonnegative and all entries below the diagonal
are nonpositive.  Equivalently, $B = B(\Gamma)$ where $\Gamma$ is an
acyclic directed graph (as in Theorem~\ref{th:buan-reiten}).

\begin{corollary}
  The mutation class~$\Bcal\SM$ contains
  an acyclic matrix if and only if $\Bcal\SM$ is of finite or affine
  type.
\end{corollary}

\begin{proof}
  This follows from Corollary~\ref{cor:mut-finite} and
  Theorem~\ref{th:buan-reiten}.
\end{proof}

\begin{problem}[\emph{Recognizing finite mutation classes}] 
\label{problem:finite-mut-class}
Find an effective criterion/al\-gorithm 
for determining whether a given skew-symmetric 
(or, more generally, skew-symmetrizable) integer
matrix defines a finite mutation class. 
\end{problem}

Theorem~\ref{th:buan-reiten} and Corollary~\ref{cor:mut-finite}
provide two sources of finite mutation classes: the
simply-laced Dynkin or extended Dynkin diagrams, and the signed
adjacency matrices of triangulations. 
We next discuss three additional examples of this kind. 

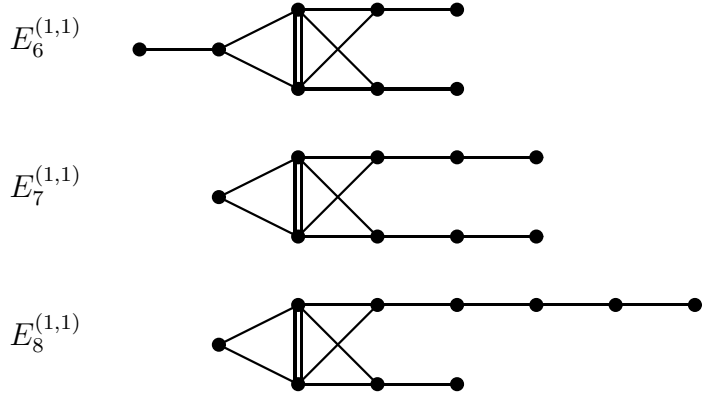
\begin{figure}[htbp] 
\vspace{-.2in} 
\[
\begin{array}{ll}
E_6^{(1,1)}\quad 
&
\setlength{\unitlength}{1.5pt} 
\begin{picture}(100,20)(0,10) 
\thicklines
\put(0,10){\circle*{3}} 
\put(20,10){\circle*{3}} 
\multiput(40,0)(20,0){3}{\circle*{3}} 
\multiput(40,20)(20,0){3}{\circle*{3}} 

\put(0,10){\line(1,0){20}} 
\put(20,10){\line(2,1){20}} 
\put(20,10){\line(2,-1){20}} 

\put(39.2,0){\line(0,1){20}} 
\put(40.8,0){\line(0,1){20}} 

\put(40,0){\line(1,1){20}} 
\put(40,20){\line(1,-1){20}} 

\put(40,0){\line(1,0){40}} 
\put(40,20){\line(1,0){40}} 
\end{picture} 
\\[.3in]
E_7^{(1,1)}
&
\setlength{\unitlength}{1.5pt} 
\begin{picture}(100,20)(0,10) 
\thicklines
\put(20,10){\circle*{3}} 
\multiput(40,0)(20,0){4}{\circle*{3}} 
\multiput(40,20)(20,0){4}{\circle*{3}} 

\put(20,10){\line(2,1){20}} 
\put(20,10){\line(2,-1){20}} 

\put(39.2,0){\line(0,1){20}} 
\put(40.8,0){\line(0,1){20}} 

\put(40,0){\line(1,1){20}} 
\put(40,20){\line(1,-1){20}} 

\put(40,0){\line(1,0){60}} 
\put(40,20){\line(1,0){60}} 

\end{picture} 
\\[.3in]
E_8^{(1,1)}
&
\setlength{\unitlength}{1.5pt} 
\begin{picture}(100,20)(0,10) 
\thicklines
\put(20,10){\circle*{3}} 
\multiput(40,0)(20,0){3}{\circle*{3}} 
\multiput(40,20)(20,0){6}{\circle*{3}} 

\put(20,10){\line(2,1){20}} 
\put(20,10){\line(2,-1){20}} 

\put(39.2,0){\line(0,1){20}} 
\put(40.8,0){\line(0,1){20}} 

\put(40,0){\line(1,1){20}} 
\put(40,20){\line(1,-1){20}} 

\put(40,0){\line(1,0){40}} 
\put(40,20){\line(1,0){100}} 

\end{picture} 
\end{array}
\]
\vspace{.2in} 
\caption{Extended affine diagrams $E_6^{(1,1)}$--$E_8^{(1,1)}$.
All triangles are oriented. 
} 
\label{fig:extended-affine} 
\end{figure}

\begin{definition}
\label{def:extended-affine}
A directed graph (quiver) $\Gamma$ (or the corresponding
skew-symmetric matrix $B=B(\Gamma)$) is of  
\emph{extended affine type} $E_6^{(1,1)}$,
$E_7^{(1,1)}$, or $E_8^{(1,1)}$, if $\Gamma$ is an orientation of the
corresponding graph in Figure~\ref{fig:extended-affine} in which all
triangles are oriented. 
More broadly, this nomenclature applies to the mutation equivalence
class of~$\Gamma$ or~$B$, or any quiver therein. 
\end{definition}

The quivers for the extended affine exceptional types 
shown in Figure~\ref{fig:extended-affine} are orientations 
of the Dynkin diagrams for extended affine root systems 
first described by Saito~\cite[Table~1]{saito}.  
The connection between cluster combinatorics and extended affine root
systems was first noticed by Geiss, Leclerc, and
Schr\"oer~\cite{gls}. 

Neither of these three mutation classes is associated with triangulated
surfaces. This can be shown using block decompositions of
Section~\ref{sec:block-decomp}. 

\begin{proposition}
\label{prop:extended-affine}
The mutation equivalence class of a quiver (or matrix) 
of extended affine type $E_6^{(1,1)}$--$E_8^{(1,1)}$
is finite. 
\end{proposition}

\begin{proof}
This statement was verified by exhaustive computer search 
with the help of \texttt{Java} applets for matrix
mutations written by Bernhard Keller~\cite{keller-applet}
and Lauren Williams. 
\end{proof}

\begin{remark}
Dynkin diagrams for
  extended affine root systems of types $A$ and $D$ do not generally 
define a
  finite mutation class (for any orientation). On the other hand, 
various subdiagrams of these diagrams 
  do come from surfaces.    
One nice example is a cluster algebra associated with a $4$-punctured sphere,
which has extended affine type~$D_4^{(1,1)}$.
\end{remark}

\begin{figure}[htbp] 
\vspace{-.2in} 
\begin{center}
\setlength{\unitlength}{1.5pt} 
\begin{picture}(60,47)(0,-2) 
\thicklines
\put(0,10){\circle*{3}} 
\put(0,-10){\circle*{3}} 
\put(20,0){\circle*{3}} 
\put(40,0){\circle*{3}} 

\put(20,40){\circle*{3}} 
\put(20,20){\circle*{3}} 
\put(40,30){\circle*{3}} 
\put(60,30){\circle*{3}} 

\put(20,0){\line(1,0){20}} 
\put(0,10){\line(2,-1){20}} 
\put(0,-10){\line(2,1){20}} 

\put(40,30){\line(1,0){20}} 
\put(20,40){\line(2,-1){20}} 
\put(20,20){\line(2,1){20}} 

\put(6,7){\vector(-2,1){2}} 
\put(6,-7){\vector(-2,-1){2}} 
\put(35,0){\vector(1,0){2}} 

\put(32,34){\vector(2,-1){2}} 
\put(32,26){\vector(2,1){2}} 
\put(45,30){\vector(-1,0){2}} 

\thinlines

\put(0,10){\line(2,3){20}} 
\put(0,-10){\line(2,3){20}} 
\put(20,0){\line(2,3){20}} 
\put(40,0){\line(2,3){20}} 

\put(12,28){\vector(2,3){2}} 
\put(12,8){\vector(2,3){2}} 
\put(50,15){\vector(2,3){2}} 
\put(32,18){\vector(-2,-3){2}}

\end{picture} 
\end{center}
\vspace{.1in} 
\caption{$A_2\times D_4$ diagram} 
\label{fig:A2D4} 
\end{figure}
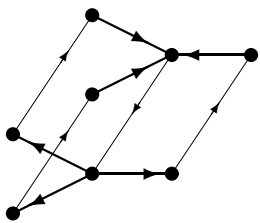

\begin{remark}
\label{rem:direct-product}
  The extended affine types $E_6^{(1,1)}$, $E_7^{(1,1)}$, and $E_8^{(1,1)}$ 
can be alternatively described using a construction of a ``direct
product'' of Dynkin quivers 
(discovered independently by F.~Chapoton~\cite{chapoton-leicester},
who also posed a version of
Problem~\ref{pr:finite-excep-finite-mutation-class} \cite{chapoton-luminy}).  
The definition of this construction should be clear from the example
of the quiver $A_2 \times D_4$
  shown in Figure~\ref{fig:A2D4}; this quiver is mutation equivalent to the
  extended affine quiver of type~$E_6^{(1,1)}$. 
Similarly, $A_3 \times A_3$ yields extended affine type~$E_7^{(1,1)}$
(the cluster type of $\operatorname{Gr}_{4,8}$), 
while $A_2 \times A_5$ gives extended affine type~$E_8^{(1,1)}$
(the cluster type of~$\operatorname{Gr}_{3,9}$).  
\end{remark}

\begin{problem}
\label{pr:finite-excep-finite-mutation-class}
Are there finitely many finite mutation classes of
indecomposable skew-sym\-metric $n\times n$ matrices ($n\geq2$) 
which are not associated with triangulated surfaces?
Are there any at all, besides the finite exceptional types $E_6$--$E_8$, 
the corresponding affine types $\widetilde E_6$--$\widetilde E_8$, 
and the extended affine types $E_6^{(1,1)}$--$E_8^{(1,1)}$?
\end{problem}

As an application, we classify the Grassmannians of
``finite mutation type.'' 
Recall that a cluster algebra structure on the
homogeneous coordinate ring $\CC[X]$ of a 
Grassmann manifold~$X=\operatorname{Gr}_{k,k+\ell}$ 
(with respect to its Pl\"ucker embedding) was described by
J.~Scott~\cite{scott} who showed, in the notation of
Remark~\ref{rem:direct-product}, that $\CC[X]$ has cluster type
$A_{k-1}\times A_{\ell-1}\,$. 

\begin{proposition}
\label{prop:G39-G48}
Let $\Acal$ be the canonical cluster algebra structure on the
homogeneous coordinate ring $\CC[X]$ of a 
Grassmannian~$\operatorname{Gr}_{k,k+\ell}$.  
The mutation equivalence class $\Bcal(\Acal)$ is finite 
if and only if $(k-2)(\ell-2)\leq 4$. 
\end{proposition}

A more concrete version of Proposition~\ref{prop:G39-G48} is presented in
Table~\ref{table:grassmann-finite}. 

\begin{table}[htbp]
\begin{tabular}{ccc}\toprule
Grassmannian~$X$ & & Cluster type of $\CC[X]$\\ \midrule
$\operatorname{Gr}_{2,n+3}\cong \operatorname{Gr}_{n+1,n+3}$ && $A_n$\\
$\operatorname{Gr}_{3,6}$ && $D_4$\\
$\operatorname{Gr}_{3,7}\cong \operatorname{Gr}_{4,7}$ && $E_6$\\
$\operatorname{Gr}_{3,8}\cong \operatorname{Gr}_{5,8}$ && $E_8$\\
$\operatorname{Gr}_{3,9}\cong \operatorname{Gr}_{6,9}$ && $E_8^{(1,1)}$\\
$\operatorname{Gr}_{4,8}$ && $E_7^{(1,1)}$\\
 \bottomrule
\end{tabular}
\caption{Grassmannians of finite mutation type}
\label{table:grassmann-finite}
\end{table}

\begin{proof}
The ``only if'' part is proved by checking that in the special cases
$\{k,\ell\}=\{3,7\}$ and $\{k,\ell\}=\{4,5\}$, the mutation class
$\Bcal(\Acal)$ is infinite. This is done by identifying, with the help
of a computer search, a sequence of mutations producing a
graph~$\Gamma'$ some of whose edges have multiplicity $\ge3$. 
(Any such graph with at least $3$ vertices is easily seen to generate
an infinite mutation class.) 

For the ``if'' part, the corresponding Grassmannians
$\operatorname{Gr}_{k,k+\ell}$ are listed in 
Table~\ref{table:grassmann-finite} on the left. 
It remains to show that in each of these cases,
the quiver $A_{k-1}\times A_{\ell-1}$ has cluster type given in the
right column. 
This is done by exhibiting an appropriate sequence of mutations
(which can be found by hand or with the help of a computer). 
\end{proof}

We note that the strict inequality
\hbox{$\!(k\!-\!2)(\ell\!-\!2)\!<\!4$} corresponds precisely to 
the situations where $\CC[\operatorname{Gr}_{k,k+\ell}]$ has 
finite cluster type.

\begin{problem} 
\label{problem:concrete-extended-affine}
Provide a concrete combinatorial description for cluster complexes of extended affine types
$E_6^{(1,1)}$--$E_8^{(1,1)}$. 
\end{problem}

\begin{remark}
It is not hard to show that the coordinate ring $\CC[\SL_4]$
(viewing~$\SL_4$ as a hypersurface 
$\det=1$ in $\CC^{16}$; 
cf.\ Example~\ref{example:cluster-type-A-affine}) 
has cluster type $E_7^{(1,1)}$, with respect to its usual cluster
algebra structure closely related to the one described in~\cite{ca3} for an
open double Bruhat cell. 
Similarly, the coordinate ring $\CC[SL_6/N]$ (here $N$ is the subgroup of unipotent
upper-triangular matrices in~$SL_6$)
has cluster type~$E_8^{(1,1)}$. 
Thus, a solution of Problem~\ref{problem:concrete-extended-affine} is likely to lead 
to a better understanding of the cluster combinatorics underlying 
the special linear group $\SL_4(\CC)$ and the affine base space
$SL_6(\CC)$, 
and the corresponding dual canonical bases. 
The cluster algebras~$\CC[\SL_n]$ ($n\geq 5$) and $\CC[\SL_n/N]$ ($n\geq 7$) 
are of infinite mutation type (``wild''), 
and therefore are unlikely to be understood in sufficient detail. 
\end{remark}

\section{Block decompositions}
\label{sec:block-decomp}

In view of Proposition~\ref{pr:mut-class-ideal}, the following classes
of skew-symmetric integer matrices coincide:
\begin{itemize}
\item
exchange matrices of cluster algebras associated with bordered
surfaces with marked points; 
\item
signed adjacency matrices of ideal triangulations; 
\item
signed adjacency matrices of tagged triangulations. 
\end{itemize}
In this section, we describe, in concrete combinatorial terms,
the set of matrices satisfying any of these equivalent descriptions. 

\begin{figure}[htbp] 
\begin{center}
\setlength{\unitlength}{1.5pt} 
\begin{picture}(20,33)(0,-12) 
\thicklines
\put(2,0){\line(1,0){16}} 
\put(7,0){\vector(1,0){6}} 
\multiput(0,0)(20,0){2}{\circle{4}}
\put(10,-10){\makebox(0,0){I}}
\end{picture} 
\qquad\quad
\begin{picture}(20,33)(0,-12) 
\thicklines
\put(2,0){\line(1,0){16}} 
\put(18,0){\vector(-1,0){11}} 
\put(1,1.5){\line(2,3){8}} 
\put(1,1.5){\vector(2,3){6}} 
\put(11,13.5){\vector(2,-3){6}} 
\put(19,1.5){\line(-2,3){8}} 
\multiput(0,0)(20,0){2}{\circle{4}}
\put(10,15){\circle{4}}
\put(10,-10){\makebox(0,0){II}}
\end{picture} 
\qquad\quad
\begin{picture}(20,33)(0,-12) 
\thicklines
\put(0,0){\line(2,3){9}} 
\put(20,0){\vector(-2,3){7}} 
\put(20,0){\line(-2,3){9}} 
\put(0,0){\vector(2,3){7}} 
\multiput(0,0)(20,0){2}{\circle*{4}}
\put(10,15){\circle{4}}
\put(10,-10){\makebox(0,0){IIIa}}
\end{picture} 
\qquad\quad
\begin{picture}(20,33)(0,-12) 
\thicklines
\put(0,0){\line(2,3){9}} 
\put(11,13.5){\vector(2,-3){5}} 
\put(20,0){\line(-2,3){9}} 
\put(9,13.5){\vector(-2,-3){5}} 
\multiput(0,0)(20,0){2}{\circle*{4}}
\put(10,15){\circle{4}}
\put(10,-10){\makebox(0,0){IIIb}}
\end{picture} 
\qquad\quad
\begin{picture}(20,48)(0,-27) 
\thicklines
\multiput(1,1.5)(10,-15){2}{\line(2,3){8}} 
\put(1,1.5){\vector(2,3){5}} 
\put(11,13.5){\vector(2,-3){6}} 
\multiput(1,-1.5)(10,15){2}{\line(2,-3){8}} 
\put(1,-1.5){\vector(2,-3){5}} 
\put(11,-13.5){\vector(2,3){6}} 
\multiput(10,15)(0,-30){2}{\circle*{4}}
\multiput(0,0)(20,0){2}{\circle{4}}
\put(2,0){\line(1,0){16}} 
\put(18,0){\vector(-1,0){11}} 
\put(10,-25){\makebox(0,0){IV}}
\end{picture} 
\qquad\quad
\begin{picture}(30,33)(0,-12) 
\thicklines
\multiput(0,0)(0,30){2}{\line(1,0){30}} 
\multiput(0,0)(30,0){2}{\line(0,1){30}} 
\multiput(0,0)(0,30){2}{\circle*{4}}
\multiput(30,0)(0,30){2}{\circle*{4}}
\put(15,15){\circle{4}}
\put(0,0){\vector(1,0){20}} 
\put(30,30){\vector(-1,0){20}} 
\put(30,30){\vector(0,-1){20}} 
\put(0,0){\vector(0,1){20}} 

\put(0,0){\line(1,1){13.5}}
\put(30,30){\line(-1,-1){13.5}}
\put(0,30){\line(1,-1){13.5}}
\put(30,0){\line(-1,1){13.5}}
\put(16.5,16.5){\vector(1,1){7.5}}
\put(13.5,13.5){\vector(-1,-1){7.5}}
\put(0,30){\vector(1,-1){10}}
\put(30,0){\vector(-1,1){10}}

\put(15,-10){\makebox(0,0){V}}
\end{picture} 
\end{center}
\vspace{-.1in} 
\caption{Blocks of types I--V} 
\label{fig:blocks-I-V} 
\end{figure}
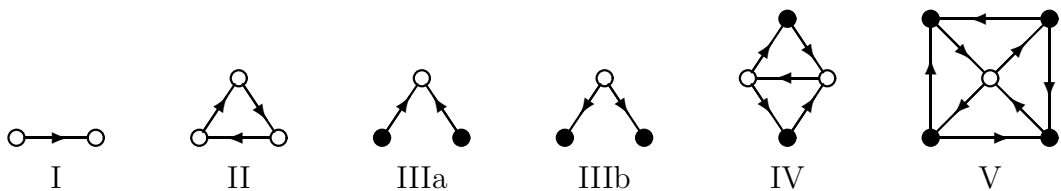

\begin{definition}
\label{def:block-decomposable}
A \emph{block} is a directed graph isomorphic to one of the graphs
shown in Figure~\ref{fig:blocks-I-V}. 
Depending on which graph it is, we call it a block of \emph{type} I,
II, III, IV, or~V. 
The vertices marked by unfilled circles in Figure~\ref{fig:blocks-I-V}
are called \emph{outlets}.
A directed graph $\Gamma$ is called \emph{block-decomposable} if it can
be obtained from a collection of disjoint blocks by the following
procedure. Take a partial matching of the combined set of outlets; 
matching an outlet to itself or to another outlet from the same block is
not allowed. Identify (or ``glue'') the vertices within
each pair of the matching. We require that the resulting graph~$\Gamma'$ be
connected. 
If $\Gamma'$ contains a pair of edges connecting the same pair of 
vertices but going in opposite directions, then remove each such a
pair of edges.
The result is a block-decomposable graph~$\Gamma$. 
\end{definition}

By design, a block-decomposable graph has no loops, and all edge
multiplicities are $1$ or~$2$. 

\begin{example}
\label{example:gluing-blocks}
Here are some examples illustrating
Definition~\ref{def:block-decomposable}. 
Two blocks of type~I can be glued to create any of the
following block-decomposable graphs: 
\begin{itemize}
\item
two isolated vertices; 
\item
a two-vertex graph with two edges of the same orientation;  
\item
an arbitrary orientation of a type~$A_3$ Coxeter-Dynkin diagram. 
\end{itemize}
\end{example}

\begin{theorem}
\label{th:B-criterion}
An integer matrix $B$ is a signed adjacency matrix of an ideal
triangulation of a bordered surface with marked points if and only if
$B$ is encoded by a block-decomposable graph~$\Gamma$
(i.e., $B=B(\Gamma)$ in the sense of
Definition~\ref{def:matrices-via-graphs}). 
\end{theorem}

\begin{proof}
Let $B=B(T)$ be the signed adjacency matrix of an ideal
triangulation~$T$. 
As explained in Remark~\ref{rem:puzzle-pieces}, with one exception the
matrix $B(T)$
can be computed by taking a puzzle-piece decomposition of~$T$
and adding up the contributions of the puzzle pieces that make up~$T$. 
It is routine to verify that this leads to a block decomposition of
the graph~$\Gamma$ associated with~$B$, as follows: 
\begin{itemize}
\item
the puzzle pieces shown in Figure~\ref{fig:puzzle-pieces} 
whose outer sides do not lie on the boundary of~$T$
(that is, are matched to outer sides of other puzzle pieces)
contribute blocks of types~II, IV, and~V, respectively; 
\item
a triangle with two sides on the boundary does not contribute 
to the block decomposition; 
\item
a triangle with one side on the boundary contributes a block of
type~I; 
\item
a punctured-digon puzzle piece with both sides on the boundary yields two disjoint
vertices (see Example~\ref{example:gluing-blocks}); 
\item
a punctured-digon puzzle piece with one side on the boundary 
contributes a block of type~III; 
\item
a twice-punctured monogon puzzle piece with the side on the boundary yields a graph
of type $\widetilde A(2,2)$ which can be obtained by gluing four
blocks of type~I.
\end{itemize}
We also need to consider the case of the exceptional triangulation of a
$4$-punctured sphere that cannot be decomposed into puzzle pieces. 
The corresponding $6\times 6$ matrix $B(T)$ is encoded by a
$6$-vertex graph~$\Gamma$ which is an orientation of
the $1$-skeleton of an octahedron.
This graph $\Gamma$ can be obtained by gluing four blocks of type~II. 

It remains to prove that for each block-decomposable graph~$\Gamma$, 
we have $B(\Gamma)=B(T)$ for some tagged triangulation~$T$. 
This is done by reversing the argument in the previous paragraph: 
each block in Figure~\ref{fig:blocks-I-V} 
corresponds to (a piece of) a triangulated bordered surface 
in which the signed adjacencies are described by the arrows in the
block, and whose exposed sides are precisely the outlets.  Glue the
pieces of the triangulated surface according to the matching of the
outlets.  Glue an additional triangle to the side corresponding to any
unmatched outlet, leaving the other two sides of this new triangle
exposed.
\end{proof}

\begin{remark}
A block decomposition of a given graph~$\Gamma$ may not be unique. 
For example, the graph $\Gamma$~of type $\widetilde A(2,2)$
(see Figure~\ref{fig:tilde-a22-a13})
can be obtained by gluing either four blocks of type~I,
or a block of type~I and a block of type~IV. 
\end{remark}

\begin{remark}
\label{rem:block-closed-mutation}
It can be proved by a direct case-by-case combinatorial argument that
the class of block-decomposable graphs is closed under mutations.
(This leads to an alternative proof of mutation finiteness.) 
It would be interesting to find a way to enlarge the set of blocks 
while keeping this invariance property. 
This would yield a better understanding (or even a complete
description) of the class of diagrams whose mutation equivalence class
is finite; cf.\ Problem~\ref{problem:finite-mut-class}. 

It can be shown by tedious case-by-case analysis that a graph~$\Gamma$
obtained by arbitrarily orienting the edges of a Coxeter-Dynkin
diagram of type~$E_n$ ($n\in\{6,7,8\}$),
or an affine extension thereof, is \emph{not} block
decomposable. 
On the other hand, all these graphs have finite mutation classes
(either provably or conjecturally); it would be nice to find a
conceptual explanation of these phenomena. 
\end{remark}

\section{Recovering topology from combinatorial data}
\label{sec:recover-topology}

Let $B$ be a skew-symmetric matrix satisfying the criterion in
Theorem~\ref{th:B-criterion}, so that there exists a (non-unique)
bordered surface $\SM$ and a (non-unique) ideal (equivalently, tagged)
triangulation~$T$ of~$\SM$ such that $B=B(T)$. 
How much information about $\SM$ can we recover from~$B$? 

Examples~\ref{example:monogon2} and~\ref{example:triangle-hexagon}
show that the topology of~$\Surf$ is not in general determined
by~$B(T)$, and even if it is, then the number of punctures is not
uniquely determined. 

Furthermore, there exist topologically inequivalent triangulations
$T_1$ and $T_2$ of the same marked surface $\SM$ for which
$B(T_1)=B(T_2)$.
Example: let $\SM$ be a quadrilateral $abcd$ with a puncture~$t$
(type~$D_4$). 
Triangulation $T_1$ is obtained by connecting $t$ to $a,b,c,d$,  
while $T_2$ is obtained by connecting $t$ to $a,c$, 
and connecting $a$ with $c$ by two arcs. 
Additional examples can be found by considering a twice-punctured digon 
(type $\widetilde D_4$) or a sphere with $4$ punctures (mentioned in
the proof of Theorem~\ref{th:B-criterion}).

On the positive side, we have the following result.

\begin{proposition}
\label{pr:block-unique}
  Let~$T$ be an ideal triangulation of $\SM$ such that $B(T) =
  B(\Gamma)$ where the graph~$\Gamma$ has a unique block
  decomposition.  Then the topology of $\SM$ is uniquely determined by
  $B(T)$.  More precisely, for any $\SMprime \not\cong \SM$ and any
  ideal triangulation~$T'$ of $\SMprime$, the matrix~$B(T')$ is
  mutation inequivalent to $B(T)$.
\end{proposition}

\begin{proof}
  It is enough to show that $B(T') \ne B(T)$.  Suppose not. By
  Remark~\ref{rem:puzzle-pieces}, every triangulation~$T$ has a
  puzzle-piece decomposition, except for one triangulation of a sphere
  with $4$ punctures.  By the proof of Theorem~\ref{th:B-criterion},
  every puzzle-piece decomposition gives a block decomposition, except
  for one triangulation of a twice-punctured monogon.  In both
  exceptional cases the graph~$\Gamma$ associated with~$T$ has more
  than one block decomposition.  As a result, if the hypotheses of the
  theorem are satisfied, the triangulations~$T$ and $T'$ can be
  decomposed into puzzle pieces leading to block decompositions of
  $\Gamma$, which must coincide by hypothesis.  By the argument at the
  end of the proof of Theorem~\ref{th:B-criterion}, the puzzle pieces
  of~$T$ (or $T'$) and the way they are glued to each other are
  uniquely determined by the block decomposition of~$\Gamma$.  Thus
  $T = T'$ as desired.
\end{proof}

\begin{proof}[Proof of Lemma~\ref{lem:affine-a}] ($\widetilde
  A(n_1,n_2) \ne \widetilde A(n_1', n_2')$)
For $\{n_1, n_2\} = \{1,1\}$ or $\{1,2\}$ there is nothing to prove.
Otherwise we may suppose $n_1 \ge 3$ (possibly switching $\{n_1,n_2\}$
and $\{n_1', n_2'\}$ and/or $n_1$ and $n_2$).
Let~$\Gamma$ be the graph of type $\widetilde A(n_1,n_2)$ obtained by
orienting a cycle of length $n_1 + n_2$ so that $n_1$ consecutive
edges are oriented the same way.
A straightforward inspection shows that~$\Gamma$ has a unique block
decomposition, and the lemma follows by Proposition~\ref{pr:block-unique}.
\end{proof} 

Suppose that $\SM\mapsto f(\SM)$ is a surface invariant 
that can be determined from the signed adjacency matrix $B(T)$ of an ideal
triangulation of~$\SM$;
that is, there exists a function $B\mapsto F(B)$ such that 
$f(\SM)=F(B(T))$. 
Then $F(B)$ must be invariant under matrix mutations. 
Consequently, the problem at hand is closely related to the following
fundamental problem in cluster combinatorics. 

\begin{problem}[\emph{Mutation invariants}] 
Describe functions on skew-symmetric (or skew-symmetrizable) matrices
that are invariant under matrix mutations. 
Use such invariants to develop efficient algorithms that distinguish
between mutation-in\-equi\-valent matrices. 
\end{problem}

Ideally, one would like to have a complete system of invariants for
matrices of given order, so that if $B$ and $B'$ are not mutation
equivalent, then one of these invariants takes different values at $B$
and~$B'$. 

Very few mutation invariants have been identified so far.  The best
known example is the (ordinary) \emph{rank} of a matrix; see
\cite[Lemma~1.2]{gsv1} \cite[Lemma~3.2]{ca3}.  
Note that for $B$ a skew-symmetric matrix, $\rank(B)$ is
even~\cite[Section 21]{prasolov}.
The following theorem is stated in terms of the \emph{corank}, the
dimension~$n$ minus the rank.

\begin{theorem}
\label{th:corank}
For an arbitrary ideal triangulation~$T$ of a marked surface~$\SM$, the
corank of $B(T)$ equals the number of punctures plus the number of
boundary components having an even number of marked points.
\end{theorem}

The statement of Theorem~\ref{th:corank} becomes more elegant if one
follows Fock and Goncharov~\cite{fg-dual-teich} and treats a puncture
as a boundary component with no marked points; then the corank is just
the number of boundary components with an even number of marked
points (possibly none).

\begin{proof}
We argue by joint induction on the size~$n$ of the
matrix and the number of boundary components.
Since $\rank(B)$ is invariant under mutations (flips),
we are free to choose the triangulation~$T$ of~$\SM$ as we please.

The following lemma lets us trade one puncture for two marked points
on the boundary without changing $\rank(B)$.

\begin{lemma}
  \label{lem:trade-puncture-marked-pts}
  Let $\SM$ be a surface with a puncture~$b$ and another marked
  point~$a$.  Form a surface $\SMprime$ by taking $\SM$ and deleting
  the puncture~$b$ and either:
  \begin{itemize}
  \item if $a$ is on the boundary of~$\Surf$, add two marked points to
    the boundary component containing~$a$; or
  \item if $a$ is a puncture, replace $a$ by a boundary component with
    two marked points.
  \end{itemize}
  Then for any triangulations~$T$ of $\SM$ and $T'$ of $\SMprime$, the
  rank of $B(T)$ equals the rank of $B(T')$.
\end{lemma}

\begin{proof}
  Construct a triangulation~$T$ of $\SM$ that includes an
  arc~$\beta$ connecting $a$ to~$b$ and a loop $\gamma$ encircling $\beta$
  and based at~$a$, making a self-folded triangle (see
  Figure~\ref{fig:cutting-open}).  By
  Definition~\ref{def:signed-adj-matrix}, the
  two rows (resp.~columns) corresponding to~$\beta$ and~$\gamma$ in~$B(T)$ are
  identical.  Thus the rank is unchanged (and the corank goes down
  by~$1$) if we delete the row and column corresponding to~$\beta$ from 
  $B(T)$.  The new matrix is the same as $B(T')$ for a triangulation
  $T'$ of $\SMprime$ obtained by cutting~$T$ open along~$\beta$.
\end{proof}

\begin{figure}[htbp]
\begin{center}
\setlength{\unitlength}{2pt}
\begin{picture}(40,32)(0,-2)
\thicklines


\qbezier(20,0)(2,32)(20,32)
\qbezier(20,0)(38,32)(20,32)
\put(20,0){\line(0,1){20}}
\put(0,0){\line(1,0){40}}

\multiput(20,0)(0,20){2}{\circle*{2}}
\put(8,22){\makebox(0,0){$\gamma$}}
\put(17.5,15){\makebox(0,0){$\beta$}}
\put(25,3){\makebox(0,0){$a$}}
\put(20,24){\makebox(0,0){$b$}}

\multiput(0,-2)(2,0){21}{\circle*{.5}}
\multiput(0,-4)(2,0){21}{\circle*{.5}}
\multiput(0,-6)(2,0){21}{\circle*{.5}}
\end{picture}
\begin{picture}(10,32)(0,-2)
\put(5,15){\makebox(0,0){$\leadsto$}}
\end{picture}
\begin{picture}(40,32)(0,-2)
\thicklines


\qbezier(15,0)(2,32)(20,32)
\qbezier(25,0)(38,32)(20,32)
\put(15,0){\line(1,4){5}}
\put(25,0){\line(-1,4){5}}
\put(0,0){\line(1,0){15}}
\put(25,0){\line(1,0){15}}

\put(20,20){\circle*{2}}
\multiput(15,0)(10,0){2}{\circle*{2}}
\put(6,22){\makebox(0,0){$\gamma$}}

\multiput(0,-2)(2,0){21}{\circle*{.5}}
\multiput(0,-4)(2,0){21}{\circle*{.5}}
\multiput(0,-6)(2,0){21}{\circle*{.5}}
\multiput(18,0)(2,0){3}{\circle*{.5}}
\multiput(18,2)(2,0){3}{\circle*{.5}}
\multiput(18,4)(2,0){3}{\circle*{.5}}
\multiput(18,6)(2,0){3}{\circle*{.5}}
\multiput(18,8)(2,0){3}{\circle*{.5}}
\multiput(18,10)(2,0){3}{\circle*{.5}}
\multiput(20,12)(0,2){4}{\circle*{.5}}

\end{picture}
\qquad\qquad
\begin{picture}(30,32)(5,-2)
\thicklines


\qbezier(20,0)(2,32)(20,32)
\qbezier(20,0)(38,32)(20,32)
\put(20,0){\line(0,1){20}}

\multiput(20,0)(0,20){2}{\circle*{2}}
\put(8,22){\makebox(0,0){$\gamma$}}
\put(17.5,15){\makebox(0,0){$\beta$}}
\put(24,0){\makebox(0,0){$a$}}
\put(20,24){\makebox(0,0){$b$}}

\end{picture}
\begin{picture}(10,32)(0,-2)
\put(5,15){\makebox(0,0){$\leadsto$}}
\end{picture}
\begin{picture}(30,32)(5,-2)
\thicklines


\qbezier(20,0)(2,32)(20,32)
\qbezier(20,0)(38,32)(20,32)

\qbezier(20,0)(13,15)(20,20)
\qbezier(20,0)(27,15)(20,20)
\multiput(20,1)(0,2){10}{\circle*{.5}}
\multiput(18,9)(0,2){4}{\circle*{.5}}
\multiput(22,9)(0,2){4}{\circle*{.5}}

\multiput(20,0)(0,20){2}{\circle*{2}}
\put(8,22){\makebox(0,0){$\gamma$}}

\end{picture}

\end{center}
\caption{Making a cut inside a self-folded triangle}
\label{fig:cutting-open}
\end{figure}
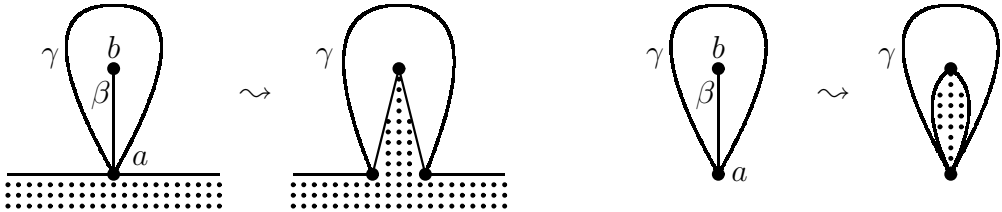

  If there is a puncture and at least one other marked point, we
  may remove all the punctures by repeated applications of
  Lemma~\ref{lem:trade-puncture-marked-pts}.  At each application, we
  leave $\rank(B)$ unchanged and decrease~$n$ by one, thus decreasing the
  corank by~1, as desired. Thus, we may henceforth assume that $\SM$
  has no punctures, or else is closed with a single puncture. 

  Suppose that there is more than one boundary component.
  If a boundary component has at least $3$ marked points, we can
  reverse the transformation in Lemma~\ref{lem:trade-puncture-marked-pts}
  to temporarily trade
  two of them for a puncture in the interior (increasing $n$ by 1 and
  not changing $\rank(B)$), and then trade the puncture for two
  marked points on another boundary component, again without changing
  $\rank(B)$.  Thus we can move pairs of marked points from one boundary
  component to another without changing either $n$ or the rank, or the
  number of even boundary components.

  If a boundary component has $2$ marked points, 
  we can again reverse Lemma~\ref{lem:trade-puncture-marked-pts} 
  to trade it for $2$ punctures, and then trade them back for
  $4$ marked points on another boundary component. As a result, we decrease $n$ by~1,
  leave $\rank(B)$ unchanged, and remove one component with an
  even number of punctures. 

  So we only need to treat the cases where at most one boundary
  component has more than one marked point. 
  If there are at least two boundary components, then take a component
  with $1$ marked point, and construct a triangulation~$T$ that
  includes the arcs  shown in
  Figure~\ref{fig:remove-component1} on the left.

\begin{figure}[htbp]
\begin{center}
\setlength{\unitlength}{2pt}
\begin{picture}(40,40)(0,-2)
\thicklines


\qbezier(20,0)(-15,42)(20,42)
\qbezier(20,0)(55,42)(20,42)

\qbezier(20,0)(2,37)(20,32)
\qbezier(20,0)(38,37)(20,32)
\put(0,0){\line(1,0){40}}
\put(20,26){\circle{12}}

\multiput(20,0)(0,32){2}{\circle*{2}}
\put(0,35){\makebox(0,0){$3$}}
\put(9,30){\makebox(0,0){$1$}}
\put(32,30){\makebox(0,0){$2$}}

\multiput(15,26)(2,0){6}{\circle*{.5}}

\multiput(15,28)(2,0){6}{\circle*{.5}}
\multiput(15,24)(2,0){6}{\circle*{.5}}

\multiput(17,30)(2,0){4}{\circle*{.5}}
\multiput(17,22)(2,0){4}{\circle*{.5}}

\multiput(0,-2)(2,0){21}{\circle*{.5}}
\multiput(0,-4)(2,0){21}{\circle*{.5}}
\multiput(0,-6)(2,0){21}{\circle*{.5}}
\end{picture}
\begin{picture}(10,38)(0,-2)
\put(5,15){\makebox(0,0){$\leadsto$}}
\end{picture}
\begin{picture}(40,38)(0,-2)
\thicklines


\qbezier(15,0)(2,32)(20,32)
\qbezier(25,0)(38,32)(20,32)
\put(15,0){\line(1,4){5}}
\put(25,0){\line(-1,4){5}}
\put(0,0){\line(1,0){15}}
\put(25,0){\line(1,0){15}}

\put(20,20){\circle*{2}}
\multiput(15,0)(10,0){2}{\circle*{2}}

\multiput(0,-2)(2,0){21}{\circle*{.5}}
\multiput(0,-4)(2,0){21}{\circle*{.5}}
\multiput(0,-6)(2,0){21}{\circle*{.5}}
\multiput(18,0)(2,0){3}{\circle*{.5}}
\multiput(18,2)(2,0){3}{\circle*{.5}}
\multiput(18,4)(2,0){3}{\circle*{.5}}
\multiput(18,6)(2,0){3}{\circle*{.5}}
\multiput(18,8)(2,0){3}{\circle*{.5}}
\multiput(18,10)(2,0){3}{\circle*{.5}}
\multiput(20,12)(0,2){4}{\circle*{.5}}

\end{picture}

\end{center}
\caption{Removing a connected component with one marked point}
\label{fig:remove-component1}
\end{figure}
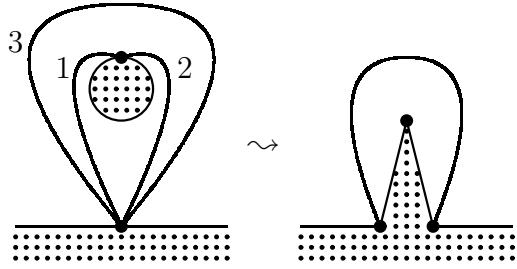

  The signed adjacency matrix~$B(T)$ has the form
  \[
  \begin{bmatrix}
    0 & -2 &  1 & 0 & 0 & \cdots \\
    2 &  0 & -1 & 0 & 0 & \cdots \\
   -1 &  1 &  0 & * & * & \cdots \\
    0 &  0 &  * & * & * & \cdots \\
    0 &  0 &  * & * & * & \cdots\\
    \vdots&\vdots&\vdots&\vdots&\vdots&\ddots
  \end{bmatrix}\,. 
  \]
  Adding the half-sum of columns~$1$ and~$2$ to column~$3$, and doing the
  same for the rows, we obtain the matrix
  \[
  \begin{bmatrix}
    0 & -2 &  0 & 0 & 0 & \cdots \\
    2 &  0 &  0 & 0 & 0 & \cdots \\
    0 &  0 &  0 & * & * & \cdots \\
    0 &  0 &  * & * & * & \cdots \\
    0 &  0 &  * & * & * & \cdots\\
    \vdots&\vdots&\vdots&\vdots&\vdots&\ddots
  \end{bmatrix}\,. 
  \]
  This matrix is a block sum of two matrices.  One is of dimension~$2$
  and rank~$2$, while the other is~$B(T')$, where $T'$ is the
  triangulation shown on the right of
  Figure~\ref{fig:remove-component1}.  Thus we can remove this (odd) 
  boundary component, decrease $n$ by~$2$,
  and add two points to another boundary component
  without changing the corank, as desired.

We are left with the cases where $\SM$ is a surface of genus~$g \ge 0$
with either one boundary component or one puncture.  If $g > 0$ and
$\SM$ is not the once-punctured torus, we can encircle one handle of
$\Surf$ with an arc based at a puncture and triangulate the enclosed
torus with $4$ additional arcs, forming a triangulation~$T$.  Again by
row and column operations we can reduce the matrix to a block-diagonal
form, with one block of dimension~$4$ and rank~$4$ and the other 
block of the form~$B(T')$, where $T'$ is a triangulation of a surface of genus
$g-1$ and two more marked points on the boundary.  In particular, the
corank is unchanged.

In this way, we can reduce everything to the base cases: a once-punctured
torus and an unpunctured $c$-gon ($c > 3$).  
In these cases, the theorem is verified by direct calculation: 
for the once-punctured torus, $\operatorname{corank}(B)=1$,
while for the $c$-gon, $\operatorname{corank}(B)=0$ if $c$ is
odd and $\operatorname{corank}(B)=1$ if $c$ is even.
\end{proof}

\begin{corollary}
\label{cor:genus-from-rank}
If $\Surf$ is closed (so that $\Mark$ consists exclusively of
punctures), then its genus~$g$ and the number of
punctures~$p$ are determined by the size $n$ and the rank $r$ of
the matrix~$B(T)$ associated with any ideal triangulation~$T$
of~$\SM$: 
\[
g=\frac{3r-2n+6}{6}\,,\qquad
p=n-r\,. 
\]
\end{corollary}

\begin{proof}
Immediate from Theorem~\ref{th:corank}. 
\end{proof} 

\section*{Acknowledgments}

The authors thank 
Vladimir Fock, Alexander Goncharov, Nikolai Ivanov, Robert Penner,
David Speyer, Lauren Williams,
and Andrei Zelevinsky for helpful advice and stimulating discussions.
We thank Bernhard Keller and Lauren Williams for providing software
that computes matrix mutations,
and Nathan Reading for letting us borrow 
Figures~\ref{fig:A3assoc_dual} and~\ref{fig:A2assoc_basic}
from~\cite{fomin-reading-pcmi}.

We are grateful to the anonymous referee and to Daniel Labardini Fragoso 
for their thorough readings of the original submitted version,
and for a number of valuable editorial suggestions whose
implementation improved the quality of the paper. 

S.~F.\ acknowledges the hospitality of the Mittag-Leffler Institute in
April-June, 2005; some of our main results were first presented there.
Much of the work was completed while D.~T.\ was a Benjamin Peirce
Assistant Professor at Harvard University.

\end{document}